\documentclass[12pt,a4paper]{article}
\usepackage[utf8]{inputenc}
\usepackage[english]{babel}
\usepackage[OT1]{fontenc}
\usepackage{amsmath}
\usepackage{amsfonts}
\usepackage{amssymb}
\usepackage[dvipsnames]{xcolor}
\usepackage{hyperref}
\usepackage{cleveref}
\usepackage{tikz}
\usepackage{pgfplots}
\usepackage{float}
\floatstyle{boxed}
\restylefloat{figure}

\def\R{\mathbb R}

\def\1{\mathbf 1}

\def\1{\bold 1}

\usepackage{amsthm}
\theoremstyle{theorem}
\newtheorem{theorem}{Theorem}[section]
\newtheorem{proposition}[theorem]{Proposition}
\newtheorem{lemma}[theorem]{Lemma}
\newtheorem{condition}[theorem]{Condition}
\newtheorem{remark}[theorem]{Remark}
\newtheorem{corollary}[theorem]{Corollary}
\theoremstyle{plain}
\newtoks\thehProclaim
\newtheorem*{Proclaim}{\the\thehProclaim}

\newtheorem{example}[theorem]{Example}

\theoremstyle{definition}
\newtoks\thehDefinition
\newtheorem*{Definition}{\the\thehDefinition}

\numberwithin{equation}{section}

\providecommand{\keywords}[1]{\textbf{{Key words:}} #1}
\providecommand{\subjclass}[1]{\textit{{2020 Mathematics Subject Classification.}} #1}
\begin{document}
\title{Variations on the theme of the Trotter-Kato theorem for homogenization of periodic hyperbolic systems
\footnote{\subjclass{Primary 35B27. Secondary 35L52.}\break Research was supported by «Native towns», a social investment program of PJSC «Gazprom Neft», and  
by the Swedish Research Council under grant no. 2016-06596 while the author was in residence at Institut Mittag-Leffler in Djursholm, Sweden during the research program ,,Spectral Methods in Mathematical Physics.'' A final touch to the paper (such as changing of notation)
 was  given under support of the European Research Council (ERC) under the European Union’s Horizon 2020 research and innovation programme (grant agreement No.~818437).}
}
\author{Yulia Meshkova \footnote{Chebyshev Laboratory, St. Petersburg State University, 14th Line V.O., 29B, Saint Petersburg 199178 Russia; Department of Mathematics and Statistics, P.O. Box 68 (Gustaf H\"allstr\"omin katu 2), FI-00014 University of Helsinki, Finland.  E-mail: {\tt{y.meshkova@spbu.ru}, \tt{iuliia.meshkova@helsinki.fi}}. }}

\maketitle

\begin{abstract}
In $L_2(\mathbb{R}^d;\mathbb{C}^n)$, we consider a matrix elliptic second order differential operator $B_\varepsilon >0$. Coefficients of the operator $B_\varepsilon$ are periodic with respect to some lattice in $\mathbb{R}^d$ and depend on $\mathbf{x}/\varepsilon$. We study the quantitative homogenization for the solutions of the hyperbolic system $\partial _t^2\mathbf{u}_\varepsilon =-B_\varepsilon\mathbf{u}_\varepsilon$. In operator terms, we are interested in approximations of the operators $\cos (tB_\varepsilon ^{1/2})$ and $B_\varepsilon ^{-1/2}\sin (tB_\varepsilon ^{1/2})$ in suitable operator norms. Approximations for the resolvent $B_\varepsilon ^{-1}$ have been already obtained by T.~A.~Suslina. So, we rewrite hyperbolic equation as a system for the vector with components $\mathbf{u}_\varepsilon $ and $\partial _t\mathbf{u}_\varepsilon$, and consider the corresponding unitary group. For this group, we adapt the proof of the Trotter-Kato theorem by introduction of some correction term and derive hyperbolic results from elliptic ones.
\end{abstract}


\section*{Introduction}

The paper is devoted to homogenization of periodic differential operators (DO's). More precisely, we are interested in the so-called operator error estimates, i.~e., in quantitative homogenization results, admitting formulation in the uniform operator topology.

For elliptic and parabolic problems, estimates of such type are very well studied, see, e.~g., books \cite[Chapter 14]{CDGr} and \cite{Sh}, the survey \cite{ZhPasUMN}, the papers \cite{BSu,SuAA} and references therein. For hyperbolic problems the situation is different. In the present paper, we deal with the hyperbolic systems in $\mathbb{R}^d$. 

\subsection{Problem setting} Let $\Gamma\subset\mathbb{R}^d$ be a lattice and let $\Omega$ be the elementary cell of the lattice $\Gamma$. For a $\Gamma$-periodic function $\psi$ in $\mathbb{R}^d$, we denote $\psi ^\varepsilon (\mathbf{x}):=\psi (\mathbf{x}/\varepsilon)$, where $\varepsilon >0$, and $\overline{\psi}:=\vert \Omega\vert ^{-1}\int_\Omega \psi (\mathbf{x})\,d\mathbf{x}$.

In $L_2(\mathbb{R}^d;\mathbb{C}^n)$, we study a selfadjoint matrix strongly elliptic second order DO ${B}_{\varepsilon}$, $0<\varepsilon\leqslant 1$.
The principal part of the operator ${B}_{\varepsilon}$ is given in a factorized form
$A_{\varepsilon}=b(\mathbf{D})^*g^\varepsilon (\mathbf{x})b(\mathbf{D})$,
where $b(\mathbf{D})$ is a matrix homogeneous first order DO, and $g(\mathbf{x})$ is a $\Gamma$-periodic bounded and positive definite matrix-valued function in  $\mathbb{R}^d$. (The precise assumptions on $b(\mathbf{D})$ and $g(\mathbf{x})$ are given below in Subsection~\ref{Subsection Operator A_eps}.) The simplest example of the operator $A_\varepsilon$ is the scalar elliptic operator $A_\varepsilon =-\mathrm{div}\,g^\varepsilon (\mathbf{x})\nabla$ (acoustics operator).  The matrix-valued example is the operator of elasticity theory. It can be written as $b(\mathbf{D})^*g^\varepsilon (\mathbf{x})b(\mathbf{D})$, see \cite[Chapter 2, \S 5]{BSu}.

The operator ${B}_{\varepsilon}$ is given by the differential expression
\begin{equation}
\label{B_eps in introduction}
{B}_{\varepsilon}=b(\mathbf{D})^* g^\varepsilon (\mathbf{x}) b(\mathbf{D})
+\sum_{j=1}^d\bigl(a_j^\varepsilon (\mathbf{x})D_j+D_ja_j^\varepsilon(\mathbf{x})^*\bigr)
+Q^\varepsilon (\mathbf{x}) +\lambda Q_0^\varepsilon (\mathbf{x}).
\end{equation}
Here $a_j(\mathbf{x})$, $j=1,\dots,d$, and $Q(\mathbf{x})$ are $\Gamma$-periodic matrix-valued functions, in general, unbounded;
a $\Gamma$-periodic matrix-valued function $Q_0(\mathbf{x})$ is such that $Q_0(\mathbf{x})>0$ and $Q_0, Q_0^{-1}\in L_\infty$.
The constant $\lambda$ is chosen so that the operator $B_{\varepsilon}$ is positive definite.
(The precise assumptions on the coefficients are given below in Subsection~\ref{Subsection lower order terms}.) An example of operator \eqref{B_eps in introduction} is the  Schrödinger operator with the singular potential of the form $\varepsilon ^{-1} v_\varepsilon$, see \cite[\S 11]{SuAA}.

The coefficients of the operator \eqref{B_eps in introduction} oscillate rapidly for small $\varepsilon$. \textit{Our goal} is to study the behaviour of the solution $\mathbf{u}_\varepsilon$ of the Cauchy problem for the hyperbolic system
\begin{equation}
\label{u_eps problem intr}
\begin{cases}
\partial _t^2 \mathbf{u}_\varepsilon (\mathbf{x},t) =-B_{\varepsilon}\mathbf{u}_\varepsilon (\mathbf{x},t),\quad\mathbf{x}\in\mathbb{R}^d,\;t\in\mathbb{R},
\\
\mathbf{u}_\varepsilon (\mathbf{x},0)=\boldsymbol{\phi}(\mathbf{x}),\quad \partial _t\mathbf{u}_\varepsilon(\mathbf{x},0)=\boldsymbol{\psi}(\mathbf{x}),
\end{cases}
\end{equation}
and to give an approximation for $\mathbf{u}_\varepsilon$ with the quantitative error estimate with explicit dependence of the norms of the initial data 
$\boldsymbol{\phi}\in H^2(\mathbb{R}^d;\mathbb{C}^n)$ and  $\boldsymbol{\psi}\in H^1(\mathbb{R}^d;\mathbb{C}^n)$. The solution $\mathbf{u}_\varepsilon$ is given by 
\begin{equation}
\label{u_eps= intr}
\mathbf{u}_\varepsilon (\cdot ,t)=\cos (tB_{\varepsilon}^{1/2})\boldsymbol{\phi}+B_{\varepsilon}^{-1/2}\sin(tB_{\varepsilon}^{1/2})\boldsymbol{\psi}.
\end{equation}
So, we are interested in approximations for the cosine $\cos (tB_{\varepsilon}^{1/2})$ and sine functions  $B_{\varepsilon}^{-1/2}\sin(tB_{\varepsilon}^{1/2})$ in suitable operator norms. (The terms ,,cosine'' and ,,sine functions'' used here are common in semigroup theory, see \cite[Definition 3.14.2 and (3.93)]{ABaHi}.)

\subsection{Main results}

It turns out that the solution $\mathbf{u}_\varepsilon$ of the problem \eqref{u_eps problem intr} with rapidly oscillating coefficients behaves in the small period limit as the solution of the so-called \textit{effective problem}:
\begin{equation*}
\begin{cases}
\partial _t^2 \mathbf{u}_0(\mathbf{x},t) =-B^0\mathbf{u}_0(\mathbf{x},t),\quad\mathbf{x}\in\mathbb{R}^d,\;t\in\mathbb{R},
\\
\mathbf{u}_0 (\mathbf{x},0)=\boldsymbol{\phi}(\mathbf{x}),\quad \partial _t\mathbf{u}_0(\mathbf{x},0)=\boldsymbol{\psi}(\mathbf{x}).
\end{cases}
\end{equation*}
Here $B^0$ is the~\textit{effective operator} with constant coefficients. The precise definition of $B^0$ can be found in Subsection~\ref{Subsection Effective operator}.

Let us formulate the main results of the paper. The principal term of approximation for the solution $\mathbf{u}_\varepsilon$ is given by
\begin{align*}
\Vert \mathbf{u}_\varepsilon (\cdot ,t)-\mathbf{u}_0(\cdot ,t)\Vert _{L_2(\mathbb{R}^d)}\leqslant C\varepsilon (1+\vert t\vert )\Vert \boldsymbol{\phi}\Vert _{H^2(\mathbb{R}^d)}+C\varepsilon\vert t\vert \Vert \boldsymbol{\psi}\Vert _{H^1(\mathbb{R}^d)}.
\end{align*}
It is also possible to approximate the time derivative $\partial _t\mathbf{u}_\varepsilon$ of the solution:
\begin{align*}
\Vert (\partial _t\mathbf{u}_\varepsilon) (\cdot ,t)-(\partial _t\mathbf{u}_0)(\cdot ,t)\Vert _{H^{-1}(\mathbb{R}^d)}\leqslant C\varepsilon (1+\vert t\vert )\Vert \boldsymbol{\phi}\Vert _{H^2(\mathbb{R}^d)}+C\varepsilon\vert t\vert \Vert \boldsymbol{\psi}\Vert _{H^1(\mathbb{R}^d)}.
\end{align*}

According to \eqref{u_eps= intr} and the similar identity for the solution $\mathbf{u}_0$ of the effective problem, these estimates can be rewritten in operator terms:
\begin{align}
\label{cos L2 intr}
\Vert &\cos(t B_\varepsilon ^{1/2})-\cos(t (B^0)^{1/2})\Vert _{H^2(\mathbb{R}^d)\rightarrow L_2(\mathbb{R}^d)}
\leqslant C\varepsilon(1+\vert t\vert),
\\
\label{sin L2 intr}
\Vert &B_\varepsilon ^{-1/2}\sin(t B_\varepsilon^{1/2})-(B^0)^{-1/2}\sin (t (B^0)^{1/2})\Vert _{H^1(\mathbb{R}^d)\rightarrow L_2(\mathbb{R}^d)}
\leqslant C \varepsilon \vert t\vert  ,
\\
\label{sin H-1 intr}
\begin{split}
\Vert &B_\varepsilon ^{1/2}\sin (tB_\varepsilon ^{1/2})-(B^0)^{1/2}\sin (t(B^0)^{1/2})\Vert _{H^2(\mathbb{R}^d)\rightarrow H^{-1}(\mathbb{R}^d)}
\leqslant C\varepsilon (1+\vert t\vert),
\end{split}
\\
\label{cos H-1 intr}
\begin{split}
\Vert &\cos(t B_\varepsilon ^{1/2})-\cos(t (B^0)^{1/2})\Vert _{H^1(\mathbb{R}^d)\rightarrow H^{-1}(\mathbb{R}^d)}
\leqslant C\varepsilon \vert t\vert .
\end{split}
\end{align}

In accordance with the classical result \cite{BrOFMu}, it is impossible to give $H^1$-approximation for the solution $\mathbf{u}_\varepsilon$ of the problem \eqref{u_eps problem intr} with an arbitrary initial data $\boldsymbol{\phi}$. In \cite{BrOFMu}, it was observed that such approximation can be constructed only for very special choices of the initial data. The argument is the following: convergence of the energy of the solution $\mathbf{u}_\varepsilon$ to the energy of $\mathbf{u}_0$ does not occur in the general situation. But the solution $\mathbf{u}_\varepsilon$  can be splitted into two parts: the first one is designed so that the corresponding energy converges to the energy for the effective equation and the second part tends to zero in some sense (but not in the uniform operator topology). In our considerations, we deal only with the first part ${\mathbf{v}}_\varepsilon $: 
\begin{equation}
\label{v problem introduction}
\begin{cases}
\partial _t^2 {\mathbf{v}}_\varepsilon  (\mathbf{x},t) =-B_{\varepsilon}{\mathbf{v}}_\varepsilon (\mathbf{x},t),\quad\mathbf{x}\in\mathbb{R}^d,\;t\in\mathbb{R},
\\
{\mathbf{v}}_\varepsilon  (\mathbf{x},0)=B_\varepsilon ^{-1}\boldsymbol{\phi}(\mathbf{x}),\quad \partial _t{\mathbf{v}}_\varepsilon  (\mathbf{x},0)=\boldsymbol{\psi}(\mathbf{x}).
\end{cases}
\end{equation}
Here $\boldsymbol{\phi}\in H^1(\mathbb{R}^d;\mathbb{C}^n)$ and  $\boldsymbol{\psi}\in H^2(\mathbb{R}^d;\mathbb{C}^n)$. 
In operator terms, this case corresponds to consideration of the operator $\cos (tB_\varepsilon ^{1/2})B_\varepsilon ^{-1}$ instead of $\cos (tB_\varepsilon ^{1/2})$ (see discussion in Subsection~\ref{Subsection Discussion} below). We have
\begin{align}
\label{cos corr intr}
\begin{split}
\Vert &\cos (tB_\varepsilon ^{1/2})B_\varepsilon ^{-1}-\cos(t (B^0)^{1/2})(B^0)^{-1}-\varepsilon \mathcal{K}_1(\varepsilon ;t)\Vert _{H^1(\mathbb{R}^d)\rightarrow H^1(\mathbb{R}^d)}
\\
&\leqslant C\varepsilon (1+\vert t\vert),
\end{split}
\\
\label{sin corr intr}
\begin{split}
\Vert &B_\varepsilon ^{-1/2}\sin(t B_\varepsilon ^{1/2})-(B^0)^{-1/2}\sin(t (B^0)^{1/2})-\varepsilon \mathcal{K}_2(\varepsilon ;t)\Vert _{H^2(\mathbb{R}^d)\rightarrow H^1(\mathbb{R}^d)}
\\
&\leqslant C \varepsilon (1+\vert t\vert).
\end{split}
\end{align}
Here $\mathcal{K}_1(\varepsilon ;t)$ and $\mathcal{K}_2(\varepsilon ;t)$ are correctors. They contain rapidly oscillating factors and so depend on $\varepsilon$. In the general case, correctors contain some smoothing operator. We distinguish the cases, when the smoothing operator can be removed from the correctors. In particular, if $d\leqslant 4$, then the smoothing operator always can be replaced by the identity operator.

Let ${\mathbf{v}}_0$ be the solution of the effective problem for \eqref{v problem introduction}:
\begin{equation*}
\begin{cases}
\partial _t^2 {\mathbf{v}}_0(\mathbf{x},t) =-B^0 {\mathbf{v}}_0(\mathbf{x},t),\quad\mathbf{x}\in\mathbb{R}^d,\;t\in\mathbb{R},
\\
{\mathbf{v}}_0 (\mathbf{x},0)=(B^0)^{-1}\boldsymbol{\phi}(\mathbf{x}),\quad \partial _t {\mathbf{v}}_0(\mathbf{x},0)=\boldsymbol{\psi}(\mathbf{x}).
\end{cases}
\end{equation*}
By ${\mathbf{w}}_\varepsilon$ we denote the first order approximation for the solution ${\mathbf{v}}_\varepsilon$:  $ {\mathbf{w}}_\varepsilon(\cdot ,t)= {\mathbf{v}}_0(\cdot,t)+\varepsilon \mathcal{K}_1(\varepsilon ;t)\boldsymbol{\phi}+\varepsilon \mathcal{K}_2(\varepsilon ;t)\boldsymbol{\psi}$. Then
\begin{align*}
&\Vert (\partial _t {\mathbf{v}}_\varepsilon)(\cdot ,t)-(\partial _t {\mathbf{v}}_0 )(\cdot ,t)\Vert _{L_2(\mathbb{R}^d)}
\leqslant C\varepsilon\vert t\vert \Vert \boldsymbol{\phi}\Vert _{H^1(\mathbb{R}^d)}+C\varepsilon (1+\vert t\vert)\Vert \boldsymbol{\psi}\Vert _{H^2(\mathbb{R}^d)},
\\
&\Vert {\mathbf{v}}_\varepsilon (\cdot ,t)-{\mathbf{w}}_\varepsilon(\cdot ,t)\Vert _{H^1(\mathbb{R}^d)}
\leqslant C\varepsilon (1+\vert t\vert )\Vert \boldsymbol{\phi}\Vert _{H^1(\mathbb{R}^d)}+C\varepsilon (1+\vert t\vert )\Vert \boldsymbol{\psi}\Vert _{H^2(\mathbb{R}^d)}.
\end{align*}
Note that our estimates grow with time as $O(\vert t\vert)$. Some growth seems natural because of the dispersion of waves in inhomogeneous media, see discussion in introduction to the paper \cite{BG} and references therein.

\subsection{Survey} 

The interest on homogenization results admitting a formulation in the uniform operator topology was stimulated by the work of M.~Sh.~Birman and T.~A.~Suslina \cite{BSu0}, where operator error estimates appeared in an explicit form at the first time. In \cite{BSu0}, the spectral theory approach to homogenization problems was applied. The method is based on the scaling transformation, the Floquet-Bloch theory, and the analytic perturbation theory. For further conceptual development of spectral method to the fibre homogenization, see \cite{CooW}. Another approach to obtaining operator error estimates for the problems in $\mathbb{R}^d$ was suggested by V.~V.~Zhikov \cite{Zh1} and developed by him together with S.~E.~Pastukhova \cite{ZhPas}. Both methods were applied to elliptic and parabolic problems. 
For the class of operators $B_\varepsilon$ under consideration, approximations for $B_\varepsilon ^{-1}$ and $e^{-tB_\varepsilon}$ were obtained in \cite{SuAA,SuAA14,M0}. The effective operator $B^0$ was calculated in \cite{Bo}.

Hyperbolic systems were studied only via the spectral approach, see \cite{BSu08,M1,M2,DSu17} and \cite{CheW}. In \cite{BSu08}, the non-stationary Schr\"odinger equation also was considered. It was obtained that
\begin{align}
\label{Th BSu cos introduction}
&\Vert \cos (t {A}_\varepsilon ^{1/2})-\cos (t( {A}^0)^{1/2})\Vert _{H^2(\mathbb{R}^d)\rightarrow L_2(\mathbb{R}^d)}\leqslant C\varepsilon (1+\vert t\vert) ,
\\
\label{Th BSu}
\begin{split}
&\Vert   {A}_\varepsilon^{-1/2}\sin( t  {A}^{1/2}_\varepsilon)-( {A}^0)^{-1/2}\sin ( t ( {A}^0)^{1/2})\Vert_{H^2(\mathbb{R}^d)\rightarrow L_2(\mathbb{R}^d)} \leqslant C\varepsilon (1+\vert t\vert )^2,
\end{split}
\\
\label{Th BSu Schrodinger}
&\Vert e^{-itA_\varepsilon}-e^{-it A^0}\Vert _{H^3(\mathbb{R}^d)\rightarrow L_2(\mathbb{R}^d)}\leqslant C\varepsilon (1+\vert t\vert ),
\end{align}
$t\in\mathbb{R}$. In \cite{D}, estimate \eqref{Th BSu Schrodinger} was generalized to the case of the operator $B_\varepsilon$. 
In \cite{M1,M2}, the estimate \eqref{Th BSu} was refined with respect to the type of the operator norm:
\begin{equation}
\begin{split}
\label{Th sin M principal term introduction}
\Vert {A}_\varepsilon ^{-1/2}\sin( t  {A}^{1/2}_\varepsilon)-( {A}^0)^{-1/2}\sin ( t ( {A}^0)^{1/2})\Vert_{H^1(\mathbb{R}^d)\rightarrow L_2(\mathbb{R}^d)} \leqslant C\varepsilon (1+\vert t\vert ), 
\end{split}
\end{equation}
$t\in\mathbb{R}$, 
and $(H^2\rightarrow H^1)$-approximation for the operator ${A}_\varepsilon^{-1/2}\sin( t  {A}^{1/2}_\varepsilon)$ was obtained
\begin{align}
\label{Th sin M corr introduction}
\begin{split}
&\left\Vert  {A}_\varepsilon ^{-1/2}\sin ( t  {A}_\varepsilon ^{1/2})-( {A}^0)^{-1/2}\sin ( t ( {A}^0)^{1/2})-\varepsilon\mathrm{K}(\varepsilon ; t)\right\Vert _{H^2(\mathbb{R}^d)\rightarrow H^1(\mathbb{R}^d)}
\\
&\leqslant C \varepsilon(1+\vert t\vert),\quad t\in\mathbb{R}.
\end{split}
\end{align}
Here $\mathrm{K}(\varepsilon ; t)$ is the corrector. 
In \cite{DSu17}, the sharpness of estimates \eqref{Th BSu cos introduction} and \eqref{Th sin M principal term introduction} with respect to the type of the norm was proven in the general case. By ,,sharpness'' we mean that it is  impossible to obtain the error estimate of the order $O(\varepsilon)$ for approximation of $\cos(tA_\varepsilon ^{1/2})$ in $(H^{2-\delta}\rightarrow L_2)$-norm with some $\delta >0$  and $(H^{1-\delta}\rightarrow L_2)$-estimate for $A_\varepsilon ^{-1/2}\sin (t A_\varepsilon ^{1/2})$. M.~Dorodnyi and T.~Suslina \cite{DSu20} showed that estimate \eqref{Th sin M corr introduction} is also sharp with respect to the type of the norm. Moreover, estimates \eqref{Th BSu cos introduction}, \eqref{Th sin M principal term introduction}, and \eqref{Th sin M corr introduction} are also sharp with respect to time, i.~e., the order $O(\vert t\vert)$, $\vert t\vert \rightarrow \infty$, cannot be refined in the general case. The proof available at \cite{DSu19}. For completeness of the survey, we mention the paper \cite{CheW}, where only the case $d=1$ was considered. For the problem of fractional elasticity, some order sharp estimate in $L_2$ was obtained in non-stationary setting as well as resolvent estimates (see \cite[Theorems 3.2 and 7.1]{CheW}).

Our estimates \eqref{cos L2 intr}, \eqref{sin L2 intr}, and \eqref{sin corr intr} transfer results from \cite{BSu08,M2} to the more general class of operators. So, we believe that estimates \eqref{cos L2 intr}, \eqref{sin L2 intr}, and \eqref{sin corr intr} are sharp with respect to the type of the norm.  
Note that inequalities \eqref{sin L2 intr} and \eqref{Th sin M principal term introduction} have different behavior for small $\vert t\vert$. 
Estimates of the form \eqref{sin H-1 intr}, \eqref{cos H-1 intr}, and \eqref{cos corr intr} are new even for the operator without lower order terms. It seems natural to expect that they are also sharp with respect to the type of the norm. Indeed,  
it is in accordance with usual difference in smoothness of the initial data for solution and its derivative in the setting of hyperbolic problems.

Recall that, according to the classical Trotter-Kato theorem (see, e.~g., \cite[Chapter~X, Theorem 1.1]{Sa}), the strong convergence of semigroups follows from the strong convergence of the corresponding resolvents. This result is fruitful for time-dependent homogenization problems. Let us note the recent work \cite{ChEl}, where the transfer of the Trotter-Kato theorem to  weak and uniform operator topologies was studied and the results were applied to homogenization of parabolic equations (without operator error estimates). Let us also mention the paper \cite{CooSav}, where the homogenization of the attractors of the quasi-linear damped wave equation was derived from the $(L_2\rightarrow L_2)$-approximation for the resolvent. But the results of \cite{CooSav} can not be written in the uniform operator topology: the error estimates for homogenization are written in terms of the attractors, not the solutions\footnote{I.~e., there are no operator-error analogues for Theorem 1.3 from \cite{FiVi}, where the quantitative homogenization of attaractors was studied in a more classical manner.}, and the dependence of constants from the norms of the initial data was not traced. For operators, acting in a bounded domain, hyperbolic results were derived from elliptic ones via the Laplace transform, see \cite{M4}. But the corresponding results do not look optimal with respect to the type of the operator norm. Let us also mention the work \cite{W}, where Laplace transform technique also was applied, but the rate of convergence was not traced explicitly. Finally, we highlight the paper \cite{CheW}, where the non-stationary homogenization result with the precise order error estimate was derived from elliptic one with the help of the Fourier–Laplace transform. 

The method of the present paper also can be considered as a variation on the theme of the Trotter-Kato theorem.
 
\subsection{Method} 
\label{Subsection Method}

The proof relies on the rewriting of the hyperbolic equation as a system for the solution and its time derivative, the explicit formulas for the corresponding group of operators and the resolvent of the generator of the group, and the known results on homogenization of the resolvent $B_\varepsilon ^{-1}$ from \cite{SuAA}. Our key observation is that the quantitative homogenization results for the hyperbolic systems can be derived from the approximation of the resolvent by analogy with the proof of the Trotter-Kato theorem. 

First, we rewrite hyperbolic problem \eqref{u_eps problem intr} in the following form:
\begin{equation*}
\partial _t \begin{pmatrix}
  \mathbf{u}_\varepsilon\\
  \partial _t\mathbf{u}_\varepsilon
\end{pmatrix}
=\begin{pmatrix}
  0& I\\
 -B_{\varepsilon}& 0
\end{pmatrix}
\begin{pmatrix}
  \mathbf{u}_\varepsilon\\
  \partial _t\mathbf{u}_\varepsilon
\end{pmatrix},
\quad
\begin{pmatrix}
  \mathbf{u}_\varepsilon (\cdot ,0)\\
  \partial _t\mathbf{u}_\varepsilon(\cdot ,0)
\end{pmatrix}
=
\begin{pmatrix}
  \boldsymbol{\phi}\\
  \boldsymbol{\psi}
\end{pmatrix}.
\end{equation*}
Denote $\mathfrak{A}_\varepsilon = \begin{pmatrix}
  0& I\\
 -B_{\varepsilon}& 0
\end{pmatrix}$. Then, according to \eqref{u_eps= intr}, 
\begin{equation*}
\begin{pmatrix}
  \mathbf{u}_\varepsilon\\
  \partial _t\mathbf{u}_\varepsilon
\end{pmatrix}
=
e^{t\mathfrak{A}_\varepsilon}\begin{pmatrix}
  \boldsymbol{\phi}\\
  \boldsymbol{\psi}
\end{pmatrix}
=
\begin{pmatrix}
  \cos(tB_{\varepsilon}^{1/2})& B_{\varepsilon}^{-1/2}\sin (tB_{\varepsilon}^{1/2})\\
  -B_{\varepsilon}^{1/2}\sin (tB_{\varepsilon}^{1/2})& \cos(tB_{\varepsilon}^{1/2})
\end{pmatrix}
\begin{pmatrix}
  \boldsymbol{\phi}\\
  \boldsymbol{\psi}
\end{pmatrix}.
\end{equation*}

Our second step can be considered as technical modification of the proof of the classical Trotter-Kato theorem in the variant from \cite[Lemma~19.1]{BSu08}. The new observation is the following: we need to use approximation for the resolvent with corrector to obtain the principal term of approximation for the solution of hyperbolic problem. We use this trick to make estimates optimal with respect to the type of the operator norm. 

Actually, estimates \eqref{cos L2 intr}--\eqref{cos H-1 intr} are neither else then approximations for $e^{t\mathfrak{A}_\varepsilon}$ in $H^2\times H^1\rightarrow L^2\times H^{-1}$-norm written entrywise. While in the context of works \cite{BSu08}, \cite{DSu17}, \cite{DSu19} results in norms of operators acting to $H^{-1}$ look new, the $L_2\times H^{-1}$-norm appeared in a more classical homogenization context, see \cite[(1.13), (1.22)]{FiVi}.

For the operator $e^{-itB_\varepsilon}$, the method simplifies significantly. It allows us to give a short proof of $(H^3\rightarrow L_2)$-approximation of the unitary group $e^{-itB_\varepsilon}$ (the original proof belongs to M.~Dorodnyi \cite{D}) and approximation with corrector for the operator $e^{-itB_\varepsilon}B_\varepsilon ^{-1}$.
 
\subsection{Plan of the paper} The paper consists of five sections and introduction. Section~\ref{Section Preliminaries} contains definitions of the operators $B_\varepsilon$ and $B^0$ and known results on homogenization of the resolvent $B_\varepsilon ^{-1}$. In Section~\ref{Section Problem setting. Main results}, the main results of the paper are formulated and proven. In Section~\ref{Section Removal of the smoothing operator from the corrector}, we discuss the possibility to remove the smoothing operator from the corrector. In Section~\ref{Section Homogenization of solutions of inhomogeneous hyperbolic systems}, we apply results in operator terms to homogenization of the solutions of hyperbolic systems. Finally, in Section~\ref{Section Schrodinger}, we give alternative proof for approximation of the operator $e^{-itB_\varepsilon}$ originally obtained in \cite{D} and, for completeness of the presentation, approximation in the energy norm for $e^{-itB_\varepsilon}B_\varepsilon ^{-1}$.

\subsection{Notation} Let $\mathfrak{H}$ and $\mathfrak{H}_\bullet$ be complex separable Hilbert spaces. The symbols $(\cdot ,\cdot)_\mathfrak{H}$ and $\Vert \cdot\Vert _\mathfrak{H}$ denote the inner product and the norm in   $\mathfrak{H}$, respectively; the symbol $\Vert \cdot\Vert _{\mathfrak{H}\rightarrow\mathfrak{H}_\bullet}$ means the norm of the linear continuous operators from $\mathfrak{H}$ to $\mathfrak{H}_\bullet$. 
The algebra of all bounded linear operators acting in $\mathfrak{H}$ is denoted by $\mathcal{B}(\mathfrak{H})$. 

The symbols $\langle \cdot ,\cdot\rangle$ and $\vert \cdot\vert$ stand for the inner product and the norm in   $\mathbb{C}^n$, respectively, $\mathbf{1}_n$ is the identity $(n\times n)$-matrix. If $a$ is an~$(m\times n)$-matrix, then the symbol $\vert a\vert$ denotes the norm of the matrix $a$ as the operator from  $\mathbb{C}^n$ to $\mathbb{C}^m$. 

We use the notation $\mathbf{x}=(x_1,\dots , x_d)\in\mathbb{R}^d$, $iD_j=\partial _j =\partial /\partial x_j$, $j=1,\dots,d$, $\mathbf{D}=-i\nabla=(D_1,\dots ,D_d)$. The classes $L_p$ of vector-valued functions in a domain $\mathcal{O}\subset\mathbb{R}^d$ with values in $\mathbb{C}^n$ are denoted by $L_p(\mathcal{O};\mathbb{C}^n)$,\break $1\leqslant p\leqslant \infty$. The Sobolev spaces of $\mathbb{C}^n$-valued functions in a domain $\mathcal{O}\subset\mathbb{R}^d$ are denoted by $H^s(\mathcal{O};\mathbb{C}^n)$. 
For $n=1$, we simply write  $L_p(\mathcal{O})$, $H^s(\mathcal{O})$ and so on, but, sometimes, if this does not lead to confusion, we use such simple notation for the spaces of vector-valued or matrix-valued functions. The symbol $L_p((0,T);\mathfrak{H})$, $1\leqslant p\leqslant\infty$, denotes the $L_p$-space of $\mathfrak{H}$-valued functions on the interval $(0,T)$.

Various constants in estimates are denoted by $c$, $\mathfrak{c}$, $C$, and $\mathfrak{C}$ 
(possibly, with indices and marks).

\subsection*{Acknowledgement} 
The author is happy to thank Institut Mittag-Leffler 
for financial support and hospitality 
during the research program ,,Spectral Methods in Mathematical Physics.'' Stimulating environment contributed to completing this work.

Particular gratitude goes to T.~A.~Suslina for drawing the author's attention to Lemma~19.1 from \cite{BSu08}. The author is appreciate to T.~A.~Suslina for the comments improved the quality of presentation. 

The author is also grateful to an anonymous colleague who pointed out that the notation $\widetilde{\quad}$ is overloaded\footnote{To simplify the notation originally coming from the cited papers of M.~Sh.~Birman and T.~A.~Suslina, the replacement of $\widetilde{\Gamma}$ by $\Gamma^*$, $\widetilde{\Omega}$ by $\Omega ^*$, and $\widetilde{H}^1$ by $H^1_\mathrm{per}$  
was carried out.} and found some inaccuracies in the author's English. 

\section{Preliminaries. Known results}
\label{Section Preliminaries}

The main material of the present section contained in \cite{SuAA,SuAA14}. For the reader convenience, we repeat it in details.

\subsection{Lattices in $\mathbb{R}^d$} 
\label{Subsection lattices}
Let $\Gamma \subset \mathbb{R}^d$ be a lattice generated by a basis  $\mathbf{a}_1,\dots ,\mathbf{a}_d \in \mathbb{R}^d$, i.~e., 
$$
\Gamma =\left\lbrace
\mathbf{a}\in \mathbb{R}^d : \mathbf{a}=\sum _{j=1}^d \nu _j \mathbf{a}_j,\; \nu _j\in \mathbb{Z}
\right\rbrace ,
$$
and let $\Omega$ be the elementary cell of $\Gamma$:
$$
\Omega =
\left \lbrace
\mathbf{x}\in \mathbb{R}^d :\mathbf{x}=\sum _{j=1}^d \tau _j \mathbf{a}_j , \; -\frac{1}{2}<\tau _j<\frac{1}{2}
\right\rbrace .
$$
By $\vert \Omega \vert $ we denote the Lebesgue measure of $\Omega$: $\vert \Omega \vert =\mathrm{mea
s}\,\Omega$.

The basis $\mathbf{b}_1,\dots ,\mathbf{b}_d$ in $\mathbb{R}^d$ dual to the basis $\mathbf{a}_1,\dots ,\mathbf{a}_d$
is defined by the relations $\langle \mathbf{b}_i ,\mathbf{a}_j \rangle =2\pi \delta _{ij}$. 
The lattice ${\Gamma}^*$ generated by the dual basis is called the lattice dual to $\Gamma $. 
We consider the first Brillouin zone
 $$
 {\Omega}^*=\lbrace \mathbf{k}\in \mathbb{R}^d : \vert \mathbf{k}\vert < \vert \mathbf{k}-\mathbf{b}\vert , \; 0\neq \mathbf{b}\in {\Gamma}^* \rbrace,
 $$
as a fundamental domain of the lattice ${\Gamma}^*$. By $r_0$ we denote the radius of the ball inscribed in $\mathrm{clos}\, {\Omega}^*$. Then $2r_0=\min _{0\neq\mathbf{b}\in {\Gamma}^*}\vert \mathbf{b}\vert $.

By ${H}^1_{\mathrm{per}}(\Omega)$ we denote the subspace of all functions in 
 $H^1(\Omega)$ whose $\Gamma$-periodic extension to $\mathbb{R}^d$ belongs to $H^1_{\mathrm{loc}}(\mathbb{R}^d)$. If $f (\mathbf{x})$~is a~$\Gamma$-periodic measurable matrix-valued function in $\mathbb{R}^d$, we put $f ^\varepsilon (\mathbf{x}):=f (\mathbf{x}/\varepsilon)$, $\varepsilon >0$; $\overline{f}:=\vert \Omega\vert ^{-1}\int _\Omega f(\mathbf{x})\,d\mathbf{x}$, and $\underline{f}:=\left(\vert \Omega\vert ^{-1}\int _\Omega f(\mathbf{x})^{-1}\,d\mathbf{x}\right)^{-1}$. Here, in the definition of  $\overline{f}$ it is assumed that $f\in L_{1,\mathrm{loc}}(\mathbb{R}^d)$, and in the definition of $\underline{f}$ it is assumed that the matrix $f(\mathbf{x})$ is square and nondegenerate, and  $f^{-1}\in L_{1,\mathrm{loc}}(\mathbb{R}^d)$. By $[f^\varepsilon ]$ we denote the operator of multiplication by the matrix-valued function $f^\varepsilon (\mathbf{x})$.
 
Let $\mathcal{U}$ be the Gelfand transformation associated with the lattice $\Gamma$. Initially, $\mathcal{U}$ is defined on the Schwartz class $\mathcal{S}(\mathbb{R}^d;\mathbb{C}^n)$ by the formula 
\begin{equation*}
(\mathcal{U}\mathbf{v})(\mathbf{k},\mathbf{x})=\vert {\Omega}^*\vert ^{-1/2}\sum _{\mathbf{a}\in\Gamma }\exp (-i\langle \mathbf{k},\mathbf{x}+\mathbf{a}\rangle )\mathbf{v}(\mathbf{x}+\mathbf{a}),\quad
\mathbf{x}\in\Omega,\quad\mathbf{k}\in {\Omega}^*.
\end{equation*}
Herewith, $\int _{{\Omega}^*}\int _\Omega \vert (\mathcal{U}\mathbf{v})(\mathbf{k},\mathbf{x})\vert ^2\,d\mathbf{x}\,d\mathbf{k}=\int _{\mathbb{R}^d}\vert \mathbf{v}(\mathbf{x})\vert ^2\,d\mathbf{x}$, and $\mathcal{U}$ extends by continuity to a unitary operator 
$$\mathcal{U}:L_2(\mathbb{R}^d;\mathbb{C}^n)\rightarrow \int _{{\Omega}^*}\oplus L_2(\Omega;\mathbb{C}^n)\,d\mathbf{k}$$ 
(for more details, see \cite[Chapter~2, Subsection 1.3]{BSu}). 
 
 \subsection{The class of operators $A_\varepsilon$}
\label{Subsection Operator A_eps}
In $L_2(\mathbb{R}^d;\mathbb{C}^n)$, we consider the operator $A_\varepsilon$ given by the differential expression  $A_\varepsilon =b(\mathbf{D})^*g^\varepsilon (\mathbf{x})b(\mathbf{D})$. Here $g(\mathbf{x})$ is a $\Gamma$-periodic $(m\times m)$-matrix-valued function (in general, with complex entries). We assume that $g(\mathbf{x})>0$ and $g,g^{-1}\in L_\infty (\mathbb{R}^d)$. Next, $b(\mathbf{D})$ is the differential operator given  by 
\begin{equation}
\label{b(D)=}
b(\mathbf{D})=\sum _{j=1}^d b_jD_j, 
\end{equation}
where $b_j$, $j=1,\dots ,d$, are constant $(m\times n)$-matrices (in general, with complex entries). It is assumed that  $m\geqslant n$ and that the symbol $b(\boldsymbol{\xi})=\sum _{j=1}^d b_j\xi_j$ of the operator $b(\mathbf{D})$ has maximal rank:
$$
\mathrm{rank}\,b(\boldsymbol{\xi})=n,\quad 0\neq \boldsymbol{\xi}\in\mathbb{R}^d.
$$
This condition is equivalent to the estimates 
\begin{equation}
\label{<b^*b<}
\alpha _0\mathbf{1}_n \leqslant b(\boldsymbol{\theta})^*b(\boldsymbol{\theta})
\leqslant \alpha _1\mathbf{1}_n,\quad
\boldsymbol{\theta}\in \mathbb{S}^{d-1},\quad
0<\alpha _0\leqslant \alpha _1<\infty,
\end{equation}
with some positive constants $\alpha _0$ and $\alpha _1$.
From \eqref{<b^*b<} it follows that  
\begin{equation}
\label{b_l <=}
\vert b_j\vert \leqslant \alpha _1^{1/2},\quad j=1,\dots ,d.
\end{equation}

The precise definition of the operator $A_\varepsilon$ is given in terms of the quadratic form 
\begin{equation*}
\mathfrak{a}_\varepsilon [\mathbf{u},\mathbf{u}]=\int _{\mathbb{R}^d}\langle g^\varepsilon (\mathbf{x})b(\mathbf{D})\mathbf{u},b(\mathbf{D})\mathbf{u}\rangle \,d\mathbf{x},\quad \mathbf{u}\in H^1(\mathbb{R}^d;\mathbb{C}^n).
\end{equation*}
Under the above assumptions, this form is closed and nonnegative. 
Using the Fourier transformation and condition \eqref{<b^*b<}, it is easy to check that
\begin{equation}
\label{<a_eps<}
\alpha _0\Vert g^{-1}\Vert ^{-1}_{L_\infty}\Vert \mathbf{D}\mathbf{u}\Vert ^2_{L_2(\mathbb{R}^d)}\leqslant \mathfrak{a}_\varepsilon [\mathbf{u},\mathbf{u}]\leqslant \alpha _1\Vert g\Vert _{L_\infty}\Vert \mathbf{D}\mathbf{u}\Vert^2_{L_2(\mathbb{R}^d)},\quad
 \mathbf{u}\in H^1(\mathbb{R}^d;\mathbb{C}^n).
\end{equation}
Let $c_1:=\alpha _0 ^{-1/2}\Vert g^{-1}\Vert ^{1/2}_{L_\infty}$. Then the lower estimate 
 \eqref{<a_eps<} can be written as 
 \begin{equation}
\label{Du <= c_1^2a}
\Vert\mathbf{D}\mathbf{u}\Vert ^2_{L_2(\mathbb{R}^d)}\leqslant c_1^2 \mathfrak{a}_\varepsilon [\mathbf{u},\mathbf{u}],\quad \mathbf{u}\in H^1(\mathbb{R}^d;\mathbb{C}^n).
\end{equation}

\subsection{Lower order terms}
\label{Subsection lower order terms}
We study a selfadjoint operator $B_\varepsilon$ whose principal part coincides with $A_\varepsilon$.
To define the lower order terms, we introduce $\Gamma$-periodic $(n\times n)$-matrix-valued functions $a_j$, $j=1,\dots ,d$, (in general, with complex entries) such that 
\begin{align*}
a_j \in L_\rho (\Omega ), \quad\rho =2 \;\mbox{if}\;d=1,\quad\rho >d\;\mbox{if}\;d\geqslant 2,\quad  j=1,\dots ,d.
\end{align*}

Let $\mathcal{Q}_0$ be the operator in $L_2(\mathbb{R}^d;\mathbb{C}^n)$ 
that acts as multiplication by the $\Gamma$-periodic positive definite and 
bounded $(n\times n)$-matrix-valued function $\mathcal{Q}_0(\textbf{x})$. We factorise this matrix as 
$\mathcal{Q}_0(\textbf{x})=f(\textbf{x})^*f(\textbf{x})$, where an $(n\times n)$-matrix-valued function $f(\textbf{x})$ is assumed to be $\Gamma$-periodic. 

Suppose that $d\mu (\textbf{x})$ is a $\Gamma$-periodic Borel
$\sigma $-finite measure in $\mathbb{R}^d$ with va\-lues in the class
of Hermitian  $(n\times n)$-matrices. Then
$d\mu (\textbf{x})=\lbrace d\mu 
_{jl}(\textbf{x})\rbrace ,\,j,l=1,\dots ,n$. 
In other words, $d\mu _{jl}(\textbf{x})$ is
a complex-valued $\Gamma$-periodic measure in $\mathbb{R}^d$, and $d\mu 
_{jl}=d\mu _{lj}^*$. Suppose that the measure $d\mu$ is such that a function 
$\vert v (\textbf{x})\vert ^2$ is integrable with respect to each measure $d\mu _{jl}$ for 
any $v\in H^1(\mathbb{R}^d)$.

In $L_2(\mathbb{R}^d;\mathbb{C}^n)$, we consider the sesquilinear form
\begin{equation*}
\int _{\mathbb{R}^d}\langle d\mu (\mathbf{x})\mathbf{u},\mathbf{v}\rangle
=\sum _{j,l=1}^n \int _{\mathbb{R}^d}u_lv_j^*d\mu _{jl}(\mathbf{x}),\quad \mathbf{u},\mathbf{v}\in H^1(\mathbb{R}^d;\mathbb{C}^n).
\end{equation*}
 
Assume that the measure $d\mu $ is subject to 
the following condition.

\begin{condition}
\label{Condition on mu}
$1^\circ$. There exist constants $\widetilde{c}_2\geqslant 0$ and $c_3\geqslant 0$ such that for any $\mathbf{u},\mathbf{v}\in H^1(\Omega;\mathbb{C}^n)$ we have
\begin{equation}
\label{cond mu 1}
\Bigl\vert \int _\Omega\langle d\mu (\mathbf{x})\mathbf{u},\mathbf{v}\rangle\Bigr\vert
\leqslant
(\widetilde{c}_2\Vert \mathbf{D}\mathbf{u}\Vert ^2 _{L_2(\Omega)}+c_3\Vert \mathbf{u}\Vert ^2_{L_2(\Omega)})^{1/2}
(\widetilde{c}_2\Vert \mathbf{D}\mathbf{v}\Vert ^2 _{L_2(\Omega)}+c_3\Vert \mathbf{v}\Vert ^2_{L_2(\Omega)})^{1/2}.
\end{equation}

$2^\circ$. We have
\begin{equation}
\label{cond mu 2}
\int _\Omega \langle d\mu (\mathbf{x})\mathbf{u},\mathbf{u}\rangle\geqslant -\widetilde{c}\Vert \mathbf{D}\mathbf{u}\Vert ^2 _{L_2(\Omega)}-c_0\Vert \mathbf{u}\Vert ^2 _{L_2(\Omega)},\quad\mathbf{u}\in H^1(\Omega;\mathbb{C}^n),
\end{equation}
with constants ${c}_0\in\mathbb{R}$  and $\widetilde{c}$ such that 
$0\leqslant \widetilde{c}<\alpha _0\Vert g^{-1}\Vert ^{-1}_{L_\infty}$.
\end{condition}

Examples of the forms satisfying Condition~\ref{Condition on mu} can be found in \cite[Subsection 5.5]{SuAA}. Here we provide only the main example.

\begin{example}
\label{Example of mu}
\textnormal{
Suppose that the measure $d\mu$ is absolutely continuous with respect to Lebesgue measure, i.~e., $d\mu(\mathbf{x})=Q(\mathbf{x})\,d\mathbf{x}$, where $Q(\mathbf{x})$ is a $\Gamma$-periodic Hermitian $(n\times n)$-matrix-valued function in $\mathbb{R}^d$ such that
\begin{equation*}
Q\in L_\varrho(\Omega),\quad \varrho=1\;\mbox{for}\;d=1,\quad \varrho>\frac{d}{2}\;\mbox{for}\;d\geqslant 2.
\end{equation*}
Then, by the Sobolev embedding theorem, for any $\nu >0$ there exists a positive constant $C_Q(\nu)$ such that
\begin{equation*}
\int _\Omega \vert Q(\mathbf{x})\vert \vert \mathbf{u}\vert ^2\,d\mathbf{x}
\leqslant \nu \int _\Omega \vert \mathbf{D}\mathbf{u}\vert ^2\,d\mathbf{x}
+C_Q(\nu)\int _\Omega \vert \mathbf{u}\vert ^2\,d\mathbf{x},\quad
\mathbf{u}\in H^1(\Omega ;\mathbb{C}^n).
\end{equation*}
Then Condition~\ref{Condition on mu} is satisfied with the constants $\widetilde{c}_2=1$, $c_3=C_Q(1)$, $\widetilde{c}=\nu$, and $c_0=C_Q(\nu)$, where $2\nu=\alpha _0\Vert g^{-1}\Vert ^{-1}_{L_\infty}$.
}
\end{example}

\subsection{The operator $B_\varepsilon$} 

Define the measure $d\mu ^\varepsilon(\mathbf{x})$ as follows. For any Borel set $\Delta\subset\mathbb{R}^d$, we consider the set 
$\varepsilon ^{-1}\Delta =\lbrace\mathbf{y}=\varepsilon ^{-1}\mathbf{x} : \mathbf{x}\in\Delta\rbrace$ and put $\mu ^\varepsilon (\Delta)=\varepsilon ^d\mu (\varepsilon ^{-1}\Delta )$. Consider the sesquilinear form $\mathfrak{q}_\varepsilon$ defined by
\begin{equation*}
\mathfrak{q}_\varepsilon [\mathbf{u},\mathbf{v}]=\int _{\mathbb{R}^d}
\langle d\mu ^\varepsilon (\mathbf{x})\mathbf{u},\mathbf{v}\rangle 
,\quad \mathbf{u},\mathbf{v}\in H^1(\mathbb{R}^d;\mathbb{C}^n).
\end{equation*}

Consider the following quadratic form
\begin{equation}
\label{b_eps =}
\begin{split}
\mathfrak{b}_\varepsilon [\mathbf{u},\mathbf{u}]&=\mathfrak{a}_\varepsilon [\mathbf{u},\mathbf{u}]+2\mathrm{Re}\,\sum _{j=1}^d(a_j^\varepsilon D_j \mathbf{u},\mathbf{u})_{L_2(\mathbb{R}^d)}\\
&+ \mathfrak{q}_\varepsilon[\mathbf{u},\mathbf{u}]+\lambda (Q_0^\varepsilon \mathbf{u},\mathbf{u})_{L_2(\mathbb{R}^d)},\quad
 \mathbf{u}\in H^1(\mathbb{R}^d;\mathbb{C}^n).
\end{split}
\end{equation}
We fix a constant $\lambda$ so that the form $\mathfrak{b}_\varepsilon$ is nonnegative (see \eqref{lambda condition} below).
Let us check that the form $\mathfrak{b}_\varepsilon$ is closed. 
By the H\"older inequality and the Sobolev embedding theorem, it is easily seen  (see  \cite[(5.11)--(5.14)]{SuAA}) that 
for any $\nu>0$ there exist constants $C_j(\nu)>0$ such that 
\begin{equation*}
\Vert a_j^*\mathbf{u}\Vert ^2 _{L_2(\mathbb{R}^d)}\leqslant \nu \Vert \mathbf{D}\mathbf{u}\Vert ^2_{L_2(\mathbb{R}^d)}+C_j(\nu )\Vert \mathbf{u}\Vert ^2 _{L_2(\mathbb{R}^d)},
\quad \mathbf{u}\in H^1(\mathbb{R}^d;\mathbb{C}^n),\;j=1,\dots ,d.
\end{equation*}
Using the change of variable $\mathbf{y}:=\varepsilon ^{-1}\mathbf{x}$ and denoting $\mathbf{u}(\mathbf{x})=:\mathbf{v}(\mathbf{y})$, we deduce 
\begin{equation}
\label{a_j^eps u<=}
\begin{split}
\Vert &  (a_j^\varepsilon )^*\mathbf{u}\Vert ^2_{L_2(\mathbb{R}^d)}=\int _{\mathbb{R}^d}\vert a_j(\varepsilon ^{-1}\mathbf{x})^*\mathbf{u}(\mathbf{x})\vert ^2\,d\mathbf{x}
=\varepsilon ^d\int _{\mathbb{R}^d}\vert a_j(\mathbf{y})^*\mathbf{v}(\mathbf{y})\vert ^2\,d\mathbf{y}\\
&\leqslant \varepsilon ^d\nu \int _{\mathbb{R}^d}\vert \mathbf{D}_{\mathbf{y}}\mathbf{v}(\mathbf{y})\vert ^2\,d\mathbf{y}
+\varepsilon ^d C_j(\nu)\int _{\mathbb{R}^d}\vert \mathbf{v}(\mathbf{y})\vert ^2\,d\mathbf{y}\\
&\leqslant \nu \Vert \mathbf{D}\mathbf{u}\Vert ^2_{L_2(\mathbb{R}^d)}+C_j(\nu)\Vert \mathbf{u}\Vert ^2_{L_2(\mathbb{R}^d)},\quad \mathbf{u}\in H^1(\mathbb{R}^d;\mathbb{C}^n),\quad 0<\varepsilon\leqslant 1.
\end{split}
\end{equation}
Hence, by \eqref{<a_eps<}, for any $\nu >0$ there exists a constant $C(\nu)>0$ such that 
\begin{equation}
\label{sum a-j u}
\begin{split}
\sum _{j=1}^d \Vert (a_j^\varepsilon)^*\mathbf{u}\Vert ^2 _{L_2(\mathbb{R}^d)}
\leqslant \nu
\mathfrak{a}_\varepsilon [\mathbf{u},\mathbf{u}]
+C(\nu)\Vert \mathbf{u}\Vert ^2_{L_2(\mathbb{R}^d)},
\\
\mathbf{u}\in H^1(\mathbb{R}^d;\mathbb{C}^n),
\quad 0<\varepsilon\leqslant 1.
\end{split}
\end{equation}
If $\nu$ is fixed, then $C(\nu)$ depends only on $d$, $\rho$, $\alpha _0$,  the norms $\Vert g^{-1}\Vert _{L_\infty}$, $\Vert a_j\Vert _{L_\rho (\Omega)}$, $j=1,\dots ,d$, and the parameters of the lattice $\Gamma$.

For functions in $H^1(\mathbb{R}^d;\mathbb{C}^n)$, we write inequalities \eqref{cond mu 1}, \eqref{cond mu 2} over shifted cells $\Omega +\mathbf{a}$, $\mathbf{a}\in\Gamma$, and sum up, obtaining similar inequalities with integration over $\mathbb{R}^d$. Using these arguments, changing variables $\mathbf{y}:=\varepsilon ^{-1}\mathbf{x}$ similarly to \eqref{a_j^eps u<=}, and taking \eqref{<a_eps<} into account, we can obtain that
\begin{align}
\label{q_eps<=}
&\vert \mathfrak{q} _\varepsilon [\mathbf{u},\mathbf{v}]\vert \leqslant (c_2 \mathfrak{a}_\varepsilon [\mathbf{u},\mathbf{u}]+c_3\Vert \mathbf{u}\Vert ^2_{L_2(\mathbb{R}^d)})^{1/2}
(c_2\mathfrak{a}_\varepsilon [\mathbf{v},\mathbf{v}]+c_3\Vert \mathbf{v}\Vert ^2_{L_2(\mathbb{R}^d)})^{1/2},
\\
\label{q-eps>=}
&\mathfrak{q }_\varepsilon [\mathbf{u},\mathbf{u}]\geqslant -(1-\kappa)\mathfrak{a}_\varepsilon [\mathbf{u},\mathbf{u}] -c_0\Vert \mathbf{u}\Vert ^2 _{L_2(\mathbb{R}^d)},\quad \mathbf{u},\mathbf{v}\in H^1(\mathbb{R}^d;\mathbb{C}^n),
\end{align}
$0<\varepsilon\leqslant 1$. Here
\begin{equation*}
c_2=\widetilde{c}_2\alpha _0^{-1}\Vert g^{-1}\Vert _{L_\infty},\quad
\kappa =1-\widetilde{c}\alpha _0^{-1}\Vert g^{-1}\Vert _{L_\infty},\quad 0<\kappa\leqslant 1.
\end{equation*}

Combining \eqref{Du <= c_1^2a}, \eqref{sum a-j u} with $\nu=1$, and \eqref{q_eps<=}, we get
\begin{equation}
\label{b_eps <=}
\begin{split}
\mathfrak{b}_\varepsilon [\mathbf{u},\mathbf{u}]
\leqslant (2+c_1^2+c_2)\mathfrak{a}_\varepsilon [\mathbf{u},\mathbf{u}]+(C(1)+c_3+\vert\lambda\vert \Vert Q_0\Vert _{L_\infty})\Vert \mathbf{u}\Vert ^2 _{L_2(\mathbb{R}^d)},
\\
\mathbf{u}\in H^1(\mathbb{R}^d;\mathbb{C}^n),\quad 0<\varepsilon\leqslant 1.
\end{split}
\end{equation}

From \eqref{Du <= c_1^2a} and \eqref{sum a-j u} it follows that
\begin{equation}
\label{2Re sum j <=}
\begin{split}
2\biggl\vert \mathrm{Re}\,\sum _{j=1}^d (D_j\mathbf{u},(a_j^\varepsilon)^*\mathbf{u})_{L_2(\mathbb{R}^d)}\biggr\vert
\leqslant\frac{\kappa }{2}\mathfrak{a}_\varepsilon [\mathbf{u},\mathbf{u}]+c_4\Vert \mathbf{u}\Vert ^2_{L_2(\mathbb{R}^d)},
\\
\mathbf{u}\in H^1(\mathbb{R}^d;\mathbb{C}^n),\  0<\varepsilon\leqslant 1,
\end{split}
\end{equation}
where $c_4=4 \kappa ^{-1} c_1^2C(\nu _0)$ with $\nu _0=16^{-1}c_1^{-2}\kappa ^2$.

Assume that the parameter $\lambda $ is subject to the following restriction:
\begin{equation}
\label{lambda condition}
\begin{split}
&\lambda >\Vert Q_0^{-1}\Vert _{L_\infty}(c_0+c_4)\;\mbox{if}\;\lambda\geqslant 0,
\\
&\lambda >\Vert Q_0\Vert ^{-1}_{L_\infty}(c_0+c_4)\;\mbox{if}\;\lambda <0\;(\mbox{and}\;c_0+c_4<0).
\end{split}
\end{equation}
This condition ensures that
\begin{equation}
\label{lambda Q_0 est}
\lambda (Q_0\mathbf{u},\mathbf{u})_{L_2(\mathbb{R}^d)}
\geqslant (c_0+c_4+\beta)\Vert \mathbf{u}\Vert ^2_{L_2(\mathbb{R}^d)},\quad
\mathbf{u}\in L_2(\mathbb{R}^d;\mathbb{C}^n),
\end{equation}
where $\beta >0$ is defined in terms of $\lambda$ as follows
\begin{equation*}
\begin{split}
&\beta =\lambda \Vert Q_0^{-1}\Vert ^{-1}_{L_\infty}-c_0-c_4\;\mbox{if}\;\lambda\geqslant 0,
\\
&\beta =\lambda \Vert Q_0\Vert _{L_\infty}-c_0-c_4\;\mbox{if}\;\lambda <0\;(\mbox{and}\;c_0+c_4<0).
\end{split}
\end{equation*}
Then from \eqref{b_eps =}, \eqref{q-eps>=}, \eqref{2Re sum j <=}, and \eqref{lambda Q_0 est} it follows that
\begin{equation}
\label{b_eps >=}
\mathfrak{b}_\varepsilon [\mathbf{u},\mathbf{u}]\geqslant \frac{\kappa}{2}\mathfrak{a}_\varepsilon [\mathbf{u},\mathbf{u}]+\beta \Vert \mathbf{u}\Vert ^2 _{L_2(\mathbb{R}^d)},\quad \mathbf{u}\in H^1(\mathbb{R}^d;\mathbb{C}^n),\quad 0<\varepsilon\leqslant 1.
\end{equation}

Thus, the form $\mathfrak{b}_\varepsilon$ is closed and positive definite. 
The selfadjoint operator in $L_2(\mathbb{R}^d;\mathbb{C}^n)$ generated by this form is denoted by $B_\varepsilon$. Formally, we have 
\begin{equation}
\label{B_eps}
B_\varepsilon = b(\mathbf{D})^* g^\varepsilon (\mathbf{x})b(\mathbf{D})+\sum _{j=1}^d \left(
a_j^\varepsilon (\mathbf{x})D_j +D_j a_j^\varepsilon (\mathbf{x})^*\right)
+Q^\varepsilon (\mathbf{x}) +\lambda Q_0^\varepsilon (\mathbf{x}).
\end{equation}
Note that, by \eqref{<a_eps<} and \eqref{b_eps <=},
\begin{align}
\label{b_eps<= H1-norm}
&\mathfrak{b}_\varepsilon [\mathbf{u},\mathbf{u}]\leqslant c_5^2\Vert \mathbf{u}\Vert ^2 _{H^1(\mathbb{R}^d)},\quad \mathbf{u}\in H^1(\mathbb{R}^d;\mathbb{C}^n);\\ 
&c_5^2:=\max
\lbrace 
(2+c_1^2+c_2)\alpha _1\Vert g\Vert _{L_\infty};
C(1)+c_3+\vert\lambda\vert \Vert Q_0\Vert _{L_\infty}
\rbrace .
\nonumber
\end{align}
From \eqref{<a_eps<} and \eqref{b_eps >=} it follows that
\begin{equation}
\label{b_eps=> H1-norm}
\begin{split}
&\Vert \mathbf{u}\Vert  _{H^1(\mathbb{R}^d)}\leqslant c_6\Vert B_\varepsilon ^{1/2}\mathbf{u}\Vert _{L_2(\mathbb{R}^d)}
,\quad \mathbf{u}\in H^1(\mathbb{R}^d;\mathbb{C}^n);
\\
&c_6:=(\min\lbrace 2^{-1}\kappa\alpha _0\Vert g^{-1}\Vert ^{-1}_{L_\infty};\beta\rbrace )^{-1/2}.
\end{split}
\end{equation}

Note that  differential expression \eqref{B_eps} also can be treated as a bounded operator $\mathfrak{B}_\varepsilon$ from $H^1(\mathbb{R}^d;\mathbb{C}^n)$ to $H^{-1}(\mathbb{R}^d;\mathbb{C}^n)$. It is easily seen that $\mathfrak{B}_\varepsilon \vert _{\mathrm{Dom}\,B_\varepsilon}=B_\varepsilon$. So, we write simply $B_\varepsilon$ instead of $\mathfrak{B}_\varepsilon$. By \eqref{b_eps<= H1-norm} and the duality arguments, $\Vert B_\varepsilon\Vert _{H^1\rightarrow H^{-1}}\leqslant c_5^2$. 

By the duality arguments, the operator $B_\varepsilon ^{1/2}:H^1\rightarrow L_2$ can be extended to bounded operator from $L_2$ to $H^{-1}$. Indeed, by \eqref{b_eps<= H1-norm}, 
\begin{equation}
\label{B1/2 eps on H1}
\Vert B_\varepsilon^{1/2}\mathbf{u}\Vert_{L_2(\mathbb{R}^d)}\leqslant c_5\Vert\mathbf{u}\Vert_{H^1(\mathbb{R}^d)},\quad \mathbf{u}\in H^1(\mathbb{R}^d;\mathbb{C}^n),
\end{equation}
and
\begin{equation*}
\begin{split}
c_5\Vert\mathbf{u}\Vert _{H^1(\mathbb{R}^d)}\Vert\mathbf{v}\Vert _{L_2(\mathbb{R}^d)}\geqslant\vert (B_\varepsilon ^{1/2}\mathbf{u},\mathbf{v})_{L_2(\mathbb{R}^d)}\vert =\vert (\mathbf{u},B_\varepsilon ^{1/2}\mathbf{v})_{L_2(\mathbb{R}^d)}\vert ,\\
\mathbf{u},\mathbf{v}\in H^1(\mathbb{R}^d;\mathbb{C}^n).
\end{split}
\end{equation*}
So, 
\begin{equation}
\label{B1/2 to H-1 star}
\Vert B_\varepsilon ^{1/2}\mathbf{v}\Vert _{H^{-1}(\mathbb{R}^d)}\leqslant c_5\Vert\mathbf{v}\Vert _{L_2(\mathbb{R}^d)},\quad\mathbf{v}\in H^1(\mathbb{R}^d;\mathbb{C}^n).
\end{equation}
By continuity, we extend $B_\varepsilon ^{1/2}$ to a bounded operator from $L_2$ to $H^{-1}$. By construction, the restriction of this extended operator onto $H^1$ coincides with the original operator $B_\varepsilon^{1/2}$. So, we denote this extended operator by $B_\varepsilon ^{1/2}$. 
By \eqref{B1/2 to H-1 star},
\begin{equation}
\label{B_eps 1/2 to H-1}
\Vert B_\varepsilon ^{1/2}\Vert _{L_2(\mathbb{R}^d)\rightarrow H^{-1}(\mathbb{R}^d)}\leqslant c_5.
\end{equation}

For convenience of further references, in what follows by the ,,problem data'' we mean the following set of parameters:
\begin{equation}
\label{problem data}
\begin{split}
&d,m,n,\rho ;\alpha _0,\alpha _1,\Vert g\Vert _{L_\infty}, \Vert g^{-1}\Vert _{L_\infty}, \Vert a_j\Vert _{L_\rho (\Omega)},j=1,\dots,d;
\\
&\widetilde{c},c_0,\widetilde{c}_2, c_3 \;\mbox{from Condition~\ref{Condition on mu}};
\\
&\lambda, \Vert Q_0\Vert _{L_\infty}, \Vert Q_0^{-1}\Vert _{L_\infty};\;\mbox{the parameters of the lattice }\Gamma.
\end{split}
\end{equation}

\subsection{The effective matrix}  
The effective operator for $A_\varepsilon =b(\mathbf{D})^*g^\varepsilon (\mathbf{x})b(\mathbf{D})$ is given by 
 $A^0=b(\mathbf{D})^*g^0b(\mathbf{D})$, where $g^0$ is a constant $(m\times m)$-matrix called the effective matrix. 
The matrix $g^0$ is defined in terms of the solution of an auxiliary cell problem. 
Suppose that a $\Gamma$-periodic $(n\times m)$-matrix-valued function $\Lambda (\mathbf{x})$ is the (weak) solution of the problem
\begin{equation}
\label{Lambda problem}
b(\mathbf{D})^*g(\mathbf{x})(b(\mathbf{D})\Lambda (\mathbf{x})+\mathbf{1}_m)=0,\quad \int _{\Omega }\Lambda (\mathbf{x})\,d\mathbf{x}=0.
\end{equation}
Then the effective matrix is given by 
\begin{equation}
\label{g^0}
g^0=\vert \Omega \vert ^{-1}\int _{\Omega} \widetilde{g}(\mathbf{x})\,d\mathbf{x},
\end{equation}
where
\begin{equation}
\label{tilde g}
\widetilde{g}(\mathbf{x})=g(\mathbf{x})(b(\mathbf{D})\Lambda (\mathbf{x})+\mathbf{1}_m).
\end{equation}
It can be checked that $g^0$ is positive definite.

We also need the following estimates for the solution of the problem
\eqref{Lambda problem} proved in \cite[(6.28) and Subsection~7.3]{BSu05}:
\begin{align}
\label{Lambda <=}
&\Vert \Lambda \Vert _{L_2(\Omega)}\leqslant \vert \Omega \vert ^{1/2}M_1,\quad M_1=m^{1/2}(2r_0)^{-1}\alpha _0^{-1/2}\Vert g\Vert ^{1/2}_{L_\infty}\Vert g^{-1}\Vert ^{1/2}_{L_\infty},\\
\label{DLambda<=}
&\Vert \mathbf{D}\Lambda \Vert _{L_2(\Omega)}\leqslant \vert \Omega \vert ^{1/2}M_2,\quad M_2=m^{1/2}\alpha _0^{-1/2}\Vert g\Vert ^{1/2}_{L_\infty}\Vert g^{-1}\Vert ^{1/2}_{L_\infty}.
\end{align}

The effective matrix satisfies the estimates known as the Voigt--Reuss bracketing 
(see, e.~g., \cite[Chapter~3, Theorem~1.5]{BSu}).

\begin{proposition}
Let $g^0$ be the effective matrix \eqref{g^0}. Then 
\begin{equation}
\label{Foigt-Reiss}
\underline{g}\leqslant g^0\leqslant \overline{g}.
\end{equation}
If $m=n$, then $g^0=\underline{g}$.
\end{proposition}

Now we distinguish the cases where one of the inequalities in \eqref{Foigt-Reiss}  becomes an identity, see  \cite[Chapter 3, Propositions 1.6 and 1.7]{BSu}. 

\begin{proposition}
The identity $g^0=\overline{g}$ is equivalent to the relations
\begin{equation}
\label{overline-g}
b(\mathbf{D})^* {\mathbf g}_k(\mathbf{x}) =0,\ \ k=1,\dots,m,
\end{equation}
where ${\mathbf g}_k(\mathbf{x})$, $k=1,\dots,m,$ are the columns of the matrix $g(\mathbf{x})$.
\end{proposition}

\begin{proposition} The identity $g^0 =\underline{g}$ is equivalent to the relations
\begin{equation}
\label{underline-g}
{\mathbf l}_k(\mathbf{x}) = {\mathbf l}_k^0 + b(\mathbf{D}) {\mathbf w}_k,\ \ {\mathbf l}_k^0\in \mathbb{C}^m,\ \
{\mathbf w}_k \in {H}^1_{\mathrm{per}}(\Omega;\mathbb{C}^m),\ \ k=1,\dots,m,
\end{equation}
where ${\mathbf l}_k(\mathbf{x})$, $k=1,\dots,m,$ are the columns of the matrix $g(\mathbf{x})^{-1}$.
\end{proposition}

\subsection{The effective operator}
\label{Subsection Effective operator}
In order to define the effective operator for $B_\varepsilon$, consider the 
 $\Gamma$-periodic $(n\times n)$-matrix-valued function $\widetilde{\Lambda}(\mathbf{x})$ which is the solution of the problem
\begin{equation}
\label{tildeLambda_problem}
b(\mathbf{D})^*g(\mathbf{x})b(\mathbf{D})\widetilde{\Lambda }(\mathbf{x})+\sum \limits _{j=1}^dD_ja_j(\mathbf{x})^*=0,\quad \int _{\Omega }\widetilde{\Lambda }(\mathbf{x})\,d\mathbf{x}=0.
\end{equation}
(Here equation is understood in the weak sense.)

The following estimates for $\widetilde{\Lambda}$ were proven in \cite[(7.51) and (7.52)]{SuAA}:
\begin{align}
\label{tilde Lambda<=}
&\Vert \widetilde{\Lambda}\Vert _{L_2(\Omega)}\leqslant  (2r_0)^{-1}C_an^{1/2}\alpha _0^{-1}\Vert g^{-1}\Vert _{L_\infty}=\vert \Omega\vert ^{1/2} \widetilde{M}_1,
\\
\label{D tilde Lambda}
&\Vert \mathbf{D}\widetilde{\Lambda}\Vert _{L_2(\Omega)}\leqslant C_a n^{1/2}\alpha _0^{-1}\Vert g^{-1}\Vert _{L_\infty}=\vert \Omega\vert ^{1/2} \widetilde{M}_2,
\end{align}
where $C_a^2:=\sum _{j=1}^d \int _\Omega \vert a_j(\mathbf{x})\vert ^2\,d\mathbf{x}$, $\widetilde{M}_1:=\vert \Omega\vert ^{-1/2}(2r_0)^{-1} C_a n^{1/2}\alpha _0^{-1}\Vert g^{-1}\Vert _{L_\infty}$, and $\widetilde{M}_2:=\vert \Omega\vert ^{-1/2} C_a n^{1/2}\alpha _0^{-1}\Vert g^{-1}\Vert _{L_\infty}$.

Define constant matrices $V$ and $W$ as follows:
\begin{align*}
&V=\vert \Omega \vert ^{-1}\int _{\Omega}(b(\mathbf{D})\Lambda (\mathbf{x}))^*g(\mathbf{x})(b(\mathbf{D})\widetilde{\Lambda}(\mathbf{x}))\,d\mathbf{x},\\
&W=\vert \Omega \vert ^{-1}\int _{\Omega} (b(\mathbf{D})\widetilde{\Lambda}(\mathbf{x}))^*g(\mathbf{x})(b(\mathbf{D})\widetilde{\Lambda}(\mathbf{x}))\,d\mathbf{x}.
\end{align*}
The effective operator for the operator \eqref{B_eps} is given by 
\begin{equation*}
B^0 =b(\mathbf{D})^*g^0b(\mathbf{D})-b(\mathbf{D})^*V-V^*b(\mathbf{D})
+\sum \limits _{j=1}^d (\overline{a_j+a_j^*})D_j-W+\overline{Q}+\lambda \overline{Q_0}.
\end{equation*}
Here $\overline{Q}=\vert \Omega\vert ^{-1}\int _\Omega d\mu(\mathbf{x})$. 
The operator $B^0$ is the elliptic second order operator with constant coefficients, $\mathrm{Dom}\,B^0=H^2(\mathbb{R}^d;\mathbb{C}^n)$. According to \cite[(3.31), (4.9), Subsection 7.2]{SuAA14}, the operator $B^0$ is positive definite: 
\begin{equation}
\label{B0>=}
B^0\geqslant \check{c} I, \quad \check{c}:=\min\lbrace 2^{-1}\kappa \alpha _0\Vert g^{-1}\Vert ^{-1}_{L_\infty};\beta \rbrace .
\end{equation}
 Moreover, its symbol $L(\boldsymbol{\xi})$ satisfies the estimate 
\begin{equation}
\label{L(xi)>=}
L(\boldsymbol{\xi})\geqslant \check{c}(\vert \boldsymbol{\xi}\vert ^2 +1)\mathbf{1}_n. 
\end{equation}
So,
\begin{equation}
\label{H1-norm <= B0 1/2}
\Vert \mathbf{u}\Vert _{H^1(\mathbb{R}^d)}\leqslant\check{c}^{-1/2}\Vert (B^0)^{1/2}\mathbf{u}\Vert _{L_2(\mathbb{R}^d)},\quad\mathbf{u}\in H^1(\mathbb{R}^d;\mathbb{C}^n),
\end{equation}
and
\begin{align}
\label{B0 -1/2}
&\Vert  (B^0)^{-1/2}\Vert _{L_2(\mathbb{R}^d)\rightarrow L_2(\mathbb{R}^d)}\leqslant \check{c}^{-1/2},
\\
\label{D B0 -1/2}
&\Vert \mathbf{D}(B^0)^{-1/2} \Vert _{L_2(\mathbb{R}^d)\rightarrow L_2(\mathbb{R}^d)}\leqslant \sup _{\boldsymbol{\xi}\in\mathbb{R}^d}\vert \boldsymbol{\xi}\vert \vert L(\boldsymbol{\xi})\vert ^{-1/2}\leqslant \check{c}^{-1/2},
\\
\label{1.41a}
&\Vert (B^0)^{-1/2} \Vert _{L_2(\mathbb{R}^d)\rightarrow H^1(\mathbb{R}^d)}\leqslant \sup _{\boldsymbol{\xi}\in\mathbb{R}^d}(1+\vert \boldsymbol{\xi}\vert ^2)^{1/2}\vert \vert L(\boldsymbol{\xi})\vert ^{-1/2}\leqslant \check{c}^{-1/2},
\\
\label{B0^2 L2 to H2}
&\Vert (B^0)^{-1}\Vert _{L_2(\mathbb{R}^d)\rightarrow H^2(\mathbb{R}^d)}\leqslant \sup _{\boldsymbol{\xi}\in\mathbb{R}^d}
(1+\vert \boldsymbol{\xi}\vert ^2)\vert L(\boldsymbol{\xi})\vert ^{-1}
\leqslant \check{c}^{-1}.
\end{align}

For the case, when the measure $d\mu$ is as in Example~\ref{Example of mu}, in \cite[Lemma 1.6]{MSuPOMI} it was proven that
\begin{equation}
\label{L(xi)<=}
L(\boldsymbol{\xi})\leqslant C_L(\vert \boldsymbol{\xi}\vert ^2+1)\mathbf{1}_n
\end{equation}
for some positive constant $C_L$ depending only on the problem data \eqref{problem data}. To modify the proof from \cite{MSuPOMI} to our case, we need only to estimate the term $\overline{Q}$. It turns out that $\overline{Q}\leqslant c_3\mathbf{1}_n$. Indeed, 
by Condition~\ref{Condition on mu}($1^\circ$), 
\begin{equation*}
\begin{split}
&\vert (\overline{Q}\mathbf{u},\mathbf{u})_{L_2(\mathbb{R}^d)}\vert
=
\vert \Omega\vert ^{-1}\Bigl\vert \int _{\mathbb{R}^d}\langle \int _\Omega d\mu (\mathbf{x})\mathbf{u}(\mathbf{y}),\mathbf{u}(\mathbf{y})\rangle \,d\mathbf{y}\Bigr\vert 
\\
&\leqslant \vert \Omega\vert ^{-1}\sum _{\mathbf{a}\in\Gamma}\int _{\Omega _\mathbf{a}}\Bigl\vert \langle\int _\Omega d\mu (\mathbf{x})\mathbf{u}(\mathbf{y}),\mathbf{u}(\mathbf{y})\rangle \Bigr\vert \,d\mathbf{y}
\leqslant
c_3\Vert \mathbf{u}\Vert ^2 _{L_2(\mathbb{R}^d)},\quad\mathbf{u}\in H^1(\mathbb{R}^d;\mathbb{C}^n).
\end{split}
\end{equation*}

Now, from \eqref{L(xi)<=} it follows that
\begin{equation}
\label{B01/2 from L2 to H-1}
\Vert (B^0)^{1/2}\Vert _{L_2(\mathbb{R}^d)\rightarrow H^{-1}(\mathbb{R}^d)}\leqslant\sup _{\boldsymbol{\xi}\in\mathbb{R}^d}(1+\vert \boldsymbol{\xi}\vert ^2)^{-1/2}\vert L(\boldsymbol{\xi})\vert ^{1/2}\leqslant C_L^{1/2}
\end{equation}
and
\begin{align}
\label{B01/2 on H1}
&\Vert (B^0)^{1/2}\Vert _{H^1(\mathbb{R}^d)\rightarrow L_2(\mathbb{R}^d)}\leqslant C_L^{1/2},
\\
\label{B0 on H2}
&\Vert B^0 \Vert _{H^2(\mathbb{R}^d) \rightarrow L_2(\mathbb{R}^d)}\leqslant C_L,
\\
\label{B0 3/2 on H3}
&\Vert (B^0)^{3/2}\Vert _{H^3(\mathbb{R}^d) \rightarrow  L_2(\mathbb{R}^d)}\leqslant C_L^{3/2}.
\end{align}

\subsection{The operator $\Pi _\varepsilon$}
\label{Subsection Pi_eps}

Let $\Pi _\varepsilon$ be the pseudodifferential operator in $L_2(\mathbb{R}^d;\mathbb{C}^n)$ whose symbol is the characteristic function $\chi _{{\Omega}^*/\varepsilon}(\boldsymbol{\xi})$ of the set ${\Omega}^*/\varepsilon$:
\begin{equation}
\label{Pi eps = definition}
(\Pi _\varepsilon \mathbf{f})(\mathbf{x})=(2\pi)^{-d/2}\int _{{\Omega}^*/\varepsilon}e^{i\langle \mathbf{x},\boldsymbol{\xi}\rangle}\widehat{\mathbf{f}} (\boldsymbol{\xi})\,d\boldsymbol{\xi}.
\end{equation}
Here $\widehat{\mathbf{f}}$ is the Fourier-image of the function $\mathbf{f}$. 
Obviously, $\Pi_\varepsilon \mathbf{D}^\alpha \mathbf{u}=\mathbf{D}^\alpha \Pi _\varepsilon \mathbf{u}$ for $\mathbf{u}\in H^l(\mathbb{R}^d;\mathbb{C}^n)$, if $\vert \alpha\vert \leqslant l$. 
Note, that $\Pi _\varepsilon ^2=\Pi _\varepsilon$ and
\begin{equation}
\label{Pi eps <= 1}
\Vert \Pi _\varepsilon \Vert _{L_2(\mathbb{R}^d)\rightarrow L_2(\mathbb{R}^d)}= 1.
\end{equation}

The following property of the operator $\Pi _\varepsilon$ was proven in  \cite[Proposition 1.4]{PSu}.

\begin{proposition}
\label{Proposition Pi -I}
We have 
\begin{equation*}
\Vert \Pi _\varepsilon \mathbf{u}-\mathbf{u}\Vert _{L_2(\mathbb{R}^d)}
\leqslant \varepsilon r_0^{-1}\Vert \mathbf{D}\mathbf{u}\Vert _{L_2(\mathbb{R}^d)}, \quad\mathbf{u}\in H^1(\mathbb{R}^d;\mathbb{C}^n),\quad \varepsilon >0.
\end{equation*}
\end{proposition}

For $\varkappa =0$, the following result was obtained in \cite[Subsection~10.2]{BSu06}. For $0<\varkappa<\infty$ the proof is quite similar.

\begin{proposition}
\label{Proposition f Pi on H-kappa}
Let $f$ be a $\Gamma$-periodic function in $\mathbb{R}^d$ such that $f\in L_2(\Omega)$. Let $0\leqslant \varkappa <\infty$. Then 
\begin{equation*}
\Vert [f^\varepsilon]\Pi _\varepsilon \Vert _{H^{-\varkappa}(\mathbb{R}^d)\rightarrow H^{-\varkappa}(\mathbb{R}^d)}\leqslant\vert \Omega\vert ^{-1/2}\Vert f\Vert _{L_2(\Omega)},\quad \varepsilon >0.
\end{equation*}
\end{proposition}

\begin{proof}
We have
\begin{equation*}
\Vert [f^\varepsilon]\Pi _\varepsilon \Vert _{H^{-\varkappa}(\mathbb{R}^d)\rightarrow H^{-\varkappa}(\mathbb{R}^d)}=\Vert (\mathbf{D}^2+I)^{-\varkappa }[f^\varepsilon ] \Pi _\varepsilon (\mathbf{D}^2+I)^{\varkappa}\Vert _{L_2(\mathbb{R}^d)\rightarrow L_2(\mathbb{R}^d)}.
\end{equation*}
Let $T_\varepsilon$ be the unitary rescaling operator in $L_2(\mathbb{R}^d;\mathbb{C}^n)$ given by $(T_\varepsilon \mathbf{u})(\mathbf{x})=\varepsilon ^{d/2} \mathbf{u}(\mathbf{x})$. Then
\begin{equation*}
(\mathbf{D}^2+I)^{-\varkappa}[ f^\varepsilon ]\Pi _\varepsilon (\mathbf{D}^2+I)^{\varkappa}
=T_\varepsilon ^* (\mathbf{D}^2+\varepsilon ^2 I)^{-\varkappa}[ f]\Pi _1 (\mathbf{D}^2+\varepsilon ^2 I)^{\varkappa}T_\varepsilon.
\end{equation*}
So,
\begin{equation*}
\Vert [f^\varepsilon]\Pi _\varepsilon \Vert _{H^{-\varkappa}(\mathbb{R}^d)\rightarrow H^{-\varkappa}(\mathbb{R}^d)}=\Vert  (\mathbf{D}^2+\varepsilon ^2 I)^{-\varkappa}[ f]\Pi _1 (\mathbf{D}^2+\varepsilon ^2 I)^{\varkappa}\Vert _{L_2(\mathbb{R}^d)\rightarrow L_2(\mathbb{R}^d)}.
\end{equation*}

It turns out that under Gelfand transformation (see Subsection~\ref{Subsection lattices}) the operator $$(\mathbf{D}^2+\varepsilon ^2 I)^{-\varkappa}[ f]\Pi _1 (\mathbf{D}^2+\varepsilon ^2 I)^{\varkappa}=
(\mathbf{D}^2+\varepsilon ^2 I)^{-\varkappa}[ f](\mathbf{D}^2+\varepsilon ^2 I)^{\varkappa}\Pi _1 $$ is decomposed into the direct integral of
the operators $$((\mathbf{D}+\mathbf{k})^2+\varepsilon ^2 I)^{-\varkappa}[ f] ((\mathbf{D}+\mathbf{k})^2+\varepsilon ^2 I)^{\varkappa}P,$$ depending on the quasi-momentum $\mathbf{k}\in{\Omega}^*$ and acting in $L_2(\Omega;\mathbb{C}^n)$ with the periodic boundary conditions. Here $P$ is the operator of averaging over the cell: $P\mathbf{v}=\vert \Omega\vert ^{-1}\int _\Omega \mathbf{v}(\mathbf{x})\,d\mathbf{x}$,  $\mathbf{v}\in L_2(\Omega;\mathbb{C}^n)$, i.~e., the projection of $L_2(\Omega ;\mathbb{C}^n)$ onto the subspace of constant functions.  
Because of the presence of the projection $P$, 
$$((\mathbf{D}+\mathbf{k})^2+\varepsilon ^2 I)^{-\varkappa}[ f] ((\mathbf{D}+\mathbf{k})^2+\varepsilon ^2 I)^{\varkappa}P
=((\mathbf{D}+\mathbf{k})^2+\varepsilon ^2 I)^{-\varkappa}[ f] (\mathbf{k}^2+\varepsilon ^2  )^{\varkappa} P.
$$
Thus,
\begin{equation}
\label{end of proof of prop 1.7}
\begin{split}
&\Vert [f^\varepsilon]\Pi _\varepsilon \Vert _{H^{-\varkappa}(\mathbb{R}^d)\rightarrow H^{-\varkappa}(\mathbb{R}^d)}=\Vert  (\mathbf{D}^2+\varepsilon ^2 I)^{-\varkappa}[ f]\Pi _1 (\mathbf{D}^2+\varepsilon ^2 I)^{\varkappa}\Vert _{L_2(\mathbb{R}^d)\rightarrow L_2(\mathbb{R}^d)}
\\
&\leqslant \sup _{\mathbf{k}\in {\Omega}^*} \Vert ((\mathbf{D}+\mathbf{k})^2+\varepsilon ^2 I)^{-\varkappa}[ f]P (\mathbf{k}^2+\varepsilon ^2 I)^{\varkappa}\Vert _{L_2(\Omega)\rightarrow L_2(\Omega)}
\\
&\leqslant  \Vert [ f]P \Vert _{L_2(\Omega)\rightarrow L_2(\Omega)} \sup _{\mathbf{k}\in {\Omega}^*} (\mathbf{k}^2+\varepsilon ^2 )^{\varkappa}\Vert 
((\mathbf{D}+\mathbf{k})^2+\varepsilon ^2 )^{-\varkappa}\Vert _{L_2(\Omega)\rightarrow L_2(\Omega)} .
\end{split}
\end{equation}
Using the discrete Fourier transform and definition of the first Brillouin zone ${\Omega}^*$, we obtain
\begin{equation*}
\Vert 
((\mathbf{D}+\mathbf{k})^2+\varepsilon ^2 )^{-\varkappa}\Vert _{L_2(\Omega)\rightarrow L_2(\Omega)} 
=\sup _{\mathbf{b}\in {\Gamma}^*}(\vert \mathbf{b}+\mathbf{k}\vert ^2+\varepsilon ^2)^{-\varkappa}=(\vert \mathbf{k}\vert ^2+\varepsilon ^2)^{-\varkappa}, \quad \mathbf{k}\in {\Omega}^*.
\end{equation*}
Together with \eqref{end of proof of prop 1.7}, this implies
\begin{equation*}
\Vert [f^\varepsilon]\Pi _\varepsilon \Vert _{H^{-\varkappa}(\mathbb{R}^d)\rightarrow H^{-\varkappa}(\mathbb{R}^d)}
\leqslant \Vert [ f]P \Vert _{L_2(\Omega)\rightarrow L_2(\Omega)}
\leqslant \vert \Omega\vert ^{-1/2}\Vert f\Vert _{L_2(\Omega)}.
\end{equation*}
\end{proof}

\subsection{Approximation of the resolvent}

In the present Subsection we formulate the homogenization results for $B_\varepsilon ^{-1}$, obtained in \cite[Theorems 9.2 and 9.7]{SuAA}.

Let $K(\varepsilon)$ be the corrector
\begin{equation}
\label{K(eps)}
K(\varepsilon):=[\Lambda^\varepsilon] \Pi _\varepsilon b(\mathbf{D})(B^0)^{-1}+[\widetilde{\Lambda}^\varepsilon ]\Pi _\varepsilon (B^0)^{-1}.
\end{equation}
The continuity of the operator $K(\varepsilon)$ from $L_2(\mathbb{R}^d;\mathbb{C}^n)$ to $H^1(\mathbb{R}^d;\mathbb{C}^n)$ follows from Proposition~\ref{Proposition f Pi on H-kappa} with $\varkappa =0$ and inclusions $\Lambda,\widetilde{\Lambda}\in {H}^1_{\mathrm{per}}(\Omega)$. By Proposition~\ref{Proposition f Pi on H-kappa} and \eqref{<b^*b<}, \eqref{Lambda <=}, \eqref{tilde Lambda<=}, and \eqref{B0^2 L2 to H2},
\begin{equation}
\label{K(eps)<= L2 to L2}
\Vert K(\varepsilon )\Vert _{L_2(\mathbb{R}^d)\rightarrow L_2(\mathbb{R}^d)}\leqslant C_K:=(M_1\alpha _1^{1/2}+\widetilde{M}_1)\check{c}^{-1}.
\end{equation}

\begin{theorem}[\textnormal{\cite{SuAA}}]
\label{Theorem elliptic}
Let the assumptions of Subsections~\textnormal{\ref{Subsection lattices}--\ref{Subsection Effective operator}} be satisfied. Let $K(\varepsilon)$ be the corrector \eqref{K(eps)}. Denote
\begin{align}
\label{R1}
&\mathcal{R}_1(\varepsilon):= B_{\varepsilon}^{-1}-(B^0)^{-1},
\\
\label{R2}
&\mathcal{R}_2(\varepsilon):=B_{\varepsilon}^{-1}-(B^0)^{-1}-\varepsilon K(\varepsilon).
\end{align} 
Then for $0<\varepsilon\leqslant 1$ we have
\begin{align}
&
\Vert \mathcal{R}_1(\varepsilon)\Vert _{L_2(\mathbb{R}^d)\rightarrow L_2(\mathbb{R}^d)}
\leqslant C_1\varepsilon,
\nonumber
\\
\label{B-eps with corrector}
& 
\Vert  \mathcal{R}_2(\varepsilon)\Vert _{L_2(\mathbb{R}^d)\rightarrow H^1(\mathbb{R}^d)}
\leqslant C_2\varepsilon .
\end{align}
The constants $C_1$ and $C_2$ are controlled in terms of the problem data \eqref{problem data}.
\end{theorem}

Note that inequalities \eqref{b_eps<= H1-norm} and \eqref{B-eps with corrector} imply the estimate
\begin{align}
\label{Beps1/2 appr corr}
\Vert B_\varepsilon^{1/2}\mathcal{R}_2(\varepsilon)\Vert _{L_2(\mathbb{R}^d)\rightarrow L_2(\mathbb{R}^d)}
\leqslant c_5 C_2\varepsilon .
\end{align}

\section{Problem setting. Main results}

\label{Section Problem setting. Main results}

\subsection{Problem}

Let $\mathbf{u}_\varepsilon$ be the generalized solution of the following Cauchy problem for hyperbolic system:
\begin{equation}
\label{u_eps problem}
\begin{cases}
\partial _t^2 \mathbf{u}_\varepsilon (\mathbf{x},t) =-B_{\varepsilon}\mathbf{u}_\varepsilon (\mathbf{x},t) ,\quad \mathbf{x}\in\mathbb{R}^d,\,t\in\mathbb{R},
\\
\mathbf{u}_\varepsilon (\mathbf{x},0)=\boldsymbol{\phi}(\mathbf{x}),\quad (\partial _t\mathbf{u}_\varepsilon)(\mathbf{x},0)=\boldsymbol{\psi}(\mathbf{x}),\quad \mathbf{x}\in\mathbb{R}^d.
\end{cases}
\end{equation}
Here $\boldsymbol{\phi}\in H^1(\mathbb{R}^d;\mathbb{C}^n)$ and  $\boldsymbol{\psi}\in L_2(\mathbb{R}^d;\mathbb{C}^n)$. 

The solution of problem \eqref{u_eps problem} is given by 
\begin{equation*}
\mathbf{u}_\varepsilon (\cdot ,t)=\cos (tB_{\varepsilon}^{1/2})\boldsymbol{\phi}+B_{\varepsilon}^{-1/2}\sin(tB_{\varepsilon}^{1/2})\boldsymbol{\psi}.
\end{equation*}
Thus,
\begin{equation*}
\partial _t\mathbf{u}_\varepsilon (\cdot ,t)=-B_{\varepsilon}^{1/2}\sin(tB_{\varepsilon}^{1/2})\boldsymbol{\phi}
+\cos (t B_{\varepsilon}^{1/2})\boldsymbol{\psi}.
\end{equation*}

\textit{Our goal} is to study the behavior of the solution $\mathbf{u}_\varepsilon (\cdot ,t)$ as $\varepsilon \rightarrow 0$. In other words, to approximate the~operator functions $\cos (tB_{\varepsilon}^{1/2})$ and 
$B_{\varepsilon}^{-1/2}\sin(tB_{\varepsilon}^{1/2})$ in suitable norms.

The problem \eqref{u_eps problem} can be rewritten as follows
\begin{equation*}
\partial _t \begin{pmatrix}
  \mathbf{u}_\varepsilon\\
  \partial _t\mathbf{u}_\varepsilon
\end{pmatrix}
=\begin{pmatrix}
  0& I\\
 -B_{\varepsilon}& 0
\end{pmatrix}
\begin{pmatrix}
  \mathbf{u}_\varepsilon\\
  \partial _t\mathbf{u}_\varepsilon
\end{pmatrix},
\quad
\begin{pmatrix}
  \mathbf{u}_\varepsilon (\cdot ,0)\\
  \partial _t\mathbf{u}_\varepsilon(\cdot ,0)
\end{pmatrix}
=
\begin{pmatrix}
  \boldsymbol{\phi}\\
  \boldsymbol{\psi}
\end{pmatrix}.
\end{equation*}
Denote
\begin{equation}
\label{A_eps}
\begin{split}
\mathfrak{A}_\varepsilon :=\begin{pmatrix}
  0& I\\
 -B_{\varepsilon}& 0
\end{pmatrix}
: H^1(\mathbb{R}^d;\mathbb{C}^n)\times L_2(\mathbb{R}^d;\mathbb{C}^n)\rightarrow H^1(\mathbb{R}^d;\mathbb{C}^n)\times L_2(\mathbb{R}^d;\mathbb{C}^n)
,\\
 \mathrm{Dom}\,\mathfrak{A}_\varepsilon :=\mathrm{Dom}\,B_{\varepsilon}\times H^1(\mathbb{R}^d;\mathbb{C}^n).
\end{split}
\end{equation}
(Our choice of $ \mathrm{Dom}\,\mathfrak{A}_\varepsilon$ guaranties that 
$\mathrm{Ran}\,\mathfrak{A}_\varepsilon \subset H^1(\mathbb{R}^d;\mathbb{C}^n)\times L_2(\mathbb{R}^d;\mathbb{C}^n)$.) 
According to \cite[Chapter 2, Section 7]{Go}, operator \eqref{A_eps} generates a $C_0$-group 
\begin{equation}
\label{exp =}
\begin{split}
e^{t\mathfrak{A}_\varepsilon}=
\begin{pmatrix}
  \cos(tB_{\varepsilon}^{1/2})& B_{\varepsilon}^{-1/2}\sin (tB_{\varepsilon}^{1/2})\\
  -B_{\varepsilon}^{1/2}\sin (tB_{\varepsilon}^{1/2})& \cos(tB_{\varepsilon}^{1/2})
\end{pmatrix}
\end{split}
\end{equation}
on the space $\mathrm{Dom}\,\mathfrak{b}_\varepsilon \times L_2(\mathbb{R}^d;\mathbb{C}^n)=H^1(\mathbb{R}^d;\mathbb{C}^n)\times L_2(\mathbb{R}^d;\mathbb{C}^n)$ equipped with the graph norm of $B_\varepsilon^{1/2}$. By \eqref{b_eps<= H1-norm} and \eqref{b_eps=> H1-norm}, this norm is equivalent to the standard norm in $H^1(\mathbb{R}^d;\mathbb{C}^n)\times L_2(\mathbb{R}^d;\mathbb{C}^n)$. So,
\begin{equation}
\label{exp B eps bounded}
e^{t\mathfrak{A}_\varepsilon}\in\mathcal{B}(H^1(\mathbb{R}^d;\mathbb{C}^n)\times L_2(\mathbb{R}^d;\mathbb{C}^n)).
\end{equation}

Then
\begin{equation*}
\begin{pmatrix}
  \mathbf{u}_\varepsilon\\
  \partial _t\mathbf{u}_\varepsilon
\end{pmatrix}
=e^{t\mathfrak{A}_\varepsilon}
\begin{pmatrix}
  \boldsymbol{\phi}\\
  \boldsymbol{\psi}
\end{pmatrix}.
\end{equation*}

It is easily seen that the resolvent of the operator $\mathfrak{A}_\varepsilon$ has the form
\begin{equation}
\label{A_eps^-1}
\mathfrak{A}_\varepsilon ^{-1}=
\begin{pmatrix}
  0& -B_{\varepsilon}^{-1}\\
  I& 0
\end{pmatrix},
\quad \mathfrak{A}_\varepsilon ^{-1}\in\mathcal{B}(H^1(\mathbb{R}^d;\mathbb{C}^n)\times L_2(\mathbb{R}^d;\mathbb{C}^n)).
\end{equation}

\subsection{The operator $\mathfrak{A}_0$}

Let $\mathfrak{A}_0$ be the effective operator for \eqref{A_eps}:
\begin{equation}
\label{A0}
\begin{split}
&\mathfrak{A}_0=\begin{pmatrix}
  0& I\\
  -B^0& 0
\end{pmatrix}: H^1(\mathbb{R}^d;\mathbb{C}^n)\times L_2(\mathbb{R}^d;\mathbb{C}^n)\rightarrow H^1(\mathbb{R}^d;\mathbb{C}^n)\times L_2(\mathbb{R}^d;\mathbb{C}^n),
\\
&\mathrm{Dom}\,\mathfrak{A}_0 = \mathrm{Dom}\,B^0\times H^1(\mathbb{R}^d;\mathbb{C}^n)=H^2(\mathbb{R}^d;\mathbb{C}^n)\times H^1(\mathbb{R}^d;\mathbb{C}^n).
\end{split}
\end{equation}
Then, similarly to \eqref{exp =}--\eqref{A_eps^-1},
\begin{align}
\label{exp 0 =}
\begin{split}
&e^{t\mathfrak{A}_0}=
\begin{pmatrix}
  \cos(t(B^0)^{1/2})& (B^0)^{-1/2}\sin (t(B^0)^{1/2})\\
  -(B^0)^{1/2}\sin (t(B^0)^{1/2})& \cos(t(B^0)^{1/2})
\end{pmatrix},
\end{split}
\\ 
\label{exp B0 bounded}
&e^{t\mathfrak{A}_0}\in\mathcal{B}(H^1(\mathbb{R}^d;\mathbb{C}^n)\times L_2(\mathbb{R}^d;\mathbb{C}^n)),
\end{align}
and
\begin{equation}
\label{A0 -1}
\mathfrak{A}_0^{-1}=
\begin{pmatrix}
  0& -(B^0)^{-1}\\
  I& 0
\end{pmatrix},\quad \mathfrak{A}_0^{-1}\in \mathcal{B}(H^1(\mathbb{R}^d;\mathbb{C}^n)\times L_2(\mathbb{R}^d;\mathbb{C}^n)).
\end{equation}
Thus
\begin{equation}
\label{A_0^-2}
\mathfrak{A}_0^{-2}=
\begin{pmatrix}
  0& -(B^0)^{-1}\\
  I& 0
\end{pmatrix}
\begin{pmatrix}
  0& -(B^0)^{-1}\\
  I& 0
\end{pmatrix}
=
\begin{pmatrix}
  -(B^0)^{-1}& 0\\
  0& -(B^0)^{-1}
\end{pmatrix}
.
\end{equation}

\subsection{Principal term of approximation}

The principal term of approximation is given by the following theorem. The proof can be found in Subsection~\ref{Subsection proof of Th  principal part} below.

\begin{theorem}
\label{Theorem principal part}
Under the assumptions of Subsections~\textnormal{\ref{Subsection lattices}--\ref{Subsection Effective operator}}, 
for $t\in\mathbb{R}$ and $0<\varepsilon\leqslant 1$ we have
\begin{align}
\label{Th cos L2}
\Vert &\cos(t B_\varepsilon ^{1/2})-\cos(t (B^0)^{1/2})\Vert _{H^2(\mathbb{R}^d)\rightarrow L_2(\mathbb{R}^d)}
\leqslant C_3\varepsilon(1+\vert t\vert),
\\
\label{Th sin L2}
\begin{split}
\Vert &B_\varepsilon ^{-1/2}\sin(t B_\varepsilon ^{1/2})-(B^0)^{-1/2}\sin(t(B^0)^{1/2})\Vert _{H^1(\mathbb{R}^d)\rightarrow L_2(\mathbb{R}^d)}
\leqslant C_4\varepsilon (1+\vert t\vert),
\end{split}
\\
\label{Th sin to H-1 from L2}
\begin{split}
\Vert &B_\varepsilon ^{1/2}\sin (tB_\varepsilon ^{1/2})-(B^0)^{1/2}\sin (t(B^0)^{1/2})\Vert _{H^2(\mathbb{R}^d)\rightarrow H^{-1}(\mathbb{R}^d)}
\leqslant C_5\varepsilon (1+\vert t\vert),
\end{split}
\\
\label{Th cos to H-1 from H1}
\begin{split}
\Vert &\cos(t B_\varepsilon ^{1/2})-\cos(t (B^0)^{1/2})\Vert _{H^1(\mathbb{R}^d)\rightarrow H^{-1}(\mathbb{R}^d)}
\leqslant C_6\varepsilon \vert t\vert .
\end{split}
\end{align}
The constants $C_3$, $C_4$, $C_5$, and $C_6$ are controlled explicitly in terms of the problem data \eqref{problem data}.
\end{theorem}

\begin{remark}
The right-hand side in \eqref{Th sin L2} can be replaced by $C\varepsilon\vert t\vert$, see \eqref{Th sin no corr new} below.
\end{remark}

\begin{corollary}
For $t\in\mathbb{R}$ and $0<\varepsilon\leqslant 1$ we have
\begin{equation*}
\Vert e^{t\mathfrak{A}_\varepsilon}-e^{t\mathfrak{A}_0}\Vert _{H^2(\mathbb{R}^d)\times H^1(\mathbb{R}^d)\rightarrow L_2(\mathbb{R}^d)\times H^{-1}(\mathbb{R}^d)}\leqslant (C_3+C_4+C_5+C_6)\varepsilon (1+\vert t\vert).
\end{equation*}
\end{corollary}

By the interpolation arguments, Theorem~\ref{Theorem principal part} implies the following result.

\begin{theorem}
Let $0\leqslant r\leqslant 2$. Then, under the assumptions of Theorem~\textnormal{\ref{Theorem principal part}}, for $0<\varepsilon\leqslant 1 $ and $t\in\mathbb{R}$ we have
\begin{align}
\label{Th cos L2 interpolation}
\Vert &\cos(t B_\varepsilon ^{1/2})-\cos(t (B^0)^{1/2})\Vert _{H^r(\mathbb{R}^d)\rightarrow L_2(\mathbb{R}^d)}
\leqslant C_7\varepsilon ^{r/2}(1+\vert t\vert)^{r/2},
\\
\label{Th sin L2 interpolation}
\begin{split}
\Vert &B_\varepsilon ^{-1/2}\sin(t B_\varepsilon ^{1/2})-(B^0)^{-1/2}\sin(t(B^0)^{1/2})\Vert _{H^{r-1}(\mathbb{R}^d)\rightarrow L_2(\mathbb{R}^d)}
\\
&\leqslant C_8\varepsilon ^{r/2} (1+\vert t\vert)^{r/2},
\end{split}
\\
\label{Th dt sin H-1 interpolation}
\Vert &B_\varepsilon ^{1/2}\sin (tB_\varepsilon ^{1/2})-(B^0)^{1/2}\sin (t(B^0)^{1/2})\Vert _{H^r(\mathbb{R}^d)\rightarrow H^{-1}(\mathbb{R}^d)}
\leqslant C_9\varepsilon ^{r/2} (1+\vert t\vert)^{r/2},
\\
\label{Th dt cos H-1 interpolation}
\Vert &\cos(t B_\varepsilon ^{1/2})-\cos(t (B^0)^{1/2})\Vert _{H^{r-1}(\mathbb{R}^d)\rightarrow H^{-1}(\mathbb{R}^d)}
\leqslant C_{10}\varepsilon ^{r/2} \vert t\vert ^{r/2} .
\end{align}
The constants $C_7$, $C_8$, $C_9$, and $C_{10}$ depend only on $r$ and the problem data \eqref{problem data}.
\end{theorem}

\begin{remark}
The right-hand side in \eqref{Th sin L2 interpolation} can be replaced by $C\varepsilon ^{r/2}\vert t\vert ^{r/2}$, see \eqref{Th corr int dt} below.
\end{remark}

\begin{proof}
Interpolating between the rough estimate
\begin{equation*}
\Vert \cos(t B_\varepsilon ^{1/2})-\cos(t (B^0)^{1/2})\Vert _{ L_2(\mathbb{R}^d)\rightarrow L_2(\mathbb{R}^d)}
\leqslant 2
\end{equation*}
and \eqref{Th cos L2}, we obtain \eqref{Th cos L2 interpolation} with the constant $C_7=2^{1-r/2} C_3^{r/2}$.

Next, by the duality arguments and \eqref{b_eps=> H1-norm},
\begin{equation}
\label{sin eps from H-1}
\begin{split}
\Vert B_\varepsilon ^{-1/2}\sin (tB_\varepsilon ^{1/2})\Vert _{H^{-1}(\mathbb{R}^d)\rightarrow L_2(\mathbb{R}^d)}&=\Vert B_\varepsilon ^{-1/2}\sin (tB_\varepsilon ^{1/2})\Vert _{L_2(\mathbb{R}^d)\rightarrow H^1(\mathbb{R}^d)}
\\
&\leqslant c_6\Vert \sin (tB_\varepsilon ^{1/2})\Vert _{L_2(\mathbb{R}^d)\rightarrow L_2(\mathbb{R}^d)}=c_6.
\end{split}
\end{equation}
Similarly, according to \eqref{H1-norm <= B0 1/2},
\begin{equation*}
\Vert (B^0)^{-1/2}\sin (t (B^0)^{1/2})\Vert _{H^{-1}(\mathbb{R}^d)\rightarrow L_2(\mathbb{R}^d)}\leqslant \check{c}^{-1/2}.
\end{equation*}
Together with \eqref{sin eps from H-1}, this implies the inequality
\begin{equation}
\label{difference of sine from H-1}
\Vert B_\varepsilon ^{-1/2}\sin (tB_\varepsilon ^{1/2})-(B^0)^{-1/2}\sin (t (B^0)^{1/2})\Vert _{H^{-1}(\mathbb{R}^d)\rightarrow L_2(\mathbb{R}^d)}\leqslant c_6 +\check{c}^{-1/2}.
\end{equation}
Interpolating between \eqref{Th sin L2} and \eqref{difference of sine from H-1}, we arrive at the estimate \eqref{Th sin L2 interpolation} with the constant $C_8:=(c_6 +\check{c}^{-1/2})^{1-r/2}C_4^{r/2}$.

From \eqref{B_eps 1/2 to H-1} and \eqref{B01/2 from L2 to H-1} it follows that
\begin{equation}
\label{dt sin from L2 to H-1}
\Vert B_\varepsilon ^{1/2}\sin (tB_\varepsilon ^{1/2})-(B^0)^{1/2}\sin (t(B^0)^{1/2})\Vert _{L_2(\mathbb{R}^d)\rightarrow H^{-1}(\mathbb{R}^d)}
\leqslant c_5+C_L^{1/2}.
\end{equation}
Interpolating between \eqref{dt sin from L2 to H-1} and \eqref{Th sin to H-1 from L2}, we obtain estimate \eqref{Th dt sin H-1 interpolation}, where $C_9:=(c_5+C_L^{1/2})^{1-r/2}C_5^{r/2}$.

Note that the operator $\cos (tB_\varepsilon ^{1/2})$ can be extended to the bounded operator in $H^{-1}(\mathbb{R}^d;\mathbb{C}^n)$. Indeed, by \eqref{b_eps=> H1-norm} and \eqref{B1/2 eps on H1}, for $\mathbf{u}\in L_2(\mathbb{R}^d;\mathbb{C}^n)$ we have
\begin{equation*}
\begin{split}
\Vert \cos (tB_\varepsilon ^{1/2})\mathbf{u}\Vert _{H^{-1}(\mathbb{R}^d)}
&=\sup _{0\neq \mathbf{v}\in H^1(\mathbb{R}^d)}
\frac{\vert (\cos (t B_\varepsilon ^{1/2})\mathbf{u},\mathbf{v})_{L_2(\mathbb{R}^d)}\vert }{ \Vert \mathbf{v}\Vert _{H^1(\mathbb{R}^d)}}
\\
&=\sup _{0\neq \mathbf{v}\in H^1(\mathbb{R}^d)}
\frac{\vert (\mathbf{u},\cos (t B_\varepsilon ^{1/2})\mathbf{v})_{L_2(\mathbb{R}^d)}\vert }{ \Vert \mathbf{v}\Vert _{H^1(\mathbb{R}^d)}}
\\
&\leqslant\Vert \mathbf{u}\Vert _{H^{-1}(\mathbb{R}^d)}\sup _{0\neq \mathbf{v}\in H^1(\mathbb{R}^d)}\frac{\Vert \cos (t B_\varepsilon ^{1/2})\mathbf{v}\Vert _{H^1(\mathbb{R}^d)} }{ \Vert \mathbf{v}\Vert _{H^1(\mathbb{R}^d)}}
\\
&\leqslant c_6 \Vert \mathbf{u}\Vert _{H^{-1}(\mathbb{R}^d)}\sup _{0\neq \mathbf{v}\in H^1(\mathbb{R}^d)}\frac{\Vert  B_\varepsilon ^{1/2}\mathbf{v}\Vert _{L_2(\mathbb{R}^d)} }{ \Vert \mathbf{v}\Vert _{H^1(\mathbb{R}^d)}}
\\
&\leqslant c_5c_6 \Vert \mathbf{u}\Vert _{H^{-1}(\mathbb{R}^d)}.
\end{split}
\end{equation*}
So, by continuity we can extend the operator $\cos (tB_\varepsilon ^{1/2})$ onto $H^{-1}(\mathbb{R}^d;\mathbb{C}^n)$, and 
\begin{equation}
\label{cos eps on H-1}
\Vert \cos (tB_\varepsilon ^{1/2})\Vert _{H^{-1}(\mathbb{R}^d)\rightarrow H^{-1}(\mathbb{R}^d)}\leqslant c_5c_6.
\end{equation}
The operator $\cos (t(B^0)^{1/2})$ also can be extended onto $H^{-1}(\mathbb{R}^d;\mathbb{C}^n)$. Since the operator $B^0$ has constant coefficients, it commutes with differentiation, and
\begin{equation}
\label{cos 0 on H-1}
\begin{split}
&\Vert \cos (t(B^0)^{1/2})\Vert  _{H^{-1}(\mathbb{R}^d)\rightarrow H^{-1}(\mathbb{R}^d)}
\\
&=
\Vert (\mathbf{D}^2 +I)^{-1/2}\cos (t(B^0)^{1/2})(\mathbf{D}^2 +I)^{1/2}\Vert _{L_2(\mathbb{R}^d)\rightarrow L_2(\mathbb{R}^d)}\leqslant 1.
\end{split}
\end{equation}
Combining \eqref{cos eps on H-1} and \eqref{cos 0 on H-1}, we obtain
\begin{equation}
\label{differences of cosines on H-1}
\Vert \cos (tB_\varepsilon ^{1/2})-\cos (t(B^0)^{1/2})\Vert  _{H^{-1}(\mathbb{R}^d)\rightarrow H^{-1}(\mathbb{R}^d)}
\leqslant 1+c_5c_6.
\end{equation}
Interpolating between \eqref{differences of cosines on H-1} and \eqref{Th cos to H-1 from H1}, we arrive at estimate \eqref{Th dt cos H-1 interpolation} with the constant $C_{10}:=C_6^{r/2}(1+c_5c_6)^{1-r/2}$.
\end{proof}

\subsection{Proof of Theorem~\textnormal{\ref{Theorem principal part}}}
\label{Subsection proof of Th  principal part}

\begin{proof}[Proof of Theorem~\textnormal{\ref{Theorem principal part}}]
Denote
\begin{equation*}
\mathcal{E}(t):=e^{-t\mathfrak{A}_\varepsilon}\mathfrak{A}_\varepsilon ^{-1}\mathfrak{A}_0^{-1}e^{t\mathfrak{A}_0}.
\end{equation*}
It is bounded operator in $ H^1(\mathbb{R}^d;\mathbb{C}^n)\times L_2 (\mathbb{R}^d;\mathbb{C}^n)$, because all factors are bounded in $ H^1(\mathbb{R}^d;\mathbb{C}^n)\times L_2 (\mathbb{R}^d;\mathbb{C}^n)$. 
Then
\begin{equation*}
\frac{d\mathcal{E}(t)}{d t}=e^{- t\mathfrak{A}_\varepsilon}(\mathfrak{A}_\varepsilon ^{-1}-\mathfrak{A}_0^{-1})e^{t\mathfrak{A}_0}.
\end{equation*}
(Derivative is taken in the strong operator topology in $H^1(\mathbb{R}^d;\mathbb{C}^n)\times L_2(\mathbb{R}^d;\mathbb{C}^n)$, see \cite[Chapter 1, Theorem 2.4(c)]{Pa}.)
So,
\begin{equation}
\label{tozd epxp start}
\mathcal{E}(t)-\mathcal{E}(0)=e^{-t\mathfrak{A}_\varepsilon}\mathfrak{A}_\varepsilon ^{-1}\mathfrak{A}_0^{-1}e^{t\mathfrak{A}_0}-\mathfrak{A}_\varepsilon ^{-1}\mathfrak{A}_0^{-1}
=\int _0^t e^{-s\mathfrak{A}_\varepsilon}(\mathfrak{A}_\varepsilon ^{-1}-\mathfrak{A}_0^{-1})e^{s\mathfrak{A}_0}\,ds.
\end{equation}
(The integral is understood as a strong Riemann integral in $H^1(\mathbb{R}^d;\mathbb{C}^n)\times L_2(\mathbb{R}^d;\mathbb{C}^n)$.)
 
Let us multiply the identity \eqref{tozd epxp start} by $ e^{-t\mathfrak{A}_0}$ from the right. Then 
\begin{equation*}
e^{-t\mathfrak{A}_\varepsilon}\mathfrak{A}_\varepsilon ^{-1}\mathfrak{A}_0^{-1}-\mathfrak{A}_\varepsilon ^{-1}\mathfrak{A}_0^{-1}e^{-t\mathfrak{A}_0}
=\int _0^t e^{-s\mathfrak{A}_\varepsilon}(\mathfrak{A}_\varepsilon ^{-1}-\mathfrak{A}_0^{-1})e^{(s-t)\mathfrak{A}_0}\,ds.
\end{equation*}
So,
\begin{equation}
\label{tozd exp}
\begin{split}
(e^{-t\mathfrak{A}_\varepsilon}-e^{-t\mathfrak{A}_0})\mathfrak{A}_0^{-2}
&=-e^{-t\mathfrak{A}_\varepsilon}(\mathfrak{A}_\varepsilon ^{-1}-\mathfrak{A}_0^{-1})\mathfrak{A}_0^{-1}
+(\mathfrak{A}_\varepsilon ^{-1}-\mathfrak{A}_0^{-1})\mathfrak{A}_0^{-1}e^{-t\mathfrak{A}_0}
\\
&+\int _0^t e^{-s\mathfrak{A}_\varepsilon}(\mathfrak{A}_\varepsilon ^{-1}-\mathfrak{A}_0^{-1})e^{(s-t)\mathfrak{A}_0}\,ds.
\end{split}
\end{equation}

Denote
\begin{equation}
\label{G(eps)=}
\mathfrak{G}(\varepsilon):=\begin{pmatrix}
\Lambda ^\varepsilon \Pi_\varepsilon b(\mathbf{D})+\widetilde{\Lambda}^\varepsilon \Pi _\varepsilon &0\\
0&0
\end{pmatrix},
\quad
\mathrm{Dom}\,\mathfrak{G}(\varepsilon)= H^2(\mathbb{R}^d;\mathbb{C}^n)\times L_2(\mathbb{R}^d;\mathbb{C}^n).
\end{equation}
Then, by \eqref{A0 -1},
\begin{equation}
\label{G(eps)A0-1}
\mathfrak{G}(\varepsilon)\mathfrak{A}_0^{-1}
=\begin{pmatrix}
\Lambda ^\varepsilon \Pi_\varepsilon b(\mathbf{D})+\widetilde{\Lambda}^\varepsilon \Pi _\varepsilon &0\\
0&0
\end{pmatrix}
\begin{pmatrix}
0&-(B^0)^{-1}\\
I&0
\end{pmatrix}
=\begin{pmatrix}
0&-K(\varepsilon)\\
0&0
\end{pmatrix}.
\end{equation}
Here $K(\varepsilon)$ is the operator \eqref{K(eps)}. So, 
\begin{equation}
\label{G(eps)res bounded}
\mathfrak{G}(\varepsilon)\mathfrak{A}_0^{-1}\in\mathcal{B}(H^1(\mathbb{R}^d;\mathbb{C}^n)\times L_2(\mathbb{R}^d;\mathbb{C}^n)).
\end{equation}

Note that 
\begin{equation}
\label{I}
\mathfrak{A}_\varepsilon ^{-1}-\mathfrak{A}_0^{-1}=
\begin{pmatrix}
0&-\mathcal{R}_1(\varepsilon)\\
0&0
\end{pmatrix}
\end{equation}
and
\begin{equation}
\label{II}
\mathfrak{A}_\varepsilon ^{-1}-\mathfrak{A}_0^{-1}-\varepsilon \mathfrak{G}(\varepsilon)\mathfrak{A}_0^{-1}=
\begin{pmatrix}
0&-\mathcal{R}_2(\varepsilon)\\
0&0
\end{pmatrix},
\end{equation}
where $\mathcal{R}_1(\varepsilon)$ and $\mathcal{R}_2(\varepsilon)$ are defined by \eqref{R1} and \eqref{R2}, respectively.

By \eqref{tozd exp},
\begin{equation}
\label{tozd exp-1}
\begin{split}
(e^{-t\mathfrak{A}_\varepsilon}-e^{-t\mathfrak{A}_0})\mathfrak{A}_0^{-2}
&=-e^{-t\mathfrak{A}_\varepsilon}(\mathfrak{A}_\varepsilon ^{-1}-\mathfrak{A}_0^{-1})\mathfrak{A}_0^{-1}
+(\mathfrak{A}_\varepsilon ^{-1}-\mathfrak{A}_0^{-1})\mathfrak{A}_0^{-1}e^{-t\mathfrak{A}_0}
\\
&+\int _0^t e^{-s\mathfrak{A}_\varepsilon}(\mathfrak{A}_\varepsilon ^{-1}-\mathfrak{A}_0^{-1}-\varepsilon \mathfrak{G}(\varepsilon)\mathfrak{A}_0^{-1})e^{(s-t)\mathfrak{A}_0}\,ds
\\
&+\int _0^t e^{-s\mathfrak{A}_\varepsilon} \varepsilon \mathfrak{G}(\varepsilon)\mathfrak{A}_0^{-1}e^{(s-t)\mathfrak{A}_0}\,ds.
\end{split}
\end{equation}
Denote the consecutive summands in the right-hand side of \eqref{tozd exp-1} by $\mathcal{I}_j(\varepsilon ;t)$, $j=1,\dots ,4$.

Then, by  \eqref{exp =},  \eqref{A0 -1}, and \eqref{I},
\begin{equation}
\label{I1=}
\mathcal{I}_1(\varepsilon ;t)=
\begin{pmatrix}
\cos(t B_\varepsilon ^{1/2})\mathcal{R}_1(\varepsilon )&0\\
B_\varepsilon ^{1/2}\sin (tB_\varepsilon ^{1/2})\mathcal{R}_1(\varepsilon )&0
\end{pmatrix}.
\end{equation}
Next, the operator $\mathcal{I}_2(\varepsilon ;t)$ has the entries
\begin{align}
\label{I2 11}
&\lbrace \mathcal{I}_2(\varepsilon ;t)\rbrace _{11}=-\mathcal{R}_1(\varepsilon )\cos (t(B^0)^{1/2}),
\\
\label{I2 12}
&\lbrace \mathcal{I}_2(\varepsilon ;t)\rbrace _{12}=
\mathcal{R}_1(\varepsilon )(B^0)^{-1/2}\sin (t(B^0)^{1/2}),
\\
\label{I2 21=22}
&\lbrace \mathcal{I}_2(\varepsilon ;t)\rbrace _{21}=\lbrace \mathcal{I}_2(\varepsilon ;t)\rbrace _{22}=0.
\end{align}

By using \eqref{exp =}, \eqref{exp 0 =}, and \eqref{II}, we see that  
the entries of the operator $\mathcal{I}_3(\varepsilon ;t)$ have the form
 \begin{align}
 \label{I-3 11}
 \begin{split}
 &\lbrace \mathcal{I}_3(\varepsilon ;t)\rbrace _{11}
= \int _0 ^t \cos (s B_\varepsilon ^{1/2})\mathcal{R}_2(\varepsilon)(B^0)^{1/2}\sin ((s-t)(B^0)^{1/2})\,ds,
\end{split}
\\
\label{I-3 12}
 &\lbrace \mathcal{I}_3(\varepsilon ;t)\rbrace _{12}
= -\int _0 ^t \cos (sB_\varepsilon ^{1/2})\mathcal{R}_2(\varepsilon)\cos ((s-t)(B^0)^{1/2})\,ds,
\\
\label{I-3 21 a}
\begin{split}
 &\lbrace \mathcal{I}_3(\varepsilon ;t)\rbrace _{21}
 = \int _0 ^t 
 B_\varepsilon ^{1/2}\sin (sB_\varepsilon ^{1/2})
 \mathcal{R}_2(\varepsilon)(B^0)^{1/2}\sin ((s-t)(B^0)^{1/2})\,ds,
 \end{split}
 \\
  \label{I3 22}
 \begin{split}
 &\lbrace \mathcal{I}_3(\varepsilon ;t)\rbrace _{22}
 = -\int _0^t  B_\varepsilon ^{1/2}\sin (sB_\varepsilon ^{1/2})
 \mathcal{R}_2(\varepsilon)\cos((s-t)(B^0)^{1/2})\,ds.
 \end{split}
 \end{align}
 Here integrals are understood in the strong sense (in different norms): in \eqref{I-3 11} in $(H^1\rightarrow H^1)$-norm, in \eqref{I-3 12} in $(L_2\rightarrow H^1)$-norm, etc. 
We can understand the integral from \eqref{I-3 11} in the strong $(H^1\rightarrow L_2)$-sense and integrate by parts. Then
 \begin{align}
 \label{I3 11}
 \begin{split}
 \lbrace \mathcal{I}_3(\varepsilon ;t)\rbrace _{11}&=
 -\cos (tB_\varepsilon ^{1/2})\mathcal{R}_2(\varepsilon)
 +\mathcal{R}_2(\varepsilon)\cos(t(B^0)^{1/2})
 \\
 &-\int _0^t B_\varepsilon ^{1/2}\sin (sB_\varepsilon ^{1/2})
 \mathcal{R}_2(\varepsilon)\cos ((s-t)(B^0)^{1/2})\,ds
 \end{split}
\end{align} 
is bounded operator in $L_2(\mathbb{R}^d;\mathbb{C}^n)$. Next, 
the integral in \eqref{I-3 21 a} is understood in the strong $(H^1\rightarrow L_2)$-sense. So, it makes sense in a weaker $(H^1\rightarrow H^{-1})$-topology. Integrating by parts, we arrive at the expression
 \begin{align}
   \label{I3 21}
 \begin{split}
  \lbrace \mathcal{I}_3(\varepsilon ;t)\rbrace _{21}&=
  -B_\varepsilon ^{1/2}\sin(tB_\varepsilon ^{1/2})\mathcal{R}_2(\varepsilon)
  \\
  &+\int _0^t B_\varepsilon ^{1/2} \cos(sB_\varepsilon ^{1/2})B_\varepsilon ^{1/2}\mathcal{R}_2(\varepsilon )\cos((s-t)(B^0)^{1/2})\,ds.
 \end{split}
 \end{align}
 This operator is obviously $(L_2\rightarrow H^{-1})$-bounded.

According to \eqref{exp =}, \eqref{exp 0 =},  and \eqref{G(eps)A0-1}, the operator $\mathcal{I}_4(\varepsilon ;t)$ has the entries
\begin{align}
\label{I4 11}
&\lbrace \mathcal{I}_4(\varepsilon ;t)\rbrace _{11}
= \int _0 ^t \cos (sB_\varepsilon ^{1/2})\varepsilon K(\varepsilon)(B^0)^{1/2}\sin ((s-t)(B^0)^{1/2})\,ds,
\\
\label{I4 12}
&\lbrace \mathcal{I}_4(\varepsilon ;t)\rbrace _{12}
= -\int _0 ^t \cos (sB_\varepsilon ^{1/2})\varepsilon K(\varepsilon )\cos((s-t)(B^0)^{1/2})\,ds,
\\
\label{I4 21}
&\lbrace \mathcal{I}_4(\varepsilon ;t)\rbrace _{21}
= \int _0 ^t B_\varepsilon ^{1/2}\sin(sB_\varepsilon ^{1/2})\varepsilon K(\varepsilon)(B^0)^{1/2}\sin ((s-t)(B^0)^{1/2})\,ds,
\\
\label{I4 22}
&\lbrace \mathcal{I}_4(\varepsilon ;t)\rbrace _{22}
= -\int _0 ^t B_\varepsilon ^{1/2}\sin (sB_\varepsilon ^{1/2})\varepsilon K(\varepsilon )\cos((s-t)(B^0)^{1/2})\,ds.
\end{align}

Now, we compute the matrix entries of the left-hand side of \eqref{tozd exp-1}. 
By \eqref{exp =}, \eqref{exp 0 =}, and \eqref{A_0^-2},
\begin{align}
\label{left in tozd 11}
&\lbrace (e^{-t\mathfrak{A}_\varepsilon}-e^{-t\mathfrak{A}_0})\mathfrak{A}_0^{-2}\rbrace _{11}=-(\cos(tB_\varepsilon ^{1/2})-\cos(t(B^0)^{1/2}))(B^0)^{-1},
\\
\label{left in tozd 12}
&\lbrace (e^{-t\mathfrak{A}_\varepsilon}-e^{-t\mathfrak{A}_0})\mathfrak{A}_0^{-2}\rbrace _{12}=(B_\varepsilon ^{-1/2}\sin (tB_\varepsilon ^{1/2})-(B^0)^{-1/2}\sin(t(B^0)^{1/2}))(B^0)^{-1},
\\
\label{left in tozd 21}
&\lbrace (e^{-t\mathfrak{A}_\varepsilon}-e^{-t\mathfrak{A}_0})\mathfrak{A}_0^{-2}\rbrace _{21}
=-(B_\varepsilon ^{1/2}\sin (tB_\varepsilon ^{1/2})-(B^0)^{1/2}\sin (t(B^0)^{1/2}))(B^0)^{-1},
\\
\label{left in tozd 22}
&\lbrace (e^{-t\mathfrak{A}_\varepsilon}-e^{-t\mathfrak{A}_0})\mathfrak{A}_0^{-2}\rbrace _{22}
=-(\cos(tB_\varepsilon ^{1/2})-\cos (t(B^0)^{1/2}))(B^0)^{-1}.
\end{align}

Combining \eqref{R1}, \eqref{R2}, \eqref{tozd exp-1}--\eqref{I2 11}, \eqref{I3 11}, \eqref{I4 11}, and \eqref{left in tozd 11}, we obtain
\begin{equation}
\label{cos identity}
\begin{split}
&-(\cos(tB_\varepsilon ^{1/2})-\cos(t(B^0)^{1/2}))(B^0)^{-1}
=
\cos(t B_\varepsilon ^{1/2})\mathcal{R}_1(\varepsilon)
\\
&-\mathcal{R}_1(\varepsilon )\cos (t(B^0)^{1/2})
-\cos (tB_\varepsilon ^{1/2})\mathcal{R}_2(\varepsilon)
 +\mathcal{R}_2(\varepsilon)\cos(t(B^0)^{1/2})
 \\
 &-\int _0^t B_\varepsilon ^{1/2}\sin (sB_\varepsilon ^{1/2})
 \mathcal{R}_2(\varepsilon)\cos ((s-t)(B^0)^{1/2})\,ds
 \\
 &+\int _0 ^t \cos (sB_\varepsilon ^{1/2})\varepsilon K(\varepsilon)(B^0)^{1/2}\sin ((s-t)(B^0)^{1/2})\,ds
 \\
 &=
\cos(t B_\varepsilon ^{1/2}) \varepsilon K(\varepsilon)-\varepsilon K(\varepsilon)\cos(t(B^0)^{1/2})\\
 &-\int _0^t B_\varepsilon ^{1/2}\sin (sB_\varepsilon ^{1/2})
 \mathcal{R}_2(\varepsilon)\cos ((s-t)(B^0)^{1/2})\,ds
 \\
 &+\int _0 ^t \cos (sB_\varepsilon ^{1/2})\varepsilon K(\varepsilon)(B^0)^{1/2}\sin ((s-t)(B^0)^{1/2})\,ds
 .
\end{split}
\end{equation}
By construction, the matrix entries with indices ,,$11$'' are bounded operators in $H^1(\mathbb{R}^d;\mathbb{C}^n)$. After integrating by parts in the term $\lbrace\mathcal{I}_3(\varepsilon ;t)\rbrace _{11}$, we can understand \eqref{cos identity} only as equality for the operators acting from $H^1(\mathbb{R}^d;\mathbb{C}^n)$ to $L_2(\mathbb{R}^d;\mathbb{C}^n)$. By continuity, we extend equality from $H^1(\mathbb{R}^d;\mathbb{C}^n)$ to the whole space $L_2(\mathbb{R}^d;\mathbb{C}^n)$.

By Proposition~\ref{Proposition f Pi on H-kappa} with $\varkappa =0$, \eqref{<b^*b<}, \eqref{Lambda <=}, and \eqref{tilde Lambda<=}, 
\begin{equation}
\label{partof corr H1 to L2}
\Vert [\Lambda ^\varepsilon]\Pi _\varepsilon b(\mathbf{D})+[\widetilde{\Lambda}^\varepsilon ]\Pi _\varepsilon\Vert _{H^1(\mathbb{R}^d)\rightarrow L_2(\mathbb{R}^d)}\leqslant M_1\alpha _1^{1/2}+\widetilde{M}_1.
\end{equation}
Together with \eqref{1.41a} and \eqref{K(eps)}, this implies
\begin{equation}
\label{K(eps)B01/2}
\Vert K(\varepsilon )(B^0)^{1/2}\Vert _{L_2(\mathbb{R}^d)\rightarrow L_2(\mathbb{R}^d)}\leqslant (M_1\alpha _1^{1/2}+\widetilde{M}_1)\check{c}^{-1/2}.
\end{equation}

Combining \eqref{K(eps)<= L2 to L2}, \eqref{Beps1/2 appr corr}, \eqref{cos identity}, and \eqref{K(eps)B01/2}, we obtain
\begin{equation}
\label{cos L2 almost final}
\begin{split}
\Vert& (\cos(tB_\varepsilon ^{1/2})-\cos(t(B^0)^{1/2}))(B^0)^{-1}\Vert _{L_2(\mathbb{R}^d)\rightarrow L_2(\mathbb{R}^d)}
\\
&\leqslant 2C_K\varepsilon +c_5C_2\varepsilon \vert t\vert +\check{c}^{-1/2}(\alpha _1^{1/2}M_1+\widetilde{M}_1)\varepsilon\vert t\vert .
\end{split}
\end{equation}

Now from \eqref{B0 on H2} and \eqref{cos L2 almost final} we derive estimate \eqref{Th cos L2} with the constant
$$C_3:=C_L\max\lbrace 2C_K;c_5C_2+\check{c}^{-1/2}(\alpha _1^{1/2}M_1+\widetilde{M}_1)\rbrace .$$

We proceed to the proof of estimate \eqref{Th sin L2}. Combining \eqref{tozd exp-1}, \eqref{I1=}, \eqref{I2 12}, \eqref{I-3 12}, \eqref{I4 12}, and \eqref{left in tozd 12}, we have
\begin{equation*}
\begin{split}
&(B_\varepsilon ^{-1/2}\sin (tB_\varepsilon ^{1/2})-(B^0)^{-1/2}\sin(t(B^0)^{1/2}))(B^0)^{-1}
\\
&=\mathcal{R}_1(\varepsilon)(B^0)^{-1/2}\sin (t(B^0)^{1/2})
\\
&-\int _0 ^t \cos (sB_\varepsilon ^{1/2})\mathcal{R}_2(\varepsilon)\cos ((s-t)(B^0)^{1/2})\,ds
\\
&-\int _0 ^t \cos (sB_\varepsilon ^{1/2})\varepsilon K(\varepsilon )\cos((s-t)(B^0)^{1/2})\,ds.
\end{split}
\end{equation*}
This equality can be understood as an identity for the operators acting from $L_2$ to $L_2$. So, it makes sense on the range of the operator $(B^0)^{1/2}:H^1\rightarrow L_2$:
\begin{equation*}
\begin{split}
&(B_\varepsilon ^{-1/2}\sin (tB_\varepsilon ^{1/2})-(B^0)^{-1/2}\sin(t(B^0)^{1/2}))(B^0)^{-1/2}
\\
&=\mathcal{R}_1(\varepsilon )\sin (t(B^0)^{1/2})
\\
&-\int _0 ^t \cos (sB_\varepsilon ^{1/2})\mathcal{R}_2(\varepsilon)(B^0)^{1/2}\cos ((s-t)(B^0)^{1/2})\,ds
\\
&-\int _0 ^t \cos (sB_\varepsilon ^{1/2})\varepsilon K(\varepsilon )(B^0)^{1/2}\cos((s-t)(B^0)^{1/2})\,ds.
\end{split}
\end{equation*}
The equality here is understood in the $(H^1\rightarrow L_2)$-sense. By continuity, we extend it onto the whole space $L_2$. Integrating by parts in the first integral, we get
\begin{equation*}
\begin{split}
&(B_\varepsilon ^{-1/2}\sin (tB_\varepsilon ^{1/2})-(B^0)^{-1/2}\sin(t(B^0)^{1/2}))(B^0)^{-1/2}
\\
&=\mathcal{R}_1(\varepsilon)\sin (t(B^0)^{1/2})
-\mathcal{R}_2(\varepsilon)\sin(t(B^0) ^{1/2})
\\
&-\int _0^t\sin(sB_\varepsilon ^{1/2})B_\varepsilon ^{1/2}\mathcal{R}_2(\varepsilon)\sin ((s-t)(B^0)^{1/2})\,ds
\\
&-\int _0 ^t \cos (sB_\varepsilon ^{1/2})\varepsilon K(\varepsilon )(B^0)^{1/2}\cos((s-t)(B^0)^{1/2})\,ds.
\end{split}
\end{equation*}
Combining this with \eqref{K(eps)<= L2 to L2}--\eqref{R2},  \eqref{Beps1/2 appr corr}, and \eqref{K(eps)B01/2}, we obtain
\begin{equation*}
\begin{split}
\Vert &(B_\varepsilon ^{-1/2}\sin (tB_\varepsilon ^{1/2})-(B^0)^{-1/2}\sin(t(B^0)^{1/2}))(B^0)^{-1/2}\Vert _{L_2(\mathbb{R}^d)\rightarrow L_2(\mathbb{R}^d)}
\\
&\leqslant \varepsilon C_K+c_5C_2\varepsilon\vert t\vert +(M_1\alpha _1^{1/2}+\widetilde{M}_1)\check{c}^{-1/2}\varepsilon\vert t\vert .
\end{split}
\end{equation*}
By \eqref{B01/2 on H1}, this implies \eqref{Th sin L2} with the constant
$$C_4:=C_L^{1/2}\max\lbrace C_K;c_5C_2+\check{c}^{-1/2}(M_1\alpha _1^{1/2}+\widetilde{M}_1)\rbrace.$$

Bringing together \eqref{tozd exp-1}, \eqref{I1=}, \eqref{I2 21=22}, \eqref{I3 21}, \eqref{I4 21}, and \eqref{left in tozd 21},  we have
\begin{equation}
\label{equality for 12 principal}
\begin{split}
&-(B_\varepsilon ^{1/2}\sin (tB_\varepsilon ^{1/2})-(B^0)^{1/2}\sin (t(B^0)^{1/2}))(B^0)^{-1}
\\
&=B_\varepsilon ^{1/2}\sin (tB_\varepsilon ^{1/2})\mathcal{R}_1(\varepsilon)
-B_\varepsilon ^{1/2}\sin(tB_\varepsilon ^{1/2})\mathcal{R}_2(\varepsilon)
  \\
  &+\int _0^t B_\varepsilon ^{1/2} \cos(sB_\varepsilon ^{1/2})B_\varepsilon ^{1/2}\mathcal{R}_2(\varepsilon )\cos((s-t)(B^0)^{1/2})\,ds
  \\
  &+\int _0 ^t B_\varepsilon ^{1/2}\sin(sB_\varepsilon ^{1/2})\varepsilon K(\varepsilon)(B^0)^{1/2}\sin ((s-t)(B^0)^{1/2})\,ds.
\end{split}
\end{equation}
By construction, this equality should be understood in the $(H^1\rightarrow H^{-1})$-sense. By continuity, it make sense in the $(L_2\rightarrow H^{-1})$-topology. 
Together with \eqref{B_eps 1/2 to H-1}, \eqref{B0 on H2},  \eqref{K(eps)<= L2 to L2}--\eqref{R2}, \eqref{Beps1/2 appr corr}, and \eqref{K(eps)B01/2}, equality \eqref{equality for 12 principal} implies estimate \eqref{Th sin to H-1 from L2} with the constant 
$$C_5:=C_L\max\lbrace c_5C_K;c_5^2C_2+c_5(M_1\alpha _1^{1/2}+\widetilde{M}_1)\check{c}^{-1/2}\rbrace .$$

Combining \eqref{tozd exp-1}, \eqref{I1=}, \eqref{I2 21=22}, \eqref{I3 22}, \eqref{I4 22}, and \eqref{left in tozd 22}, we obtain 
\begin{equation*}
\begin{split}
&(\cos(tB_\varepsilon ^{1/2})-\cos (t(B^0)^{1/2}))(B^0)^{-1}
\\
&=\int _0^t  B_\varepsilon ^{1/2}\sin (sB_\varepsilon ^{1/2})
 \mathcal{R}_2(\varepsilon)\cos((s-t)(B^0)^{1/2})\,ds
 \\
 &+\int _0 ^t B_\varepsilon ^{1/2}\sin (sB_\varepsilon ^{1/2})\varepsilon K(\varepsilon )\cos((s-t)(B^0)^{1/2})\,ds.
\end{split}
\end{equation*}
This is the equality for the operators acting from $L_2$ to $L_2$. But $L_2$ is the subset of $H^{-1}$ and the $L_2$-norm is stronger than the $H^{-1}$-norm, so this identity also can be understood in the $(L_2\rightarrow H^{-1})$-sense. Applying it to functions from the space $(B^0)^{1/2}H^1(\mathbb{R}^d;\mathbb{C}^n)=L_2(\mathbb{R}^d;\mathbb{C}^n)$, we arrive at
\begin{equation*}
\begin{split}
&(\cos(tB_\varepsilon ^{1/2})-\cos (t(B^0)^{1/2}))(B^0)^{-1/2}
\\
&=\int _0^t  B_\varepsilon ^{1/2}\sin (sB_\varepsilon ^{1/2})
 \mathcal{R}_2(\varepsilon)(B^0)^{1/2}\cos((s-t)(B^0)^{1/2})\,ds
 \\
 &+\int _0 ^t B_\varepsilon ^{1/2}\sin (sB_\varepsilon ^{1/2})\varepsilon K(\varepsilon )(B^0)^{1/2}\cos((s-t)(B^0)^{1/2})\,ds.
\end{split}
\end{equation*}
Here the equality is understood in the $(H^1\rightarrow H^{-1})$-topology. By continuity, we can understand it in the $(L_2\rightarrow H^{-1})$-sense. Integrating by parts in the first integral, we obtain
\begin{equation*}
\begin{split}
&(\cos(tB_\varepsilon ^{1/2})-\cos (t(B^0)^{1/2}))(B^0)^{-1/2}
\\
&=-\int _0^t  B_\varepsilon ^{1/2} \cos (sB_\varepsilon ^{1/2})
B_\varepsilon ^{1/2} \mathcal{R}_2(\varepsilon)\sin((st)(B^0)^{1/2})\,ds
 \\
 &+\int _0 ^t B_\varepsilon ^{1/2}\sin (sB_\varepsilon ^{1/2})\varepsilon K(\varepsilon )(B^0)^{1/2}\cos((s-t)(B^0)^{1/2})\,ds.
\end{split}
\end{equation*}
Combining this with \eqref{B_eps 1/2 to H-1}, \eqref{B01/2 on H1}, \eqref{Beps1/2 appr corr}, and \eqref{K(eps)B01/2}, we arrive at estimate \eqref{Th cos to H-1 from H1} with the constant $C_6:=C_L^{1/2}(c_5^2C_2+c_5(M_1\alpha_1^{1/2}+\widetilde{M}_1)\check{c}^{-1/2})$.
\end{proof}

\subsection{Approximation with corrector}

\begin{theorem}
\label{Theorem corrector}
Let the assumptions of Subsections~\textnormal{\ref{Subsection lattices}--\ref{Subsection Effective operator}} be satisfied. 
Denote
\begin{align}
\label{K_1(t,eps)}
&\mathcal{K}_1(\varepsilon ;t):=
(\Lambda ^\varepsilon \Pi _\varepsilon b(\mathbf{D})+ \widetilde{\Lambda}^\varepsilon\Pi _\varepsilon)\cos(t (B^0)^{1/2})(B^0)^{-1},
\\
\label{K(t,eps)}
&\mathcal{K}_2(\varepsilon ;t):=
(\Lambda ^\varepsilon \Pi _\varepsilon b(\mathbf{D})+ \widetilde{\Lambda}^\varepsilon\Pi _\varepsilon)(B^0)^{-1/2}\sin(t (B^0)^{1/2}).
\end{align}
Then for $0<\varepsilon\leqslant 1$ and $t\in\mathbb{R}$ we have
\begin{align}
\label{Th cos corr}
\begin{split}
\Vert &\cos (tB_\varepsilon ^{1/2})B_\varepsilon ^{-1}-\cos(t (B^0)^{1/2})(B^0)^{-1}-\varepsilon \mathcal{K}_1(\varepsilon ;t)\Vert _{H^1(\mathbb{R}^d)\rightarrow H^1(\mathbb{R}^d)}
\\
&\leqslant C_{11} \varepsilon (1+\vert t\vert),
\end{split}
\\
\label{Th sin corr}
\begin{split}
\Vert &B_\varepsilon ^{-1/2}\sin(t B_\varepsilon ^{1/2})-(B^0)^{-1/2}\sin(t (B^0)^{1/2})-\varepsilon \mathcal{K}_2(\varepsilon ;t)\Vert _{H^2(\mathbb{R}^d)\rightarrow H^1(\mathbb{R}^d)}
\\
&\leqslant C_{12}\varepsilon (1+\vert t\vert),
\end{split}
\\
\label{Th sin no corr new}
\begin{split}
\Vert &B_\varepsilon ^{-1/2}\sin(t B_\varepsilon^{1/2})-(B^0)^{-1/2}\sin (t (B^0)^{1/2})\Vert _{H^1(\mathbb{R}^d)\rightarrow L_2(\mathbb{R}^d)}
\leqslant C_{13} \varepsilon \vert t\vert .
\end{split}
\end{align}
The constants $C_{11}$, $C_{12}$, and $C_{13}$ are controlled in terms of the problem data \eqref{problem data}.
\end{theorem}

The proof of Theorem~\ref{Theorem corrector} can be found in Subsection~\ref{Subsection proof Th sin corr} below. 
According to \eqref{exp =}, \eqref{exp 0 =}, \eqref{Th cos L2}, and \eqref{G(eps)=}, Theorem~\ref{Theorem corrector} implies the following result.

\begin{corollary}
Let $\mathfrak{A}_\varepsilon$, $\mathfrak{A}_0$, and $\mathfrak{G}(\varepsilon)$ be the operators \eqref{A_eps}, \eqref{A0}, and \eqref{G(eps)=}, respectively. Then for $0<\varepsilon\leqslant 1$ and $t\in\mathbb{R}$ we have
\begin{equation*}
\begin{split}
&\left\Vert e^{t\mathfrak{A}_\varepsilon}\begin{pmatrix}
B_\varepsilon ^{-1}&0\\0&I
\end{pmatrix}
-(I+\varepsilon\mathfrak{G}(\varepsilon))
e^{t\mathfrak{A}_0}\begin{pmatrix}
(B^0)^{-1}&0\\0&I
\end{pmatrix}
\right\Vert _{H^1(\mathbb{R}^d)\times H^2(\mathbb{R}^d)\rightarrow H^1(\mathbb{R}^d)\times L_2(\mathbb{R}^d)}
\\
&\leqslant (C_3+C_{11}+C_{12}+C_{13})\varepsilon(1+\vert t\vert).
\end{split}
\end{equation*}
\end{corollary}

By interpolation arguments, from Theorem~\ref{Theorem corrector} we derive the following result.

\begin{theorem}
\label{Theorem corrector interpolation}
Let $0\leqslant q\leqslant 1$ and let $0\leqslant r\leqslant 2$. Then for $0<\varepsilon\leqslant 1$ and $t\in\mathbb{R}$ we have
\begin{align}
\begin{split}
\label{Th cos corr interpol}
\Vert &\cos (t B_\varepsilon ^{1/2})B_\varepsilon ^{-1} -\cos (t (B^0)^{1/2})(B^0)^{-1}-\varepsilon\mathcal{K}_1(\varepsilon ;t)\Vert _{H^q(\mathbb{R}^d)\rightarrow H^1(\mathbb{R}^d)}
\\
&\leqslant C_{14}\varepsilon ^q (1+\vert t\vert)^q,
\end{split}
\\
\label{Th sin corr interpol}
\begin{split}
\Vert &B_\varepsilon ^{-1/2}\sin(t B_\varepsilon ^{1/2})-(B^0)^{-1/2}\sin(t (B^0)^{1/2})-\varepsilon \mathcal{K}_2(\varepsilon ;t)
\Vert _{H^{1+q}(\mathbb{R}^d)\rightarrow H^1(\mathbb{R}^d)}
\\
&\leqslant C_{15}\varepsilon ^q (1+\vert t\vert)^q,
\end{split}
\\
\label{Th corr int dt}
\Vert &B_\varepsilon ^{-1/2}\sin(t B_\varepsilon^{1/2})-(B^0)^{-1/2}\sin (t (B^0)^{1/2})\Vert _{H^{-1+r}(\mathbb{R}^d)\rightarrow L_2(\mathbb{R}^d)}
\leqslant C_{16} \varepsilon ^{r/2} \vert t\vert ^{r/2}.
\end{align}
The constants $C_{14}$ and $C_{15}$ depend only on $q$ and the problem data \eqref{problem data}, the constant $C_{16}$ depends on $r$ and the~problem data \eqref{problem data}.
\end{theorem}

\begin{proof}
According to \eqref{b_eps >=}, $B_\varepsilon \geqslant \beta I$. So, 
\begin{equation}
\label{B eps -1/2 est}
\Vert B_\varepsilon ^{-1/2}\Vert _{L_2(\mathbb{R}^d)\rightarrow L_2(\mathbb{R}^d)}\leqslant\beta ^{-1/2}.
\end{equation}
Using this argument and \eqref{b_eps=> H1-norm}, we obtain
\begin{equation}
\label{cos eps L2 to H1}
\begin{split}
\Vert \cos (t B_\varepsilon ^{1/2})B_\varepsilon ^{-1}\Vert _{L_2(\mathbb{R}^d)\rightarrow H^1(\mathbb{R}^d)}
\leqslant
c_6\Vert B_\varepsilon ^{-1/2}\Vert  _{L_2(\mathbb{R}^d)\rightarrow L_2(\mathbb{R}^d)}
\leqslant c_6\beta ^{-1/2}.
\end{split}
\end{equation}

By \eqref{B0^2 L2 to H2}, 
\begin{equation}
\label{cos 0 L2 to H1}
\Vert \cos (t(B^0)^{1/2})(B^0)^{-1}\Vert  _{L_2(\mathbb{R}^d)\rightarrow H^1(\mathbb{R}^d)}\leqslant \Vert (B^0)^{-1}\Vert  _{L_2(\mathbb{R}^d)\rightarrow H^1(\mathbb{R}^d)}
\leqslant  \check{c}^{-1}.
\end{equation}

Next, by Proposition~\ref{Proposition f Pi on H-kappa} with $\varkappa =0$ and \eqref{<b^*b<}, \eqref{Lambda <=}, \eqref{DLambda<=},
\begin{equation}
\label{K1 est start}
\begin{split}
\Vert & \varepsilon \Lambda ^\varepsilon \Pi _\varepsilon b(\mathbf{D})\cos (t (B^0)^{1/2})(B^0)^{-1}\Vert _{L_2(\mathbb{R}^d)\rightarrow H^1(\mathbb{R}^d)}
\\
&\leqslant
(\varepsilon M_1+M_2)\alpha _1^{1/2}\Vert \mathbf{D}(B^0)^{-1}\Vert _{L_2(\mathbb{R}^d)\rightarrow L_2(\mathbb{R}^d)}
\\
&+\varepsilon M_1\alpha _1^{1/2}\Vert \mathbf{D}^2(B^0)^{-1}\Vert _{L_2(\mathbb{R}^d)\rightarrow L_2(\mathbb{R}^d)}.
\end{split}
\end{equation}
Similarly, by Proposition~\ref{Proposition f Pi on H-kappa} and \eqref{tilde Lambda<=}, \eqref{D tilde Lambda},
\begin{equation}
\label{K1 est not final}
\begin{split}
\Vert &\varepsilon \widetilde{\Lambda}^\varepsilon \Pi _\varepsilon \cos (t (B^0)^{1/2})(B^0)^{-1}\Vert _{L_2(\mathbb{R}^d)\rightarrow H^1(\mathbb{R}^d)}
\\
&\leqslant
(\varepsilon\widetilde{M}_1+\widetilde{M}_2)\Vert (B^0)^{-1}\Vert _{L_2(\mathbb{R}^d)\rightarrow L_2(\mathbb{R}^d)}
+\varepsilon\widetilde{M}_1\Vert \mathbf{D}(B^0)^{-1}\Vert _{L_2(\mathbb{R}^d)\rightarrow L_2(\mathbb{R}^d)}.
\end{split}
\end{equation}
Bringing together \eqref{B0^2 L2 to H2}, \eqref{K_1(t,eps)}, \eqref{K1 est start}, and \eqref{K1 est not final}, we obtain
\begin{equation}
\label{K1 L2 to H1}
\Vert \mathcal{K}_1(\varepsilon ;t)\Vert _{L_2(\mathbb{R}^d)\rightarrow H^1(\mathbb{R}^d)}\leqslant C_{17},\quad C_{17}:=(2M_1+M_2)\alpha _1^{1/2}\check{c}^{-1}+(2\widetilde{M}_1+\widetilde{M}_2)\check{c}^{-1}.
\end{equation}
Now from \eqref{cos eps L2 to H1}, \eqref{cos 0 L2 to H1}, and \eqref{K1 L2 to H1} it follows that
\begin{equation}
\label{appr cos H1 on L2}
\begin{split}
\Vert &\cos (t B_\varepsilon ^{1/2})B_\varepsilon ^{-1} -\cos (t (B^0)^{1/2})(B^0)^{-1}-\varepsilon\mathcal{K}_1(\varepsilon ;t)\Vert _{L_2(\mathbb{R}^d)\rightarrow H^1(\mathbb{R}^d)}
\\
&\leqslant c_6\beta ^{-1/2}+\check{c}^{-1}+ C_{17}.
\end{split}
\end{equation}
Interpolating between \eqref{appr cos H1 on L2} and \eqref{Th cos corr}, we arrive at estimate \eqref{Th cos corr interpol} with the constant $C_{14}:=( c_6\beta ^{-1/2}+\check{c}^{-1}+ C_{17})^{1-q}C_{11}^q$.

We proceed to the proof of estimate \eqref{Th sin corr interpol}. By \eqref{b_eps=> H1-norm},
\begin{equation}
\label{sin eps H1 to H1}
\Vert \sin (tB_\varepsilon ^{1/2})B_\varepsilon ^{-1/2}\Vert _{H^1(\mathbb{R}^d)\rightarrow H^1(\mathbb{R}^d)}
\leqslant c_6.
\end{equation}
Next, by \eqref{H1-norm <= B0 1/2},
\begin{equation}
\label{sin 0 H1 to H1}
\Vert \sin (t(B^0)^{1/2})(B^0)^{-1/2}\Vert _{H^1(\mathbb{R}^d)\rightarrow H^1(\mathbb{R}^d)}
\leqslant  \check{c}^{-1/2}.
\end{equation}
By Proposition~\ref{Proposition f Pi on H-kappa}, \eqref{<b^*b<}, \eqref{Lambda <=}, \eqref{DLambda<=}, \eqref{B0 -1/2}, and \eqref{D B0 -1/2},
\begin{equation}
\begin{split}
\Vert &\varepsilon \Lambda ^\varepsilon\Pi _\varepsilon b(\mathbf{D})\sin (t (B^0)^{1/2})(B^0)^{-1/2}\Vert _{H^1(\mathbb{R}^d)\rightarrow H^1(\mathbb{R}^d)}
\\
&\leqslant
(\varepsilon M_1+M_2)\alpha _1^{1/2}\Vert \mathbf{D}(B^0)^{-1/2}\Vert _{H^1(\mathbb{R}^d)\rightarrow L_2(\mathbb{R}^d)}
\\
&+\varepsilon M_1\alpha _1^{1/2}\Vert \mathbf{D}^2 (B^0)^{-1/2}\Vert _{H^1(\mathbb{R}^d)\rightarrow L_2(\mathbb{R}^d)}
\leqslant
(2M_1+M_2)\alpha _1^{1/2}\check{c}^{-1/2}.
\end{split}
\end{equation}
Similarly, by using Proposition~\ref{Proposition f Pi on H-kappa}, \eqref{tilde Lambda<=}, \eqref{D tilde Lambda}, \eqref{B0 -1/2}, and \eqref{D B0 -1/2}, we obtain
\begin{equation}
\label{part K2 in H1}
\begin{split}
\Vert  \varepsilon\widetilde{\Lambda}^\varepsilon \Pi_\varepsilon \sin (t (B^0)^{1/2})(B^0)^{-1/2}\Vert _{H^1(\mathbb{R}^d)\rightarrow H^1(\mathbb{R}^d)}
\leqslant (2\widetilde{M}_1+\widetilde{M}_2)\check{c}^{-1/2}.
\end{split}
\end{equation}
Now from \eqref{K(t,eps)} and \eqref{sin eps H1 to H1}--\eqref{part K2 in H1} it follows that
\begin{equation}
\label{Th sin corr grubo}
\begin{split}
\Vert B_\varepsilon ^{-1/2}\sin(t B_\varepsilon ^{1/2})-(B^0)^{-1/2}\sin(t (B^0)^{1/2})-\varepsilon \mathcal{K}_2(\varepsilon ;t)
\Vert _{H^{1}(\mathbb{R}^d)\rightarrow H^1(\mathbb{R}^d)}
\leqslant C_{18},
\end{split}
\end{equation}
where $C_{18}:=c_6+\check{c}^{-1/2}+(2M_1+M_2)\alpha _1^{1/2}\check{c}^{-1/2}+(2\widetilde{M}_1+\widetilde{M}_2)\check{c}^{-1/2}$.
Interpolating between \eqref{Th sin corr grubo} and \eqref{Th sin corr}, we arrive at estimate \eqref{Th sin corr interpol} with the constant $C_{15}:=C_{18}^{1-q}C_{12}^q$.

Finally, interpolating between \eqref{difference of sine from H-1} and \eqref{Th sin no corr new}, we obtain
estimate \eqref{Th corr int dt}, where $C_{16}:=(c_6+\check{c})^{1-r/2}C_{13}^{r/2}$.
\end{proof}

\subsection{The case where the corrector is equal to zero}
\label{Subsection K=0 in operator terms}
 
Assume that $g^0=\overline{g}$, i.~e., relations \eqref{overline-g} are satisfied. 
In addition, suppose that
\begin{equation}
\label{sum Dj aj =0}
\sum _{j=1}^d D_j a_j(\mathbf{x})^* =0.
\end{equation}
Then the $\Gamma$-periodic solutions of problems \eqref{Lambda problem} and~\eqref{tildeLambda_problem} are equal to zero: $\Lambda (\mathbf{x})=0$ and $\widetilde{\Lambda}(\mathbf{x})=0$. 

\begin{proposition}
\label{Proposition K=0 for sin and cos}
Let the assumptions of  Subsections~\textnormal{\ref{Subsection lattices}--\ref{Subsection Effective operator}} be satisfied. Suppose that relations \eqref{overline-g} and \eqref{sum Dj aj =0} hold. Let $0\leqslant r\leqslant 2$. Then for $0<\varepsilon\leqslant 1$ and $t\in\mathbb{R}$ we have
\begin{align}
\label{Th cos corr=0 int}
\begin{split}
\Vert &\cos (tB_\varepsilon ^{1/2})B_\varepsilon ^{-1}-\cos(t (B^0)^{1/2})(B^0)^{-1}\Vert _{H^{-1+r}(\mathbb{R}^d)\rightarrow H^1(\mathbb{R}^d)}
\\
&\leqslant (c_5c_6^3 + \check{c}^{-1})^{1-r/2} C_{11}^{r/2} \varepsilon ^{r/2} (1+\vert t\vert)^{r/2},
\end{split}
\\
\label{Th sin corr=0 int}
\begin{split}
\Vert &B_\varepsilon ^{-1/2}\sin(t B_\varepsilon ^{1/2})-(B^0)^{-1/2}\sin(t (B^0)^{1/2})\Vert _{H^r(\mathbb{R}^d)\rightarrow H^1(\mathbb{R}^d)}
\\
&\leqslant (c_6+\check{c}^{-1/2})^{1-r/2} C_{12}^{r/2}\varepsilon ^{r/2} (1+\vert t\vert)^{r/2}.
\end{split}
\end{align}
\end{proposition}

\begin{proof}
By Theorem~\ref{Theorem corrector}, for $t\in\mathbb{R}$ and $0<\varepsilon\leqslant 1$ we have
\begin{align}
\label{Th cos corr=0}
\begin{split}
\Vert &\cos (tB_\varepsilon ^{1/2})B_\varepsilon ^{-1}-\cos(t (B^0)^{1/2})(B^0)^{-1}\Vert _{H^1(\mathbb{R}^d)\rightarrow H^1(\mathbb{R}^d)}
\leqslant C_{11} \varepsilon (1+\vert t\vert),
\end{split}
\\
\label{Th sin corr=0}
\begin{split}
\Vert &B_\varepsilon ^{-1/2}\sin(t B_\varepsilon ^{1/2})-(B^0)^{-1/2}\sin(t (B^0)^{1/2})\Vert _{H^2(\mathbb{R}^d)\rightarrow H^1(\mathbb{R}^d)}
\leqslant C_{12}\varepsilon (1+\vert t\vert).
\end{split}
\end{align}

By \eqref{b_eps=> H1-norm}, \eqref{cos eps on H-1}, and the duality arguments,
\begin{equation}
\label{cos res on H-1}
\begin{split}
\Vert & \cos (tB_\varepsilon ^{1/2})B_\varepsilon ^{-1}\Vert _{H^{-1}(\mathbb{R}^d)\rightarrow H^1(\mathbb{R}^d)}\leqslant c_5c_6\Vert B_\varepsilon ^{-1}\Vert _{H^{-1}(\mathbb{R}^d)\rightarrow H^1(\mathbb{R}^d)}
\\
&\leqslant c_5 c_6^2\Vert B_\varepsilon ^{-1/2}\Vert _{H^{-1}(\mathbb{R}^d)\rightarrow L_2(\mathbb{R}^d)}=c_5 c_6^2\Vert B_\varepsilon ^{-1/2}\Vert _{L_2(\mathbb{R}^d)\rightarrow H^1(\mathbb{R}^d)}\leqslant c_5c_6^3.
\end{split}
\end{equation}
Similarly, according to \eqref{H1-norm <= B0 1/2} and \eqref{cos 0 on H-1},
\begin{equation*}
\Vert \cos(t (B^0)^{1/2})(B^0)^{-1}\Vert _{H^{-1}(\mathbb{R}^d)\rightarrow H^1(\mathbb{R}^d)}
\leqslant \check{c}^{-1}.
\end{equation*}
Together with \eqref{cos res on H-1}, this implies
\begin{equation}
\label{cos corr=0 grubo}
\Vert \cos (tB_\varepsilon ^{1/2})B_\varepsilon ^{-1}-\cos(t (B^0)^{1/2})(B^0)^{-1}\Vert _{H^{-1}(\mathbb{R}^d)\rightarrow H^1(\mathbb{R}^d)}
\leqslant c_5c_6^3 + \check{c}^{-1}.
\end{equation}
Interpolating between \eqref{cos corr=0 grubo} and \eqref{Th cos corr=0}, 
we arrive at estimate \eqref{Th cos corr=0 int}.

By \eqref{b_eps=> H1-norm} and \eqref{H1-norm <= B0 1/2},
\begin{equation}
\label{sin corr=0 grubo}
\Vert B_\varepsilon ^{-1/2}\sin(t B_\varepsilon ^{1/2})-(B^0)^{-1/2}\sin(t (B^0)^{1/2})\Vert _{L_2(\mathbb{R}^d)\rightarrow H^1(\mathbb{R}^d)}
\leqslant c_6+\check{c}^{-1/2}.
\end{equation}
Interpolating between \eqref{sin corr=0 grubo} and \eqref{Th sin corr=0}, we get \eqref{Th sin corr=0 int}.
\end{proof}

\subsection{Proof of Theorem~\textnormal{\ref{Theorem corrector}}}
\label{Subsection proof Th sin corr}

\begin{proof}[Proof of Theorem~\textnormal{\ref{Theorem corrector}}]
Denote
\begin{equation*}
\Sigma (t):=e^{-t\mathfrak{A}_\varepsilon}\mathfrak{A}_\varepsilon ^{-1}(I+\varepsilon\mathfrak{G}(\varepsilon))\mathfrak{A}_0^{-1}e^{t\mathfrak{A}_0}.
\end{equation*}
Since all the factors are bounded operators (see \eqref{exp B eps bounded}, \eqref{A_eps^-1}, \eqref{exp B0 bounded}, \eqref{A0 -1}, and \eqref{G(eps)res bounded}), $\Sigma (t)\in \mathcal{B}(H^1(\mathbb{R}^d;\mathbb{C}^n)\times L_2(\mathbb{R}^d;\mathbb{C}^n))$. 
Then, according to semigroup theory (see \cite[Chapter 1, Theorem 2.4(c)]{Pa}),  there exists derivative in the strong topology in $H^1(\mathbb{R}^d;\mathbb{C}^n)\times L_2(\mathbb{R}^d;\mathbb{C}^n)$:
\begin{equation*}
\frac{d\Sigma(t)}{d t}=e^{-t\mathfrak{A}_\varepsilon}(\mathfrak{A}_\varepsilon ^{-1}-\mathfrak{A}_0^{-1}-\varepsilon\mathfrak{G}(\varepsilon)\mathfrak{A}_0^{-1})e^{t\mathfrak{A}_0}
+\varepsilon e^{-t\mathfrak{A}_\varepsilon}\mathfrak{A}_\varepsilon ^{-1}\mathfrak{G}(\varepsilon)e^{t\mathfrak{A}_0}.
\end{equation*}
Thus
\begin{equation*}
\begin{split}
\Sigma(t)-\Sigma(0)
&=
e^{-t\mathfrak{A}_\varepsilon}\mathfrak{A}_\varepsilon ^{-1}(I+\varepsilon\mathfrak{G}(\varepsilon))\mathfrak{A}_0^{-1}e^{t\mathfrak{A}_0}
-\mathfrak{A}_\varepsilon ^{-1}(I+\varepsilon\mathfrak{G}(\varepsilon))\mathfrak{A}_0^{-1}
\\
&=
\int _0^t e^{-s\mathfrak{A}_\varepsilon}(\mathfrak{A}_\varepsilon ^{-1}-\mathfrak{A}_0^{-1}-\varepsilon\mathfrak{G}(\varepsilon)\mathfrak{A}_0^{-1})e^{s\mathfrak{A}_0}\,ds
\\
&+\int_0^t \varepsilon e^{-s\mathfrak{A}_\varepsilon}\mathfrak{A}_\varepsilon ^{-1}\mathfrak{G}(\varepsilon)e^{s\mathfrak{A}_0}\,ds.
\end{split}
\end{equation*}
(The integrals are understood as strong Riemann integrals in the  $H^1(\mathbb{R}^d;\mathbb{C}^n)\times L_2(\mathbb{R}^d;\mathbb{C}^n)$-topology.) 
Let us multiply this identity by $e^{-t\mathfrak{A}_0}$ from the right:
\begin{equation}
\label{tozd G start}
\begin{split}
&e^{-t\mathfrak{A}_\varepsilon}\mathfrak{A}_\varepsilon ^{-1}(I+\varepsilon\mathfrak{G}(\varepsilon))\mathfrak{A}_0^{-1}
-\mathfrak{A}_\varepsilon ^{-1}(I+\varepsilon\mathfrak{G}(\varepsilon))\mathfrak{A}_0^{-1}e^{-t\mathfrak{A}_0}
\\
&=
\int _0^t e^{-s\mathfrak{A}_\varepsilon}(\mathfrak{A}_\varepsilon ^{-1}-\mathfrak{A}_0^{-1}-\varepsilon\mathfrak{G}(\varepsilon)\mathfrak{A}_0^{-1})e^{(s-t)\mathfrak{A}_0}\,ds
\\
&+\int_0^t \varepsilon e^{-s\mathfrak{A}_\varepsilon}\mathfrak{A}_\varepsilon ^{-1}\mathfrak{G}(\varepsilon)e^{(s-t)\mathfrak{A}_0}\,ds.
\end{split}
\end{equation}
We have
\begin{equation*}
\mathfrak{A}_\varepsilon ^{-1}(I+\varepsilon\mathfrak{G}(\varepsilon))\mathfrak{A}_0^{-1}
=
(I+\varepsilon\mathfrak{G}(\varepsilon))\mathfrak{A}_0^{-1}\mathfrak{A}_\varepsilon ^{-1}
+\left[
\mathfrak{A}_\varepsilon ^{-1},(I+\varepsilon\mathfrak{G}(\varepsilon))\mathfrak{A}_0^{-1}
\right].
\end{equation*}
(Here $[\cdot,\cdot]$ denotes the commutator of operators.) Combining this with \eqref{tozd G start}, we obtain
\begin{equation*}
\begin{split}
&e^{-t\mathfrak{A}_\varepsilon}(I+\varepsilon\mathfrak{G}(\varepsilon))\mathfrak{A}_0^{-1}\mathfrak{A}_\varepsilon ^{-1}
-(I+\varepsilon\mathfrak{G}(\varepsilon))\mathfrak{A}_0^{-1}\mathfrak{A}_\varepsilon ^{-1}e^{-t\mathfrak{A}_0}
\\
&=
-e^{-t\mathfrak{A}_\varepsilon}\left[
\mathfrak{A}_\varepsilon ^{-1},(I+\varepsilon\mathfrak{G}(\varepsilon))\mathfrak{A}_0^{-1}
\right]
+\left[
\mathfrak{A}_\varepsilon ^{-1},(I+\varepsilon\mathfrak{G}(\varepsilon))\mathfrak{A}_0^{-1}
\right]e^{-t\mathfrak{A}_0}
\\
&+\int _0^t e^{-s\mathfrak{A}_\varepsilon}(\mathfrak{A}_\varepsilon ^{-1}-\mathfrak{A}_0^{-1}-\varepsilon\mathfrak{G}(\varepsilon)\mathfrak{A}_0^{-1})e^{(s-t)\mathfrak{A}_0}\,ds
\\
&+\int_0^t \varepsilon e^{-s\mathfrak{A}_\varepsilon}\mathfrak{A}_\varepsilon ^{-1}\mathfrak{G}(\varepsilon)e^{(s-t)\mathfrak{A}_0}\,ds.
\end{split}
\end{equation*}
So,
\begin{equation}
\label{tozd G okonchat}
\begin{split}
&e^{-t\mathfrak{A}_\varepsilon}\mathfrak{A}_0^{-2}-(I+\varepsilon\mathfrak{G}(\varepsilon))e^{-t\mathfrak{A}_0}\mathfrak{A}_0^{-2}
=
-e^{-t\mathfrak{A}_\varepsilon}\mathfrak{A}_0^{-1}(\mathfrak{A}_\varepsilon ^{-1}-\mathfrak{A}_0^{-1})
\\
&-e^{-t\mathfrak{A}_\varepsilon}\varepsilon\mathfrak{G}(\varepsilon)\mathfrak{A}_0^{-1}\mathfrak{A}_\varepsilon^{-1}
+(I+\varepsilon\mathfrak{G}(\varepsilon))\mathfrak{A}_0^{-1}(\mathfrak{A}_\varepsilon ^{-1}-\mathfrak{A}_0^{-1})e^{-t\mathfrak{A}_0}
\\
&-e^{-t\mathfrak{A}_\varepsilon}\left[
\mathfrak{A}_\varepsilon ^{-1},(I+\varepsilon\mathfrak{G}(\varepsilon))\mathfrak{A}_0^{-1}
\right]
+\left[
\mathfrak{A}_\varepsilon ^{-1},(I+\varepsilon\mathfrak{G}(\varepsilon))\mathfrak{A}_0^{-1}
\right]e^{-t\mathfrak{A}_0}
\\
&+\int _0^t e^{-s\mathfrak{A}_\varepsilon}(\mathfrak{A}_\varepsilon ^{-1}-\mathfrak{A}_0^{-1}-\varepsilon\mathfrak{G}(\varepsilon)\mathfrak{A}_0^{-1})e^{(s-t)\mathfrak{A}_0}\,ds
\\
&+\int_0^t \varepsilon e^{-s\mathfrak{A}_\varepsilon}\mathfrak{A}_\varepsilon ^{-1}\mathfrak{G}(\varepsilon)e^{(s-t)\mathfrak{A}_0}\,ds
=:\sum _{j=1}^7\mathcal{J}_j(\varepsilon;t).
\end{split}
\end{equation}
According to \eqref{exp =}, \eqref{exp 0 =}, \eqref{A_0^-2}, and \eqref{G(eps)=}, the left-hand side of identity \eqref{tozd G okonchat} has the entries
\begin{align}
\label{left tozd 2 11}
\begin{split}
&\lbrace e^{-t\mathfrak{A}_\varepsilon}\mathfrak{A}_0^{-2}-(I+\varepsilon\mathfrak{G}(\varepsilon))e^{-t\mathfrak{A}_0}\mathfrak{A}_0^{-2}\rbrace _{11}
=-\cos(t B_\varepsilon ^{1/2})(B^0)^{-1}\\
&+(I+\varepsilon \Lambda ^\varepsilon \Pi_\varepsilon b(\mathbf{D})+\varepsilon\widetilde{\Lambda}^\varepsilon\Pi_\varepsilon)\cos(t (B^0)^{1/2})(B^0)^{-1},
\end{split}
\\
\label{left tozd 2 12}
\begin{split}
&\lbrace e^{-t\mathfrak{A}_\varepsilon}\mathfrak{A}_0^{-2}-(I+\varepsilon\mathfrak{G}(\varepsilon))e^{-t\mathfrak{A}_0}\mathfrak{A}_0^{-2}\rbrace _{12}
=
B_\varepsilon ^{-1/2}\sin(t B_\varepsilon ^{1/2})(B^0)^{-1}\\
&-(I+\varepsilon \Lambda ^\varepsilon \Pi_\varepsilon b(\mathbf{D})+\varepsilon\widetilde{\Lambda}^\varepsilon\Pi_\varepsilon)(B^0)^{-1/2}\sin(t(B^0)^{1/2})(B^0)^{-1},
\end{split}
\\
\label{left tozd 2 21}
\begin{split}
&\lbrace e^{-t\mathfrak{A}_\varepsilon}\mathfrak{A}_0^{-2}-(I+\varepsilon\mathfrak{G}(\varepsilon))e^{-t\mathfrak{A}_0}\mathfrak{A}_0^{-2}\rbrace _{21}
\\
&=-B_\varepsilon ^{1/2}\sin(t B_\varepsilon ^{1/2})(B^0)^{-1}+(B^0)^{1/2}\sin(t(B^0)^{1/2})(B^0)^{-1},
\end{split}
\\
\label{left tozd 2 22}
\begin{split}
&\lbrace e^{-t\mathfrak{A}_\varepsilon}\mathfrak{A}_0^{-2}-(I+\varepsilon\mathfrak{G}(\varepsilon))e^{-t\mathfrak{A}_0}\mathfrak{A}_0^{-2}\rbrace _{22}
=-(\cos(t B_\varepsilon ^{1/2})-\cos(t (B^0)^{1/2}))(B^0)^{-1}.
\end{split}
\end{align}

We proceed to consideration of the terms $\mathcal{J}_j(\varepsilon ;t)$, $j=1,\dots,7$. 
By  \eqref{exp =}, \eqref{A0 -1}, and \eqref{I},
\begin{equation}
\label{J1}
\begin{split}
\mathcal{J}_1(\varepsilon;t)
&=
\begin{pmatrix}
0&-B_\varepsilon ^{-1/2}\sin(t B_\varepsilon ^{1/2})\mathcal{R}_1(\varepsilon)\\
0&\cos(t B_\varepsilon ^{1/2})\mathcal{R}_1(\varepsilon)
\end{pmatrix}.
\end{split}
\end{equation}

Next, using \eqref{exp =}, \eqref{A_eps^-1}, and \eqref{G(eps)A0-1}, we obtain
\begin{equation}
\label{J2}
\begin{split}
\mathcal{J}_2(\varepsilon ;t)
&=
\begin{pmatrix}
\cos(t B_\varepsilon^{1/2})\varepsilon K(\varepsilon)&0\\
 B_\varepsilon ^{1/2}\sin(t B_\varepsilon ^{1/2})\varepsilon K(\varepsilon )&0
\end{pmatrix}.
\end{split}
\end{equation}

Now, let us consider the term $\mathcal{J}_3(\varepsilon ;t)$. According to  \eqref{G(eps)A0-1} and \eqref{I},
\begin{equation*}
\mathfrak{G}(\varepsilon)\mathfrak{A}_0^{-1}(\mathfrak{A}_\varepsilon ^{-1}-\mathfrak{A}_0^{-1})
=0.
\end{equation*}
So,
\begin{equation*}
\begin{split}
&\mathcal{J}_3(\varepsilon ;t)
=\mathfrak{A}_0^{-1}(\mathfrak{A}_\varepsilon ^{-1}-\mathfrak{A}_0^{-1})e^{-t \mathfrak{A}_0}
\end{split}
\end{equation*}
and (by \eqref{exp 0 =}, \eqref{A0 -1}, and \eqref{I})
\begin{align}
\label{J3}
&\lbrace \mathcal{J}_3(\varepsilon ;t)\rbrace _{11}=\lbrace \mathcal{J}_3(\varepsilon ;t)\rbrace _{12}=0,
\\
\label{J3 21}
&\lbrace \mathcal{J}_3(\varepsilon ;t)\rbrace _{21}=-\mathcal{R}_1(\varepsilon )(B^0)^{1/2}\sin (t (B^0)^{1/2}),
\\
&\lbrace \mathcal{J}_3(\varepsilon ;t)\rbrace _{22}=-\mathcal{R}_1(\varepsilon )\cos(t(B^0) ^{1/2}).
\nonumber
\end{align}

To evaluate the matrix entries of the terms $\mathcal{J}_4(\varepsilon ;t)$ and $\mathcal{J}_5(\varepsilon ;t)$, we need to calculate the commutator $[\mathfrak{A}_\varepsilon ^{-1},(I+\varepsilon\mathfrak{G}(\varepsilon))\mathfrak{A}_0^{-1}]$. By \eqref{A0 -1} and \eqref{G(eps)A0-1},
\begin{equation*}
(I+\varepsilon\mathfrak{G}(\varepsilon))\mathfrak{A}_0^{-1}
=
\begin{pmatrix}
0&-(B^0)^{-1}-\varepsilon K(\varepsilon)\\
I&0
\end{pmatrix}.
\end{equation*}
Together with \eqref{A_eps^-1}, this implies
\begin{align}
\label{comm 2}
&\mathfrak{A}_\varepsilon ^{-1}(I+\varepsilon\mathfrak{G}(\varepsilon))\mathfrak{A}_0^{-1}
=
\begin{pmatrix}
-B_\varepsilon ^{-1}&0\\
0&-(B^0)^{-1}-\varepsilon K(\varepsilon)
\end{pmatrix},
\\
\label{comm 3}
&(I+\varepsilon\mathfrak{G}(\varepsilon))\mathfrak{A}_0^{-1}\mathfrak{A}_\varepsilon ^{-1}
=
\begin{pmatrix}
-(B^0)^{-1}-\varepsilon K(\varepsilon)&0\\
0&-B_\varepsilon ^{-1}
\end{pmatrix}.
\end{align}
Recall notation \eqref{R2}. Combining \eqref{comm 2} and \eqref{comm 3}, we obtain
\begin{equation}
\label{commutator}
\begin{split}
\left[
\mathfrak{A}_\varepsilon ^{-1},(I+\varepsilon\mathfrak{G}(\varepsilon))\mathfrak{A}_0^{-1}
\right]
&=\mathfrak{A}_\varepsilon ^{-1}(I+\varepsilon\mathfrak{G}(\varepsilon))\mathfrak{A}_0^{-1}
-(I+\varepsilon\mathfrak{G}(\varepsilon))\mathfrak{A}_0^{-1}\mathfrak{A}_\varepsilon ^{-1}
\\
&=
\begin{pmatrix}
-\mathcal{R}_2(\varepsilon)&0\\
0&\mathcal{R}_2(\varepsilon)
\end{pmatrix}
.
\end{split}
\end{equation}

Now from \eqref{exp =} and \eqref{commutator} it follows that
\begin{equation}
\label{J4}
\begin{split}
&\mathcal{J}_4(\varepsilon;t)
=
\begin{pmatrix}
\cos(t B_\varepsilon^{1/2})\mathcal{R}_2(\varepsilon)&
B_\varepsilon ^{-1/2}\sin(t B_\varepsilon ^{1/2})\mathcal{R}_2(\varepsilon)\\
B_\varepsilon ^{1/2}\sin(t B_\varepsilon ^{1/2})\mathcal{R}_2(\varepsilon)&-\cos(t B_\varepsilon^{1/2})\mathcal{R}_2(\varepsilon)
\end{pmatrix}.
\end{split}
\end{equation}

By \eqref{exp 0 =} and \eqref{commutator},
\begin{equation}
\label{J5}
\begin{split}
\mathcal{J}_5(\varepsilon;t)&=
\begin{pmatrix}
-\mathcal{R}_2(\varepsilon)\cos(t (B^0)^{1/2})&\mathcal{R}_2(\varepsilon)(B^0)^{-1/2}\sin(t (B^0)^{1/2})\\
\mathcal{R}_2(\varepsilon)(B^0)^{1/2}\sin(t (B^0)^{1/2})&
\mathcal{R}_2(\varepsilon)\cos(t (B^0)^{1/2})
\end{pmatrix}.
\end{split}
\end{equation}

Next, by using \eqref{exp =}, \eqref{exp 0 =}, and \eqref{II}, we calculate $\mathcal{J}_6(\varepsilon;t)$. It turns out that 
\begin{align}
\label{J6 11}
&\lbrace \mathcal{J}_6(\varepsilon;t)\rbrace _{11}=\int _0^t \cos(s B_\varepsilon^{1/2})\mathcal{R}_2(\varepsilon)(B^0)^{1/2}\sin((s-t) (B^0)^{1/2})\,ds,
\\
\label{J6 12}
&\lbrace \mathcal{J}_6(\varepsilon;t)\rbrace _{12}
=-\int _0^t \cos(s B_\varepsilon^{1/2})\mathcal{R}_2(\varepsilon)\cos((s-t) (B^0)^{1/2})\,ds
,
\\
\label{J6 21}
&\lbrace \mathcal{J}_6(\varepsilon;t)\rbrace _{21}
=\int _0^t B_\varepsilon ^{1/2}\sin(s B_\varepsilon ^{1/2})\mathcal{R}_2(\varepsilon)(B^0)^{1/2}\sin((s-t) (B^0)^{1/2})\,ds,
\\
&\lbrace \mathcal{J}_6(\varepsilon;t)\rbrace _{22}
=-\int _0^t B_\varepsilon ^{1/2}\sin(s B_\varepsilon ^{1/2})\mathcal{R}_2(\varepsilon)\cos((s-t) (B^0)^{1/2})\,ds.
\nonumber
\end{align} 

Finally, by \eqref{exp =}, \eqref{A_eps^-1}, \eqref{exp 0 =}, and \eqref{G(eps)=},
\begin{align}
\label{J7 11}
\begin{split}
&\lbrace \mathcal{J}_7(\varepsilon;t)\rbrace _{11}
\\
&=-\varepsilon \int _0^t B_\varepsilon ^{-1/2}\sin(sB_\varepsilon^{1/2}) (\Lambda ^\varepsilon \Pi _\varepsilon b(\mathbf{D})+\widetilde{\Lambda}^\varepsilon\Pi _\varepsilon)\cos((s-t) (B^0)^{1/2})\,ds,
\end{split}
\\
\label{J7 12}
\begin{split}
&\lbrace \mathcal{J}_7(\varepsilon;t)\rbrace _{12}
\\
&=-\varepsilon \int _0^t B_\varepsilon ^{-1/2}\sin(sB_\varepsilon^{1/2}) (\Lambda ^\varepsilon \Pi _\varepsilon b(\mathbf{D})+\widetilde{\Lambda}^\varepsilon\Pi _\varepsilon)(B^0)^{-1/2}\sin((s-t) (B^0)^{1/2})\,ds,
\end{split}
\\
\label{J7 21}
\begin{split}
&\lbrace \mathcal{J}_7(\varepsilon;t)\rbrace _{21}
=\varepsilon \int _0^t \cos(s B_\varepsilon^{1/2}) (\Lambda ^\varepsilon \Pi _\varepsilon b(\mathbf{D})+\widetilde{\Lambda}^\varepsilon\Pi _\varepsilon)\cos((s-t) (B^0)^{1/2})\,ds,
\end{split}
\\
\begin{split}
&\lbrace \mathcal{J}_7(\varepsilon;t)\rbrace _{22}
\\
&=\varepsilon \int _0^t \cos(s B_\varepsilon^{1/2}) (\Lambda ^\varepsilon \Pi _\varepsilon b(\mathbf{D})+\widetilde{\Lambda}^\varepsilon\Pi _\varepsilon)(B^0)^{-1/2}\sin((s-t) (B^0)^{1/2})\,ds.
\end{split}
\nonumber
\end{align}

Combining \eqref{tozd G okonchat}, \eqref{left tozd 2 11}, \eqref{J1}--\eqref{J3}, \eqref{J4}--\eqref{J6 11}, and \eqref{J7 11},  we obtain
\begin{equation*}
\begin{split}
-&\cos(t B_\varepsilon ^{1/2})(B^0)^{-1}+(I+\varepsilon \Lambda ^\varepsilon \Pi_\varepsilon b(\mathbf{D})+\varepsilon\widetilde{\Lambda}^\varepsilon\Pi_\varepsilon)\cos(t (B^0)^{1/2})(B^0)^{-1}
\\
&=\cos(t B_\varepsilon^{1/2})\varepsilon K(\varepsilon)+\cos(t B_\varepsilon^{1/2})\mathcal{R}_2(\varepsilon)
-\mathcal{R}_2(\varepsilon)\cos(t (B^0)^{1/2})
\\
&+\int _0^t \cos(s B_\varepsilon^{1/2})\mathcal{R}_2(\varepsilon)(B^0)^{1/2}\sin((s-t) (B^0)^{1/2})\,ds
\\
&-\varepsilon \int _0^t B_\varepsilon ^{-1/2}\sin(sB_\varepsilon^{1/2}) (\Lambda ^\varepsilon \Pi _\varepsilon b(\mathbf{D})+\widetilde{\Lambda}^\varepsilon\Pi _\varepsilon)\cos((s-t) (B^0)^{1/2})\,ds.
\end{split}
\end{equation*}
It is the equality for the operators, acting from $H^1$ to $H^1$. Taking \eqref{R2} and \eqref{K_1(t,eps)} into account, we arrive at
\begin{equation*}
\begin{split}
&\cos (tB_\varepsilon)^{1/2}B_\varepsilon ^{-1}-\cos(t (B^0)^{1/2})(B^0)^{-1}
-\varepsilon\mathcal{K}_1(\varepsilon ;t)
\\
&=\mathcal{R}_2(\varepsilon)\cos(t (B^0)^{1/2})-\int _0^t \cos(s B_\varepsilon^{1/2})\mathcal{R}_2(\varepsilon)(B^0)^{1/2}\sin((s-t) (B^0)^{1/2})\,ds
\\
&+\varepsilon \int _0^t B_\varepsilon ^{-1/2}\sin(sB_\varepsilon^{1/2}) (\Lambda ^\varepsilon \Pi _\varepsilon b(\mathbf{D})+\widetilde{\Lambda}^\varepsilon\Pi _\varepsilon)\cos((s-t) (B^0)^{1/2})\,ds.
\end{split}
\end{equation*}
Together with \eqref{b_eps=> H1-norm}, \eqref{B01/2 on H1}, \eqref{B-eps with corrector}, \eqref{Beps1/2 appr corr}, and \eqref{partof corr H1 to L2}, this implies \eqref{Th cos corr} with the constant 
$C_{11}:=\max\lbrace C_2;c_5c_6C_L^{1/2}C_2+c_6(M_1\alpha_1 ^{1/2}+\widetilde{M}_1)\rbrace .$

Recall the notation \eqref{K(t,eps)}. Combining \eqref{tozd G okonchat}, \eqref{left tozd 2 12}, \eqref{J1}--\eqref{J3}, \eqref{J4}, \eqref{J5}, \eqref{J6 12}, and \eqref{J7 12}, 
we arrive at
\begin{equation}
\label{sin corr identity almost final}
\begin{split}
&B_\varepsilon ^{-1/2}\sin(t B_\varepsilon ^{1/2})(B^0)^{-1}-
(B^0)^{-1/2}\sin(t(B^0)^{1/2})(B^0)^{-1}-\varepsilon\mathcal{K}_2(\varepsilon ;t)(B^0)^{-1}
\\
&=
-\varepsilon B_\varepsilon ^{-1/2}\sin(t B_\varepsilon ^{1/2})K(\varepsilon)
+\mathcal{R}_2(\varepsilon)(B^0)^{-1/2}\sin(t (B^0)^{1/2})
\\
&-\int _0^t \cos(s B_\varepsilon^{1/2})\mathcal{R}_2(\varepsilon)\cos((s-t) (B^0)^{1/2})\,ds
\\
&-\varepsilon \int _0^t B_\varepsilon ^{-1/2}\sin(sB_\varepsilon^{1/2}) (\Lambda ^\varepsilon \Pi _\varepsilon b(\mathbf{D})+\widetilde{\Lambda}^\varepsilon\Pi _\varepsilon)(B^0)^{-1/2}\sin((s-t) (B^0)^{1/2})\,ds.
\end{split}
\end{equation}
Here the equality $\mathcal{R}_2(\varepsilon)-\mathcal{R}_1(\varepsilon)=-\varepsilon K(\varepsilon)$ was taken into account. 
By construction,  identity \eqref{sin corr identity almost final} should be understood in the $(H^1\to L_2)$-sense.
By continuity, we extend domains of the operators from the left- and right-hand sides of 
identity \eqref{sin corr identity almost final} onto the space $L_2$. Moreover, the ranges of all operators in \eqref{sin corr identity almost final} lie in $H^1$.

By Theorem~\ref{Theorem elliptic}, the elementary inequality $\vert \sin x\vert /\vert x\vert\leqslant 1$, $x\in\mathbb{R}$, and \eqref{b_eps=> H1-norm}, \eqref{K(eps)}, \eqref{K(eps)<= L2 to L2}, \eqref{Beps1/2 appr corr}, \eqref{K(eps)B01/2}, and \eqref{sin corr identity almost final},
\begin{equation*}
\begin{split}
\Vert &\left(B_\varepsilon ^{-1/2}\sin(t B_\varepsilon ^{1/2})-
(B^0)^{-1/2}\sin(t(B^0)^{1/2})-\varepsilon\mathcal{K}_2(\varepsilon ;t)\right)(B^0)^{-1}\Vert _{L_2(\mathbb{R}^d)\rightarrow H^1(\mathbb{R}^d)}
\\
&\leqslant
c_6C_K\varepsilon +(1+c_5c_6)C_2\varepsilon\vert t\vert +c_6\check{c}^{-1/2}(M_1\alpha _1^{1/2}+\widetilde{M}_1)\varepsilon\vert t\vert .
\end{split}
\end{equation*}
Together with \eqref{B0 on H2}, this implies estimate \eqref{Th sin corr} with the constant $C_{12}:=C_L \max\lbrace c_6C_K;(1+c_5c_6)C_2+c_6\check{c}^{-1/2}(M_1\alpha _1^{1/2}+\widetilde{M}_1)\rbrace$.

Combining \eqref{tozd G okonchat}, \eqref{left tozd 2 21}, \eqref{J1}, \eqref{J2}, \eqref{J3 21}, \eqref{J4}, \eqref{J5}, \eqref{J6 21}, and \eqref{J7 21}, 
we arrive at the $(H^1\rightarrow L_2)$-equality
\begin{equation*}
\begin{split}
-&B_\varepsilon ^{1/2}\sin(t B_\varepsilon ^{1/2})(B^0)^{-1}+(B^0)^{1/2}\sin(t(B^0)^{1/2})(B^0)^{-1}
\\
&=B_\varepsilon ^{1/2}\sin(t B_\varepsilon ^{1/2})\varepsilon K(\varepsilon)-\mathcal{R}_1(\varepsilon )(B^0)^{1/2}\sin (t (B^0)^{1/2})
\\
&+
B_\varepsilon ^{1/2}\sin(t B_\varepsilon ^{1/2})\mathcal{R}_2(\varepsilon)
+\mathcal{R}_2(\varepsilon)(B^0)^{1/2}\sin(t (B^0)^{1/2})
\\
&+\int _0^t B_\varepsilon ^{1/2}\sin(s B_\varepsilon ^{1/2})\mathcal{R}_2(\varepsilon)(B^0)^{1/2}\sin((s-t) (B^0)^{1/2})\,ds
\\
&+\varepsilon \int _0^t \cos(s B_\varepsilon^{1/2}) (\Lambda ^\varepsilon \Pi _\varepsilon b(\mathbf{D})+\widetilde{\Lambda}^\varepsilon\Pi _\varepsilon)\cos((s-t) (B^0)^{1/2})\,ds.
\end{split}
\end{equation*}
Taking \eqref{R1}, \eqref{R2}, and \eqref{K(t,eps)} into account, we get 
\begin{equation*}
\begin{split}
&B_\varepsilon ^{-1/2}\sin(t B_\varepsilon^{1/2})-(B^0)^{-1/2}\sin (t (B^0)^{1/2})-\varepsilon\mathcal{K}_2(\varepsilon ;t)
\\
&=-\int _0^t B_\varepsilon ^{1/2}\sin(s B_\varepsilon ^{1/2})\mathcal{R}_2(\varepsilon)(B^0)^{1/2}\sin((s-t) (B^0)^{1/2})\,ds
\\
&-\varepsilon \int _0^t \cos(s B_\varepsilon^{1/2}) (\Lambda ^\varepsilon \Pi _\varepsilon b(\mathbf{D})+\widetilde{\Lambda}^\varepsilon\Pi _\varepsilon)\cos((s-t) (B^0)^{1/2})\,ds.
\end{split}
\end{equation*}
(The equality here is understood in the $(H^1\rightarrow L_2)$-sense.) Together with \eqref{B01/2 on H1}, \eqref{Beps1/2 appr corr}, and \eqref{partof corr H1 to L2}, this implies
\begin{equation}
\label{Th sin new almost fin}
\begin{split}
\Vert &B_\varepsilon ^{-1/2}\sin(t B_\varepsilon^{1/2})-(B^0)^{-1/2}\sin (t (B^0)^{1/2})-\varepsilon\mathcal{K}_2(\varepsilon ;t)\Vert _{H^1(\mathbb{R}^d)\rightarrow L_2(\mathbb{R}^d)}
\\
&\leqslant (c_5C_2C_L^{1/2}+M_1\alpha _1^{1/2}+\widetilde{M}_1) \varepsilon \vert t\vert .
\end{split}
\end{equation}
Next, since the operator $B^0$ commutes with differentiation, 
by \eqref{partof corr H1 to L2} and the elementary inequality $\vert \sin x\vert /\vert x\vert \leqslant 1$,  we have
\begin{equation*}
\begin{split}
\Vert \mathcal{K}_2(\varepsilon ;t)\Vert _{H^1(\mathbb{R}^d)\rightarrow L_2(\mathbb{R}^d)}
&\leqslant (M_1\alpha _1^{1/2}+\widetilde{M}_1)
\Vert (B^0)^{-1/2}\sin(t(B^0)^{1/2})\Vert _{H^1(\mathbb{R}^d)\rightarrow H^1(\mathbb{R}^d)}
\\
&\leqslant (M_1\alpha _1^{1/2}+\widetilde{M}_1)
\vert t\vert .
\end{split}
\end{equation*}
Together with \eqref{Th sin new almost fin}, this implies estimate \eqref{Th sin no corr new} with the constant $C_{13}:=c_5C_2C_L^{1/2}+2M_1\alpha _1^{1/2}+2\widetilde{M}_1$.
\end{proof}

\begin{remark}
Consideration of element \eqref{left tozd 2 22} does not give us any new information compared with estimate \eqref{Th cos L2}.
\end{remark}

\section{Removal of the smoothing operator from the corrector}

\label{Section Removal of the smoothing operator from the corrector}

\subsection{Removal of $\Pi _\varepsilon$ under additional assumptions on $\Lambda$ and $\widetilde{\Lambda}$}

It turns out that the smoothing operator can be removed from the corrector if the matrix-valued functions  $\Lambda (\mathbf{x})$ and $\widetilde{\Lambda}(\mathbf{x})$ are subject to some additional assumptions. 

\begin{condition}
\label{Condition Lambda in L infty}
Assume that the $\Gamma$-periodic solution $\Lambda (\mathbf{x})$ of problem \eqref{Lambda problem} is bounded, i.~e., $\Lambda\in L_\infty (\mathbb{R}^d)$.
\end{condition}

The cases when the Condition~\ref{Condition Lambda in L infty} is fulfilled automatically  were distinguished in \cite[Lemma~8.7]{BSu06}.

\begin{proposition}
\label{Proposition Lambda in L infty <=}
Suppose that at least one of the following assumptions is satisfied\textnormal{:}

\noindent
$1^\circ )$ $d\leqslant 2${\rm ;}

\noindent
$2^\circ )$ the dimension $d\geqslant 1$ is arbitrary, and the differential expression $A_\varepsilon$ is given by $A_\varepsilon =\mathbf{D}^* g^\varepsilon (\mathbf{x})\mathbf{D}$, where $g(\mathbf{x})$ is a symmetric matrix with real entries{\rm ;}

\noindent
$3^\circ )$ the dimension $d$ is arbitrary, and $g^0=\underline{g}$, i.~e., relations \eqref{underline-g} are satisfied.

\noindent
Then Condition~\textnormal{\ref{Condition Lambda in L infty}}  holds.
\end{proposition}

In order to remove $\Pi_\varepsilon$ from the term involving $\widetilde{\Lambda}^\varepsilon$, it suffices to impose the following condition.

\begin{condition}
\label{Condition tilde Lambda in Lp}
Assume that the $\Gamma$-periodic solution $\widetilde{\Lambda}(\mathbf{x})$ of problem \eqref{tildeLambda_problem} is such that
\begin{equation*}
\widetilde{\Lambda}\in L_p(\Omega),\quad p=2 \;\mbox{for}\;d=1,\quad p>2\;\mbox{for}\;d=2, \quad p=d \;\mbox{for}\;d\geqslant 3.
\end{equation*}
\end{condition}

The following result was obtained in \cite[Proposition 8.11]{SuAA}.

\begin{proposition}
\label{Proposition tilde Lambda in Lp if}
Condition \textnormal{\ref{Condition tilde Lambda in Lp}} is fulfilled, if at least one of the following assumptions is satisfied\textnormal{:}

\noindent
$1^\circ )$ $d\leqslant 4${\rm ;}

\noindent
$2^\circ )$ the dimension $d$ is arbitrary, and the differential expression $A_{\varepsilon}$ has the form $A_{\varepsilon} =\mathbf{D}^* g^\varepsilon (\mathbf{x})\mathbf{D}$, where $g(\mathbf{x})$ is a~symmetric matrix-valued function with real entries. 
\end{proposition}

\begin{remark}
\label{Remark scalar problem}
If $A_{\varepsilon} =\mathbf{D}^* g^\varepsilon (\mathbf{x})\mathbf{D}$, where $g(\mathbf{x})$ is a~symmetric matrix-valued function with real entries, from \textnormal{\cite[Chapter {\rm III}, Theorem {\rm 13.1}]{LaU}} it follows that $\Lambda,\widetilde{\Lambda}\in L_\infty$ and the norm $\Vert\Lambda\Vert _{L_\infty}$ is controlled in terms of $d$, $\Vert g\Vert _{L_\infty}$, $\Vert g^{-1}\Vert _{L_\infty}$, and $\Omega$\textnormal{;} the norm $\Vert \widetilde{\Lambda}\Vert _{L_\infty }$ does not exceed a constant depending on $d$, $\rho$, $\Vert g\Vert _{L_\infty}$, $\Vert g^{-1}\Vert _{L_\infty}$, $\Vert a_j\Vert _{L_\rho (\Omega)}$, $j=1,\dots ,d$, and $\Omega$. In this case, Conditions \textnormal{\ref{Condition Lambda in L infty}} and \textnormal{\ref{Condition tilde Lambda in Lp}} are fulfilled simultaneously.
\end{remark}

\begin{theorem}
\label{Theorem removal Pi eps under conditions}
Let the assumptions of Subsections~\textnormal{\ref{Subsection lattices}--\ref{Subsection Effective operator}} be satisfied. Assume that Conditions \textnormal{\ref{Condition Lambda in L infty}} and \textnormal{\ref{Condition tilde Lambda in Lp}} hold. 
Denote
\begin{align}
\label{K_10(t,eps)}
&\mathcal{K}_1^0(\varepsilon ;t):=
(\Lambda ^\varepsilon  b(\mathbf{D})+ \widetilde{\Lambda}^\varepsilon ) \cos(t (B^0)^{1/2})(B^0)^{-1},
\\
\label{K20(t,eps)}
&\mathcal{K}_2^0(\varepsilon ;t):=
(\Lambda ^\varepsilon  b(\mathbf{D})+ \widetilde{\Lambda}^\varepsilon)(B^0)^{-1/2}\sin(t (B^0)^{1/2}).
\end{align}
Let $0\leqslant q\leqslant 1$. 
Then for $0<\varepsilon\leqslant 1$ and $t\in\mathbb{R}$ we have
\begin{align*}
\begin{split}
\Vert &\cos (tB_\varepsilon ^{1/2})B_\varepsilon ^{-1}-\cos(t (B^0)^{1/2})(B^0)^{-1}-\varepsilon \mathcal{K}_1^0(\varepsilon ;t)\Vert _{H^q(\mathbb{R}^d)\rightarrow H^1(\mathbb{R}^d)}
\\
&\leqslant {C}_{19}\varepsilon ^q (1+\vert t\vert)^q,
\end{split}
\\
\begin{split}
\Vert &B_\varepsilon ^{-1/2}\sin(t B_\varepsilon ^{1/2})-(B^0)^{-1/2}\sin(t (B^0)^{1/2})-\varepsilon \mathcal{K}_2^0(\varepsilon ;t)\Vert _{H^{1+q}(\mathbb{R}^d)\rightarrow H^1(\mathbb{R}^d)}
\\
&\leqslant  {C}_{20}\varepsilon ^q (1+\vert t\vert)^q.
\end{split}
\end{align*}
The constants ${C}_{19}$ and ${C}_{20}$ depend only on $q$, $p$, $\Vert \Lambda \Vert _{L_\infty}$, $\Vert \widetilde{\Lambda}\Vert _{L_p(\Omega)}$, and on the problem data \eqref{problem data}.
\end{theorem}

Proof of Theorem~\ref{Theorem removal Pi eps under conditions} is given below in Subsection~\ref{Proof of theorem corr no Pi}. In the proof we use some properties of the matrix-valued functions $\Lambda$ and $\widetilde{\Lambda}$ from the next subsection.

\subsection{Properties of the matrix-valued functions $\Lambda$ and $\widetilde{\Lambda}$}
\label{Subsection properties of Lambda and tilde Lambda}

The following results were obtained in~\cite[Corollary~2.4]{PSu} and \cite[Lemma 3.5 and Corollary~3.6]{MSu15}.

\begin{lemma}
\label{Lemma Lambda in L infty}
Suppose that $\Lambda $ is the $\Gamma$-periodic solution of problem \eqref{Lambda problem}.
Assume also that $\Lambda \in L_\infty$.
Then for any $u\in H^1 (\mathbb{R}^d)$ and $\varepsilon >0$ we have
\begin{equation*}
\int _{\mathbb{R}^d}\vert (\mathbf{D}\Lambda )^\varepsilon (\mathbf{x})\vert ^2\vert u(\mathbf{x})\vert ^2\,d\mathbf{x}
\leqslant\mathfrak{c} _1\Vert u\Vert ^2_{L_2(\mathbb{R}^d)}
+\mathfrak{c} _2 \varepsilon ^2 \|\Lambda\|^2_{L_\infty} \int _{\mathbb{R}^d} \vert \mathbf{D}u(\mathbf{x})\vert ^2\,d\mathbf{x}.
\end{equation*}
The constants $\mathfrak{c} _1$ and $\mathfrak{c} _2$ depend only on
$m$, $d$, $\alpha _0$, $\alpha _1$, $\Vert g\Vert _{L_\infty}$, and $\Vert g^{-1}\Vert _{L_\infty}$.
\end{lemma}

\begin{lemma}
\label{Lemma_Lambda_tilda3}
Suppose that the $\Gamma$-periodic solution $\widetilde{\Lambda}$ of problem \eqref{tildeLambda_problem} satisfies Condition~\textnormal{\ref{Condition tilde Lambda in Lp}}.
Then for $0< \varepsilon \leqslant 1$ the operator $[\widetilde{\Lambda} ^\varepsilon]$ is a continuous mapping from $H^1({\mathbb R}^d)$ to
$L_2({\mathbb R}^d)$, and
$
\| [\widetilde{\Lambda}^\varepsilon] \|_{H^1({\mathbb R}^d) \to L_2({\mathbb R}^d)} \leqslant \|\widetilde{\Lambda}\|_{L_p(\Omega)} C^{(p)}_{\Omega},
$
where $C^{(p)}_{\Omega}$ is the norm of the embedding  $H^1(\Omega )\hookrightarrow L_{2(p/2)'}(\Omega)$.
Here \hbox{$(p/2)'= \infty$} for $d=1$, and $(p/2)'= p (p-2)^{-1}$ for $d\geqslant 2$.
\end{lemma}

\begin{lemma}
\label{Lemma_Lambda_tilda4}
Suppose that $\widetilde{\Lambda}$ is the $\Gamma$-periodic solution of problem \eqref{tildeLambda_problem}.
Suppose also that $\widetilde{\Lambda}$ satisfies Condition~\textnormal{\ref{Condition tilde Lambda in Lp}}. Then for any $u \in H^2({\mathbb R}^d)$ and $0< \varepsilon \leqslant 1$ we have
\begin{equation*}
\begin{split}
\int _{\mathbb{R}^d}\vert (\mathbf{D}\widetilde{\Lambda})^\varepsilon(\mathbf{x})\vert ^2\vert u(\mathbf{x})\vert ^2\,d\mathbf{x}
\leqslant \widetilde{\mathfrak{c}} _1 \Vert u\Vert ^2_{H^1(\mathbb{R}^d)}
+\widetilde{\mathfrak{c}} _2 \varepsilon ^2
\| \widetilde{\Lambda} \|_{L_p(\Omega)}^2 (C_\Omega^{(p)})^2 \| \mathbf{D}u \|^2_{H^1({\mathbb R}^d)}.
\end{split}
\end{equation*}
The constants $\widetilde{\mathfrak{c}} _1$ and $\widetilde{\mathfrak{c}} _2$
depend only on $n$, $d$, $\alpha _0$, $\alpha _1$, $\rho$, $\Vert g\Vert _{L_\infty}$, $\Vert g^{-1}\Vert _{L_\infty}$,
the norms $\Vert a_j\Vert _{L_\rho (\Omega)}$, $j=1,\dots ,d$, and the parameters of the lattice $\Gamma$.
\end{lemma}

To prove Theorem~\ref{Theorem removal Pi eps under conditions}, we need the following lemmas. With the another smoothing operator (the Steklov smoothing), these lemmas were proven in 
\cite[Lemmas 7.7 and 7.8]{MSuPOMI}. For our smoothing operator $\Pi _\varepsilon$, the proof is quite similar. 

\begin{lemma}
\label{Lemma Lambda (S-I)}
Suppose that Condition~\textnormal{\ref{Condition Lambda in L infty}} is satisfied.
Let $\Pi_\varepsilon$ be the operator~\eqref{Pi eps = definition}. Then for~$0<\varepsilon\leqslant 1$ we have
\begin{equation}
\label{lemma Lambda (S-I)}
\Vert [\Lambda ^\varepsilon ]b(\mathbf{D})(\Pi _\varepsilon -I)\Vert _{H^2(\mathbb{R}^d)\rightarrow H^1(\mathbb{R}^d)}\leqslant\mathfrak{C}_\Lambda .
\end{equation}
The constant~$\mathfrak{C}_\Lambda$ depends only on $m$, $d$, $\alpha _0$, $\alpha _1$, $\Vert g\Vert _{L_\infty}$, $\Vert g^{-1}\Vert _{L_\infty}$, the parameters of the lattice~$\Gamma$, and the norm~$\Vert \Lambda\Vert _{L_\infty}$.
\end{lemma}

\begin{proof}
Let $\boldsymbol{\Phi}\in H^2(\mathbb{R}^d;\mathbb{C}^n)$. From  \eqref{<b^*b<}, \eqref{Pi eps <= 1}, and Condition~\ref{Condition Lambda in L infty}, it follows that
\begin{equation}
\label{1 lemma Lambda (S-I)}
\Vert \Lambda ^\varepsilon b(\mathbf{D})(\Pi _\varepsilon -I)\boldsymbol{\Phi}\Vert _{L_2(\mathbb{R}^d)}
\leqslant 2\alpha _1^{1/2}\Vert \Lambda \Vert _{L_\infty}\Vert \mathbf{D}\boldsymbol{\Phi}\Vert _{L_2(\mathbb{R}^d)}.
\end{equation}

Consider the derivatives:
\begin{equation*}
\partial _j\left( \Lambda ^\varepsilon b(\mathbf{D})(\Pi _\varepsilon -I)\boldsymbol{\Phi}\right)
=\varepsilon ^{-1}(\partial _j \Lambda )^\varepsilon (\Pi _\varepsilon -I)b(\mathbf{D})\boldsymbol{\Phi}
+\Lambda ^\varepsilon (\Pi _\varepsilon -I) b(\mathbf{D})\partial _j \boldsymbol{\Phi}.
\end{equation*}
Hence,
\begin{equation*}
\begin{split}
\Vert \mathbf{D}\left(\Lambda ^\varepsilon b(\mathbf{D})(\Pi _\varepsilon -I)\boldsymbol{\Phi}\right)\Vert ^2 _{L_2(\mathbb{R}^d)}
&\leqslant
2\varepsilon ^{-2}\Vert (\mathbf{D}\Lambda )^\varepsilon (\Pi _\varepsilon -I)b(\mathbf{D})\boldsymbol{\Phi}\Vert ^2 _{L_2(\mathbb{R}^d)}
\\
&+2\Vert \Lambda \Vert ^2_{L_\infty}\Vert (\Pi _\varepsilon -I)b(\mathbf{D})\mathbf{D}\boldsymbol{\Phi}\Vert ^2 _{L_2(\mathbb{R}^d)}.
\end{split}
\end{equation*}
By Lemma~\ref{Lemma Lambda in L infty}, this yields
\begin{equation*}
\begin{split}
\Vert \mathbf{D}\left(\Lambda ^\varepsilon b(\mathbf{D})(\Pi _\varepsilon -I)\boldsymbol{\Phi}\right)\Vert ^2 _{L_2(\mathbb{R}^d)}
&\leqslant
2\mathfrak{c} _1 \varepsilon ^{-2}\Vert (\Pi _\varepsilon -I)b(\mathbf{D})\boldsymbol{\Phi}\Vert ^2 _{L_2(\mathbb{R}^d)}
\\
&+ 2\Vert \Lambda \Vert ^2_{L_\infty}(\mathfrak{c} _2+1)\Vert (\Pi _\varepsilon -I)b(\mathbf{D})\mathbf{D}\boldsymbol{\Phi}\Vert ^2_{L_2(\mathbb{R}^d)}.
\end{split}
\end{equation*}
So, combining \eqref{<b^*b<}, \eqref{Pi eps <= 1},  and Proposition~\ref{Proposition Pi -I}, we get
\begin{equation}
\label{2 lemma Lambda (S-I)}
\begin{split}
\Vert \mathbf{D}\left(\Lambda ^\varepsilon b(\mathbf{D})(\Pi _\varepsilon -I)\boldsymbol{\Phi}\right)\Vert ^2 _{L_2(\mathbb{R}^d)}
&\leqslant \alpha _1 \left(2\mathfrak{c} _1 r_0^{-2}+8\Vert \Lambda \Vert ^2_{L_\infty}(\mathfrak{c} _2+1)\right)\Vert \mathbf{D}^2\boldsymbol{\Phi}\Vert _{L_2(\mathbb{R}^d)}^2.
\end{split}
\end{equation}
Finally, relations~\eqref{1 lemma Lambda (S-I)} and~\eqref{2 lemma Lambda (S-I)} imply the inequality
\begin{equation*}
\begin{split}
&\Vert \Lambda ^\varepsilon b(\mathbf{D})(\Pi _\varepsilon -I)\boldsymbol{\Phi}\Vert _{H^1(\mathbb{R}^d)}
\leqslant
\mathfrak{C}_\Lambda \Vert \mathbf{D}\boldsymbol{\Phi}\Vert _{H^1(\mathbb{R}^d)}
\leqslant \mathfrak{C}_\Lambda \Vert \boldsymbol{\Phi}\Vert _{H^2(\mathbb{R}^d)},
 \\
 &\boldsymbol{\Phi}\in H^2(\mathbb{R}^d;\mathbb{C}^n), \quad \mathfrak{C}_\Lambda ^2
:=\alpha _1\left(2\mathfrak{c} _1 r_0^{-2} +(8\mathfrak{c} _2+12)\Vert \Lambda \Vert ^2 _{L_\infty}\right).
 \end{split}
\end{equation*}
This is equivalent to estimate~\eqref{lemma Lambda (S-I)}.
\end{proof}

\begin{lemma}
\label{Lemma tilde Lambda(S-I)}
Suppose that Condition~\textnormal{\ref{Condition tilde Lambda in Lp}} is satisfied. Let $\Pi_\varepsilon$ be the smoothing operator~\eqref{Pi eps = definition}. 
Then for $0<\varepsilon\leqslant 1$ we have
\begin{equation*}
\Vert [\widetilde{\Lambda}^\varepsilon ](\Pi _\varepsilon -I)\Vert _{H^2(\mathbb{R}^d)\rightarrow H^1(\mathbb{R}^d)}\leqslant \mathfrak{C}_{\widetilde{\Lambda}}.
\end{equation*}
The constant $\mathfrak{C}_{\widetilde{\Lambda}}$ depends only on $n$, $d$, $\alpha _0$, $\alpha _1$, $\rho$, $p$, $\Vert g\Vert _{L_\infty}$, $\Vert g^{-1}\Vert _{L_\infty}$, the norms $\Vert a_j\Vert _{L_\rho (\Omega)}$, $j=1,\dots,d$, and 
$\Vert \widetilde{\Lambda}\Vert _{L_p(\Omega)}$, and the parameters of the lattice $\Gamma$.
\end{lemma}

\begin{proof}
Let $\boldsymbol{\Phi}\in H^2(\mathbb{R}^d;\mathbb{C}^n)$.
From \eqref{Pi eps <= 1} and Lemma~\ref{Lemma_Lambda_tilda3}
it follows that
\begin{equation}
\label{1 lemma tilde Lambda (S-I)}
\Vert \widetilde{\Lambda}^\varepsilon (\Pi _\varepsilon -I)\boldsymbol{\Phi}\Vert _{L_2(\mathbb{R}^d)}
\leqslant 2 C^{(p)}_{\Omega}\Vert \widetilde{\Lambda}\Vert _{L_p(\Omega)}\Vert \boldsymbol{\Phi}\Vert _{H^1(\mathbb{R}^d)}.
\end{equation}

Consider the derivatives:
\begin{equation*}
\partial _j \bigl(\widetilde{\Lambda}^\varepsilon (\Pi _\varepsilon -I)\boldsymbol{\Phi}\bigr)
=\varepsilon ^{-1}(\partial _j \widetilde{\Lambda})^\varepsilon (\Pi _\varepsilon -I)\boldsymbol{\Phi}
+\widetilde{\Lambda}^\varepsilon (\Pi _\varepsilon -I)\partial _j \boldsymbol{\Phi}.
\end{equation*}
Together with Lemmas~\ref{Lemma_Lambda_tilda3} and \ref{Lemma_Lambda_tilda4}, this yields
\begin{equation*}
\begin{split}
\Vert \mathbf{D}\bigl(\widetilde{\Lambda}^\varepsilon (\Pi _\varepsilon -I)\boldsymbol{\Phi}\bigr)\Vert ^2 _{L_2(\mathbb{R}^d)}
&\leqslant 2\widetilde{\mathfrak{c}}_1\varepsilon ^{-2}\Vert (\Pi _\varepsilon -I)\boldsymbol{\Phi} \Vert ^2 _{H^1(\mathbb{R}^d)}
\\
&+2(\widetilde{\mathfrak{c}}_2+1)\Vert \widetilde{\Lambda}\Vert ^2_{L_p(\Omega)}(C^{(p)}_{\Omega})^2\Vert \mathbf{D}(\Pi _\varepsilon -I)\boldsymbol{\Phi}\Vert ^2_{H^1(\mathbb{R}^d)}.
\end{split}
\end{equation*}
Combining this with \eqref{Pi eps <= 1} and Proposition~\ref{Proposition Pi -I}, we obtain  
\begin{equation}
\label{2 lemma tilde Lambda (S-I)}
\Vert \mathbf{D}\bigl(\widetilde{\Lambda}^\varepsilon (\Pi _\varepsilon -I)\boldsymbol{\Phi}\bigr)\Vert ^2 _{L_2(\mathbb{R}^d)}
\leqslant
\bigl(
2\widetilde{\mathfrak{c}}_1r_0^{-2}+8(\widetilde{\mathfrak{c}}_2+1)\Vert \widetilde{\Lambda}\Vert ^2 _{L_p(\Omega)}(C^{(p)}_{\Omega})^2\bigr)
\Vert \mathbf{D}\boldsymbol{\Phi}\Vert ^2 _{H^1(\mathbb{R}^d)}.
\end{equation}
Now, \eqref{1 lemma tilde Lambda (S-I)} and \eqref{2 lemma tilde Lambda (S-I)} imply that
\begin{equation*}
\Vert \widetilde{\Lambda}^\varepsilon (\Pi _\varepsilon -I)\boldsymbol{\Phi}\Vert _{H^1(\mathbb{R}^d)}
\leqslant \mathfrak{C}_{\widetilde{\Lambda}}\Vert \boldsymbol{\Phi}\Vert _{H^2(\mathbb{R}^d)},
\quad \boldsymbol{\Phi}\in H^2(\mathbb{R}^d;\mathbb{C}^n).
\end{equation*}
Here
$\mathfrak{C}_{\widetilde{\Lambda}}^2 :=2\widetilde{\mathfrak{c}}_1r_0^{-2}+(8 \widetilde{\mathfrak{c}}_2 + 12)
(C^{(p)}_{\Omega})^2\Vert \widetilde{\Lambda}\Vert ^2_{L_p(\Omega)}$.
\end{proof}

\subsection{Proof of Theorem~\ref{Theorem removal Pi eps under conditions}}
\label{Proof of theorem corr no Pi}

\begin{proof}[Proof of Theorem~\textnormal{\ref{Theorem removal Pi eps under conditions}}]
By \eqref{B0^2 L2 to H2},
\begin{equation}
\label{st}
\begin{split}
\Vert &\cos (t(B^0)^{1/2})(B^0)^{-1}\Vert _{H^q(\mathbb{R}^d)\rightarrow H^2(\mathbb{R}^d)}
\\
&\leqslant\Vert \cos (t(B^0)^{1/2})(B^0)^{-1}\Vert _{L_2(\mathbb{R}^d)\rightarrow H^2(\mathbb{R}^d)}\leqslant\check{c}^{-1}.
\end{split}
\end{equation}
According to estimate \eqref{L(xi)>=} for the symbol of $B^0$,
\begin{equation}
\label{st st}
\begin{split}
\Vert &(B^0)^{-1/2}\sin (t(B^0)^{1/2})\Vert _{H^{1+q}(\mathbb{R}^d)\rightarrow H^2(\mathbb{R}^d)}
\\
&\leqslant
\Vert (B^0)^{-1/2}\sin (t(B^0)^{1/2})\Vert _{H^1(\mathbb{R}^d)\rightarrow H^2(\mathbb{R}^d)}\leqslant \check{c}^{-1/2}.
\end{split}
\end{equation}
Together with Theorem~\ref{Theorem corrector interpolation}, Lemmas~\ref{Lemma Lambda (S-I)} and \ref{Lemma tilde Lambda(S-I)},  estimates \eqref{st} and \eqref{st st} imply the results of Theorem~\ref{Theorem removal Pi eps under conditions} with the constants ${C}_{19}:=C_{14}+(\mathfrak{C}_\Lambda +\mathfrak{C}_{\widetilde{\Lambda}})\check{c}^{-1}$ and  $C_{20}:=C_{15}+(\mathfrak{C}_\Lambda +\mathfrak{C}_{\widetilde{\Lambda}})\check{c}^{-1/2}$.
\end{proof}

\subsection{Removal of the smoothing operator from the corrector for $3\leqslant d\leqslant 4$}

If $d\leqslant 2$, then, according to Pro\-po\-si\-tions~\ref{Proposition Lambda in L infty <=} and \ref{Proposition tilde Lambda in Lp if}, Theorem~\ref{Theorem removal Pi eps under conditions} is applicable. So, let $d\geqslant 3$. Now we are interested in the possibility to remove the smoothing operator from the corrector without any additional assumptions on the matrix-valued functions $\Lambda$ and $\widetilde{\Lambda}$. 

If $3\leqslant d\leqslant 4$, it turns out that the smoothing operator $\Pi _\varepsilon$ can be eliminated from both terms of the corrector. But now restrictions on $q$ are stronger than ones from Theorem~\ref{Theorem corrector interpolation}. These new restrictions are caused by the multiplier properties of the matrix-valued function $\Lambda ^\varepsilon$.  The following result was obtained in \cite[Lemma 6.3]{MSuAA17}.

\begin{lemma}
\label{Lemma Lambda multiplicator properties}
Let  $\Lambda({\mathbf x})$ be the $\Gamma$-periodic matrix-valued solution of problem 
\textnormal{\eqref{Lambda problem}}. Assume that $d\geqslant 3$ and put $l=d/2$. 

\noindent $1^\circ$. For $0< \varepsilon \leqslant 1$,
the operator $[\Lambda^\varepsilon]$ is a continuous mapping from  $H^{l-1}({\mathbb R}^d;{\mathbb C}^m)$ to $L_2({\mathbb R}^d;{\mathbb C}^n)$ 
and
\begin{equation*}
\|[\Lambda^\varepsilon] \|_{ H^{l-1}(\mathbb{R}^d) \to L_2(\mathbb{R}^d)}
\leqslant C^{(0)} .
\end{equation*}

\noindent $2^\circ$. Let $0< \varepsilon \leqslant 1$. Then for the function
 ${\mathbf u} \in H^{l}({\mathbb R}^d; {\mathbb C}^m)$ we have the inclusion
$\Lambda^\varepsilon {\mathbf u} \in H^1({\mathbb R}^d;{\mathbb C}^n)$ and the estimate
\begin{equation*}
\| \Lambda^\varepsilon {\mathbf u}\|_{H^1({\mathbb R}^d)} \leqslant  C^{(1)} \varepsilon^{-1}
\| {\mathbf u}\|_{L_2({\mathbb R}^d)} + C^{(2)}  \|{\mathbf u}\|_{H^l({\mathbb R}^d)}.
\end{equation*}
The constants $C^{(0)}$, $C^{(1)}$, and $C^{(2)}$ depend on $m$, $d$, $\alpha_0$, $\alpha_1$,
$\|g\|_{L_\infty}$, $\|g^{-1}\|_{L_\infty}$, and the parameters of the lattice $\Gamma$. 
\end{lemma}

\begin{theorem}
\label{Theorem d=3,4}
Suppose that the assumptions of Subsections~\textnormal{\ref{Subsection lattices}--\ref{Subsection Effective operator}} are satisfied. Let $3\leqslant d\leqslant 4$. Let $1/2\leqslant q\leqslant 1$ for $d=3$ and $q=1$ for $d=4$. Let $\mathcal{K}_1^0(\varepsilon ;t)$ and $\mathcal{K}_2^0(\varepsilon ;t)$ be the operators \eqref{K_10(t,eps)}, \eqref{K20(t,eps)}, respectively. Then for $0<\varepsilon\leqslant 1$ and $t\in\mathbb{R}$ we have
\begin{align}
\label{Th cos corr d=3,4}
\begin{split}
\Vert &\cos (tB_\varepsilon ^{1/2})B_\varepsilon ^{-1}-\cos(t (B^0)^{1/2})(B^0)^{-1}-\varepsilon \mathcal{K}_1^0(\varepsilon ;t)\Vert _{H^q(\mathbb{R}^d)\rightarrow H^1(\mathbb{R}^d)}
\\
&\leqslant {C}_{21} \varepsilon ^q (1+\vert t\vert)^q,
\end{split}
\\
\label{Th sin corr d=3,4}
\begin{split}
\Vert &B_\varepsilon ^{-1/2}\sin(t B_\varepsilon ^{1/2})-(B^0)^{-1/2}\sin(t (B^0)^{1/2})-\varepsilon \mathcal{K}_2^0(\varepsilon ;t)\Vert _{H^{1+q}(\mathbb{R}^d)\rightarrow H^1(\mathbb{R}^d)}
\\
&\leqslant {C}_{22}\varepsilon ^q (1+\vert t\vert)^q.
\end{split}
\end{align}
The constants ${C}_{21}$ and ${C}_{22}$ depend only on the problem data \eqref{problem data} and $q$.
\end{theorem}

\begin{proof}
By Proposition~\ref{Proposition Pi -I}, Lemma~\ref{Lemma Lambda multiplicator properties}($2^\circ$), and \eqref{<b^*b<}, \eqref{B0^2 L2 to H2}, \eqref{Pi eps <= 1}, 
\begin{equation}
\label{proof of Th d=3,4 1}
\begin{split}
\Vert &\varepsilon \Lambda ^\varepsilon (\Pi _\varepsilon - I)b(\mathbf{D})\cos (t(B^0)^{1/2})(B^0)^{-1}\Vert _{H^q(\mathbb{R}^d)\rightarrow H^1(\mathbb{R}^d)}
\\
&\leqslant
C^{(1)}\Vert (\Pi_\varepsilon -I) b(\mathbf{D})(B^0)^{-1}\Vert _{H^q(\mathbb{R}^d)\rightarrow L_2(\mathbb{R}^d)}
\\
&+\varepsilon C^{(2)}\Vert (\Pi_\varepsilon -I) b(\mathbf{D})(B^0)^{-1}\Vert _{H^q(\mathbb{R}^d)\rightarrow H^l(\mathbb{R}^d)}
\\
&\leqslant \varepsilon \alpha _1^{1/2}(r_0^{-1}C^{(1)}+2C^{(2)})\check{c}^{-1},\quad l=d/2.
\end{split}
\end{equation}
Similarly, by Proposition~\ref{Proposition tilde Lambda in Lp if}, Lemma~\ref{Lemma tilde Lambda(S-I)}, and \eqref{B0^2 L2 to H2}, 
\begin{equation}
\label{proof of Th d=3,4 2}
\Vert \varepsilon\widetilde{\Lambda}^\varepsilon (\Pi _\varepsilon -I)\cos (t(B^0)^{1/2})(B^0)^{-1}\Vert _{H^q(\mathbb{R}^d)\rightarrow H^1(\mathbb{R}^d)} \leqslant\varepsilon \mathfrak{C}_{\widetilde{\Lambda}}\check{c}^{-1}.
\end{equation}
Combining \eqref{K_1(t,eps)}, \eqref{Th cos corr interpol}, \eqref{K_10(t,eps)}, \eqref{proof of Th d=3,4 1}, and \eqref{proof of Th d=3,4 2}, 
we arrive at estimate \eqref{Th cos corr d=3,4} with the constant ${C}_{21}:=C_{14}+\alpha _1^{1/2}(r_0^{-1}C^{(1)}+2C^{(2)})\check{c}^{-1}+\mathfrak{C}_{\widetilde{\Lambda}}\check{c}^{-1}$.

Now we proceed to the proof of estimate \eqref{Th sin corr d=3,4}. Combining Proposition~\ref{Proposition Pi -I}, Lemma~\ref{Lemma Lambda multiplicator properties}($2^\circ$),  estimate \eqref{L(xi)>=} for the symbol of the operator $B^0$, and \eqref{<b^*b<}, \eqref{Pi eps <= 1}, we get
\begin{equation}
\label{proof of Th d=3,4 3}
\begin{split}
\Vert & \varepsilon \Lambda ^\varepsilon (\Pi _\varepsilon -I)b(\mathbf{D})(B^0)^{-1/2}\sin (t(B^0)^{1/2})\Vert _{H^{1+q}(\mathbb{R}^d)\rightarrow H^1(\mathbb{R}^d)}
\\
&\leqslant \varepsilon \alpha _1^{1/2} (r_0^{-1}C^{(1)}+2C^{(2)})\check{c}^{-1/2}.
\end{split}
\end{equation}
Next, by Proposition~\ref{Proposition tilde Lambda in Lp if}, Lemma~\ref{Lemma tilde Lambda(S-I)}, and  estimate \eqref{L(xi)>=} for the symbol of the operator $B^0$, 
\begin{equation}
\label{proof of Th d=3,4 4}
\Vert \varepsilon\widetilde{\Lambda}^\varepsilon (\Pi _\varepsilon -I)(B^0)^{-1/2}\sin (t(B^0)^{1/2})\Vert _{H^{1+q}(\mathbb{R}^d)\rightarrow H^1(\mathbb{R}^d)} \leqslant\varepsilon \mathfrak{C}_{\widetilde{\Lambda}}\check{c}^{-1/2}.
\end{equation}
Bringing together \eqref{K(t,eps)}, \eqref{Th sin corr interpol}, \eqref{K20(t,eps)}, \eqref{proof of Th d=3,4 3}, and \eqref{proof of Th d=3,4 4}, we arrive at estimate \eqref{Th sin corr d=3,4} with the constant ${C}_{22}:=C_{15}+\alpha _1^{1/2} (r_0^{-1}C^{(1)}+2C^{(2)})\check{c}^{-1/2}+(r_0^{-1}\widetilde{C}^{(1)}+2\widetilde{C}^{(2)})\check{c}^{-1/2}$.
\end{proof}

\section{Homogenization of solutions of hyperbolic systems}

\label{Section Homogenization of solutions of inhomogeneous hyperbolic systems}

It the present section, we apply the results in operator terms to homogenization of solutions of hyperbolic systems.

\subsection{Principal term of approximation}
\label{Subsection solutions principal term}
Let $\mathbf{u}_\varepsilon$ be the generalized solution of the problem
\begin{equation}
\label{4.1a}
\begin{cases}
\partial _t^2 \mathbf{u}_\varepsilon (\mathbf{x},t) =-B_{\varepsilon}\mathbf{u}_\varepsilon (\mathbf{x},t)+\mathbf{F}(\mathbf{x},t) ,\quad \mathbf{x}\in\mathbb{R}^d,\,t\in(0,T),
\\
\mathbf{u}_\varepsilon (\mathbf{x},0)=\boldsymbol{\phi}(\mathbf{x}),\quad (\partial _t\mathbf{u}_\varepsilon)(\mathbf{x},0)=\boldsymbol{\psi}(\mathbf{x}),\quad \mathbf{x}\in\mathbb{R}^d.
\end{cases}
\end{equation}
Here $\boldsymbol{\phi}\in H^r(\mathbb{R}^d;\mathbb{C}^n)$,  $\boldsymbol{\psi}\in H^{r-1}(\mathbb{R}^d;\mathbb{C}^n)$, and $\mathbf{F}\in L_1((0,T);H^{r-1}(\mathbb{R}^d;\mathbb{C}^n))$ for some $0<T\leqslant\infty$ and $0\leqslant r\leqslant 2$. We have
\begin{align}
\label{u-eps with F=}
\begin{split}
\mathbf{u}_\varepsilon (\cdot ,t)
&=\cos (tB_{\varepsilon}^{1/2})\boldsymbol{\phi}
+B_{\varepsilon}^{-1/2}\sin(tB_{\varepsilon}^{1/2})\boldsymbol{\psi}
\\
&+\int _0^t B_\varepsilon ^{-1/2}\sin ((t-s)B_\varepsilon ^{1/2})\mathbf{F}(\cdot ,s)\,ds,
\end{split}
\\
\label{d t u-eps with F=}
\begin{split}
\partial _t \mathbf{u}_\varepsilon (\cdot ,t)&=-\sin (tB_{\varepsilon}^{1/2})B_\varepsilon ^{1/2}\boldsymbol{\phi}
+\cos (tB_{\varepsilon}^{1/2})\boldsymbol{\psi}
\\
&+\int _0^t\cos ((t-s)B_\varepsilon ^{1/2})\mathbf{F}(\cdot ,s)\,ds.
\end{split}
\end{align}

Let $\mathbf{u}_0$ be the solution of the effective problem
\begin{equation*}
\begin{cases}
\partial _t^2 \mathbf{u}_0(\mathbf{x},t) =-B^0\mathbf{u}_0 (\mathbf{x},t)+\mathbf{F}(\mathbf{x},t) ,\quad \mathbf{x}\in\mathbb{R}^d,\,t\in(0,T),
\\
\mathbf{u}_0 (\mathbf{x},0)=\boldsymbol{\phi}(\mathbf{x}),\quad (\partial _t\mathbf{u}_0)(\mathbf{x},0)=\boldsymbol{\psi}(\mathbf{x}),\quad \mathbf{x}\in\mathbb{R}^d.
\end{cases}
\end{equation*}

Combining estimates \eqref{Th cos L2 interpolation}, \eqref{Th dt sin H-1 interpolation}, \eqref{Th dt cos H-1 interpolation}, and \eqref{Th corr int dt} with identities \eqref{u-eps with F=}, \eqref{d t u-eps with F=} and similar representations for $\mathbf{u}_0$ and $\partial _t\mathbf{u}_0$, we arrive at the following result. 

\begin{theorem}
Under the assumptions of Subsections~\textnormal{\ref{Subsection lattices}--\ref{Subsection Effective operator}} and \textnormal{\ref{Subsection solutions principal term}}, for $0<\varepsilon\leqslant 1$ and $t\in (0,T)$, we have
\begin{align*}
\begin{split}
\Vert &\mathbf{u}_\varepsilon (\cdot ,t)-\mathbf{u}_0(\cdot ,t)\Vert _{L_2(\mathbb{R}^d)}\leqslant C_7\varepsilon ^{r/2} (1+\vert t\vert )^{r/2}\Vert \boldsymbol{\phi}\Vert _{H^r(\mathbb{R}^d)}
\\
&+C_{16}\varepsilon ^{r/2}\vert t\vert ^{r/2} (\Vert \boldsymbol{\psi}\Vert _{H^{r-1}(\mathbb{R}^d)}+\Vert \mathbf{F}\Vert _{L_1((0,t);H^{r-1}(\mathbb{R}^d))}),
\end{split}
\\
\begin{split}
\Vert & (\partial _t\mathbf{u}_\varepsilon) (\cdot ,t)-(\partial _t\mathbf{u}_0)(\cdot ,t)\Vert _{H^{-1}(\mathbb{R}^d)}\leqslant C_9\varepsilon ^{r/2} (1+\vert t\vert )^{r/2}\Vert \boldsymbol{\phi}\Vert _{H^r(\mathbb{R}^d)}\\
&+C_{10}\varepsilon ^{r/2}\vert t\vert ^{r/2} (\Vert \boldsymbol{\psi}\Vert _{H^{r-1}(\mathbb{R}^d)}+\Vert \mathbf{F}\Vert _{L_1((0,t);H^{r-1}(\mathbb{R}^d))}).
\end{split}
\end{align*}
The constants $C_7$, $C_9$, $C_{10}$, and $C_{16}$ are controlled explicitly in terms of $r$ and the problem data \eqref{problem data}.
\end{theorem}

\subsection{Approximation with corrector}
\label{Subsection appr with corr for sol with F}

Now, consider the following problem:
\begin{equation}
\label{u_eps tilde problem 4.2a}
\begin{cases}
\partial _t^2 {\mathbf{v}}_\varepsilon  (\mathbf{x},t) =-B_{\varepsilon} {\mathbf{v}}_\varepsilon (\mathbf{x},t)+\mathbf{F}(\mathbf{x},t),\quad\mathbf{x}\in\mathbb{R}^d,\;t\in(0,T),
\\
 {\mathbf{v}}_\varepsilon  (\mathbf{x},0)=B_\varepsilon ^{-1}\boldsymbol{\phi}(\mathbf{x}),\quad (\partial _t {\mathbf{v}}_\varepsilon ) (\mathbf{x},0)=\boldsymbol{\psi}(\mathbf{x}),\quad\mathbf{x}\in\mathbb{R}^d.
\end{cases}
\end{equation}
Here $\boldsymbol{\phi}\in H^q(\mathbb{R}^d;\mathbb{C}^n)$, $\boldsymbol{\psi}\in H^{1+q}(\mathbb{R}^d;\mathbb{C}^n)$, and $\mathbf{F}\in L_1((0,T);H^{1+q}(\mathbb{R}^d;\mathbb{C}^n))$ for some $0<T\leqslant\infty$ and $0\leqslant q\leqslant 1$. We have
\begin{align}
\label{u-eps tilde with F=}
\begin{split}
{\mathbf{v}}_\varepsilon (\cdot ,t)
&=\cos (tB_{\varepsilon}^{1/2})B_\varepsilon ^{-1}\boldsymbol{\phi}
+B_{\varepsilon}^{-1/2}\sin(tB_{\varepsilon}^{1/2})\boldsymbol{\psi}
\\
&+\int _0^t B_\varepsilon ^{-1/2}\sin ((t-s)B_\varepsilon ^{1/2})\mathbf{F}(\cdot ,s)\,ds,
\end{split}
\\
\label{d t u-eps tilde with F=}
\begin{split}
\partial _t  {\mathbf{v}}_\varepsilon (\cdot ,t)&=-\sin (tB_{\varepsilon}^{1/2})B_\varepsilon ^{-1/2}\boldsymbol{\phi}
+\cos (tB_{\varepsilon}^{1/2})\boldsymbol{\psi}
\\
&+\int _0^t\cos ((t-s)B_\varepsilon ^{1/2})\mathbf{F}(\cdot ,s)\,ds.
\end{split}
\end{align}

Let $ {\mathbf{v}}_0$ be the solution of the corresponding effective problem:
\begin{equation}
\label{u_0 tilde problem intr}
\begin{cases}
\partial _t^2  {\mathbf{v}}_0(\mathbf{x},t) =-B^0 {\mathbf{v}}_0(\mathbf{x},t)+\mathbf{F}(\mathbf{x},t),\quad\mathbf{x}\in\mathbb{R}^d,\;t\in(0,T),
\\
{\mathbf{v}}_0 (\mathbf{x},0)=(B^0)^{-1}\boldsymbol{\phi}(\mathbf{x}),\quad (\partial _t {\mathbf{v}}_0)(\mathbf{x},0)=\boldsymbol{\psi}(\mathbf{x}),\quad\mathbf{x}\in\mathbb{R}^d.
\end{cases}
\end{equation}
Then 
\begin{align}
\label{u-0 tilde with F=}
\begin{split}
{\mathbf{v}}_0 (\cdot ,t)
&=\cos (t(B^0)^{1/2})(B^0) ^{-1}\boldsymbol{\phi}
+(B^0)^{-1/2}\sin(t(B^0)^{1/2})\boldsymbol{\psi}
\\
&+\int _0^t (B^0) ^{-1/2}\sin ((t-s)(B^0) ^{1/2})\mathbf{F}(\cdot ,s)\,ds,
\end{split}
\\
\label{d t u-0 tilde with F=}
\begin{split}
\partial _t {\mathbf{v}}_0 (\cdot ,t)&=-\sin (t(B^0)^{1/2})(B^0) ^{-1/2}\boldsymbol{\phi}
+\cos (t(B^0)^{1/2})\boldsymbol{\psi}
\\
&+\int _0^t\cos ((t-s)(B^0) ^{1/2})\mathbf{F}(\cdot ,s)\,ds.
\end{split}
\end{align}

Let ${\mathbf{w}}_\varepsilon $ be the first order approximation for the solution ${\mathbf{v}}_\varepsilon$:
\begin{equation}
\label{tilde v_eps=}
{\mathbf{w}}_\varepsilon (\cdot ,t)={\mathbf{v}}_0(\cdot ,t)+\varepsilon (\Lambda ^\varepsilon \Pi _\varepsilon b(\mathbf{D})+\widetilde{\Lambda}^\varepsilon\Pi _\varepsilon){\mathbf{v}}_0(\cdot ,t).
\end{equation}

From Theorem~\ref{Theorem corrector interpolation}, \eqref{Th cos L2 interpolation} with $r=q+1$, and \eqref{u-eps tilde with F=}, \eqref{d t u-eps tilde with F=}, \eqref{u-0 tilde with F=}--\eqref{tilde v_eps=} we derive the following result.

\begin{theorem}
\label{Theorem solutions H1 interpol}
Under the assumptions of Subsections~\textnormal{\ref{Subsection lattices}--\ref{Subsection Pi_eps}} and \textnormal{\ref{Subsection appr with corr for sol with F}}, let $0\leqslant q\leqslant 1$. Then for $0<\varepsilon\leqslant 1$ and $t\in (0,T)$ we have
\begin{align}
\begin{split}
\Vert &(\partial _t {\mathbf{v}}_\varepsilon )(\cdot ,t)-(\partial _t {\mathbf{v}}_0)(\cdot ,t)\Vert _{L_2(\mathbb{R}^d)}
\leqslant  C_{16}\varepsilon ^{(q+1)/2}\vert t\vert ^{(q+1)/2} \Vert \boldsymbol{\phi}\Vert _{H^q(\mathbb{R}^d)}
\\
&+C_7\varepsilon ^{(q+1)/2} (1+\vert t\vert)^{(q+1)/2}(\Vert \boldsymbol{\psi}\Vert _{H^{1+q}(\mathbb{R}^d)}+\Vert \mathbf{F}\Vert _{L_1((0,t);H^{1+q}(\mathbb{R}^d))}),
\end{split}
\nonumber
\\
\label{appr sol in H1}
\begin{split}
\Vert &{\mathbf{v}}_\varepsilon (\cdot ,t)- {\mathbf{w}}_\varepsilon(\cdot ,t)\Vert _{H^1(\mathbb{R}^d)}
\leqslant C_{14}\varepsilon ^q (1+\vert t\vert )^q\Vert \boldsymbol{\phi}\Vert _{H^q(\mathbb{R}^d)}
\\
&+C_{15}\varepsilon ^q (1+\vert t\vert )^q(\Vert \boldsymbol{\psi}\Vert _{H^{1+q}(\mathbb{R}^d)}+\Vert \mathbf{F}\Vert _{L_1((0,t);H^{1+q}(\mathbb{R}^d))}).
\end{split}
\end{align}
The constants $C_7$, $C_{14}$, $C_{15}$, and $C_{16}$ are controlled explicitly in terms of $q$ and the problem data \eqref{problem data}.
\end{theorem}

According to  Pro\-po\-si\-tions~\textnormal{\ref{Proposition Lambda in L infty <=}} and \textnormal{\ref{Proposition tilde Lambda in Lp if}} and Theorems~\textnormal{\ref{Theorem removal Pi eps under conditions}} and \textnormal{\ref{Theorem d=3,4}}, for $d\leqslant 2$ we can always replace function \eqref{tilde v_eps=} by 
\begin{equation}
\label{check v}
\check{\mathbf{w}}_\varepsilon (\cdot ,t)={\mathbf{v}}_0(\cdot ,t)+\varepsilon (\Lambda ^\varepsilon  b(\mathbf{D})+\widetilde{\Lambda}^\varepsilon ){\mathbf{v}}_0(\cdot ,t)
\end{equation}
in approximation \eqref{appr sol in H1}, i.~e., to remove the smoothing operator from the corrector. This changes only the constants in estimate. For $d=3,4$, we also can remove the smoothing operator, but not for all $q\in (0,1]$.

\begin{theorem}
\label{Theorem solutions d<=4}
Under the assumptions of Subsections~\textnormal{\ref{Subsection lattices}--\ref{Subsection Pi_eps}} and \textnormal{\ref{Subsection appr with corr for sol with F}}, let $d\leqslant 4$. 
Let $0\leqslant q\leqslant 1$ for $d=1,2$\textnormal{;} $1/2\leqslant q\leqslant 1$ for $d=3$\textnormal{;} and let $q=1$ for $d=4$. 
Let $\check{\mathbf{w}}_\varepsilon (\cdot ,t)$ be the function \eqref{check v}. Then for $0<\varepsilon\leqslant 1$ and $t\in\mathbb{R}$ we have
\begin{equation}
\label{Th sol corr d<=4}
\begin{split}
\Vert & {\mathbf{v}}_\varepsilon (\cdot ,t)-\check{\mathbf{w}}_\varepsilon(\cdot ,t)\Vert _{H^1(\mathbb{R}^d)}
\leqslant C_{23}\varepsilon ^q (1+\vert t\vert )^q\Vert \boldsymbol{\phi}\Vert _{H^q(\mathbb{R}^d)}
\\
&+C_{24}\varepsilon ^q (1+\vert t\vert )^q(\Vert \boldsymbol{\psi}\Vert _{H^{1+q}(\mathbb{R}^d)}+\Vert \mathbf{F}\Vert _{L_1((0,t);H^{1+q}(\mathbb{R}^d))}).
\end{split}
\end{equation}
The constants $C_{23}$ and $C_{24}$ depend only on $q$ and the problem data \eqref{problem data}.
\end{theorem}

Under the additional assumption that Conditions~\textnormal{\ref{Condition Lambda in L infty}} and \textnormal{\ref{Condition tilde Lambda in Lp}} hold, we also can remove smoothing operator from the corrector due to Theorem~\ref{Theorem removal Pi eps under conditions}.

\begin{theorem}
\label{Theorem solutions no Pi eps under assumptions}
Let the assumptions of Subsections~\textnormal{\ref{Subsection lattices}--\ref{Subsection Pi_eps}} and \textnormal{\ref{Subsection appr with corr for sol with F}} be satisfied. Assume that Conditions \textnormal{\ref{Condition Lambda in L infty}} and \textnormal{\ref{Condition tilde Lambda in Lp}} hold. Let $\check{\mathbf{w}}_\varepsilon (\cdot ,t)$ be function \eqref{check v}. Then for $0<\varepsilon\leqslant 1$ and $t\in\mathbb{R}$
\begin{equation*}
\begin{split}
\Vert &{\mathbf{v}}_\varepsilon (\cdot ,t)-\check{\mathbf{w}}_\varepsilon(\cdot ,t)\Vert _{H^1(\mathbb{R}^d)}
\leqslant C_{19}\varepsilon ^q (1+\vert t\vert )^q\Vert \boldsymbol{\phi}\Vert _{H^q(\mathbb{R}^d)}
\\
&+C_{20}\varepsilon ^q (1+\vert t\vert )^q(\Vert \boldsymbol{\psi}\Vert _{H^{1+q}(\mathbb{R}^d)}+\Vert \mathbf{F}\Vert _{L_1((0,t);H^{1+q}(\mathbb{R}^d))}).
\end{split}
\end{equation*}
The constants ${C}_{19}$ and ${C}_{20}$ depend only on $q$, $p$, $\Vert \Lambda \Vert _{L_\infty}$, $\Vert \widetilde{\Lambda}\Vert _{L_p(\Omega)}$, and on the problem data \eqref{problem data}.
\end{theorem}

\subsection{Discussion of the results}

\label{Subsection Discussion}

According to the abstract theory of hyperbolic equations (see, e.~g., \cite[\S 8.2.4]{BSol}), the following law of conservation of energy holds for the solution ${\mathbf{v}}_\varepsilon$ of problem \eqref{u_eps tilde problem 4.2a} with $\mathbf{F}=0$:
\begin{equation}
\label{energy eps}
\begin{split}
\Vert &\partial _t {\mathbf{v}}_\varepsilon (\cdot ,t)\Vert ^2 _{L_2(\mathbb{R}^d)}+\Vert B_\varepsilon ^{1/2} {\mathbf{v}}_\varepsilon (\cdot ,t)\Vert ^2 _{L_2(\mathbb{R}^d)}
\\
&=
\Vert \partial _t {\mathbf{v}}_\varepsilon (\cdot ,0)\Vert ^2 _{L_2(\mathbb{R}^d)}+\Vert B_\varepsilon ^{1/2} {\mathbf{v}}_\varepsilon (\cdot ,0)\Vert ^2 _{L_2(\mathbb{R}^d)}
\\
&=\Vert \boldsymbol{\psi}\Vert ^2_{L_2(\mathbb{R}^d)}+\Vert B_\varepsilon ^{-1/2}\boldsymbol{\phi}\Vert ^2_{L_2(\mathbb{R}^d)}.
\end{split}
\end{equation}
For the solution $ {\mathbf{v}}_0$ of problem \eqref{u_0 tilde problem intr} with $\mathbf{F}=0$ we also have the conservation of energy:
\begin{equation}
\label{energy 0}
\Vert \partial _t {\mathbf{v}}_0 (\cdot ,t)\Vert ^2 _{L_2(\mathbb{R}^d)}+\Vert (B^0) ^{1/2} {\mathbf{v}}_ 0 (\cdot ,t)\Vert ^2 _{L_2(\mathbb{R}^d)}
=\Vert \boldsymbol{\psi}\Vert ^2_{L_2(\mathbb{R}^d)}+\Vert (B^0) ^{-1/2}\boldsymbol{\phi}\Vert ^2_{L_2(\mathbb{R}^d)}.
\end{equation}
As was mentioned in Introduction, we consider problem of the form \eqref{u_eps tilde problem 4.2a} instead of \eqref{4.1a} to have convergence of the corresponding energy to the energy of the solution of the effective problem. This convergence allows us to give energy norm approximation for the solution. For $\boldsymbol{\phi}\in H^2(\mathbb{R}^d;\mathbb{C}^n)$, convergence of the energy \eqref{energy eps} to the energy \eqref{energy 0} is a consequence of the following lemma.

\begin{lemma}
\label{Lemma square root homogenization}
Under the assumptions of Subsections~\textnormal{\ref{Subsection lattices}--\ref{Subsection Effective operator}}, for $0<\varepsilon\leqslant 1$ we have
\begin{equation}
\label{lemma B_eps -1/2 appr}
\Vert B_\varepsilon ^{-1/2}-(B^0)^{-1/2}\Vert _{H^2(\mathbb{R}^d)\rightarrow L_2(\mathbb{R}^d)}
\leqslant\mathfrak{C}\varepsilon .
\end{equation}
The constant $\mathfrak{C}$ depends only on the problem data \eqref{problem data}.
\end{lemma}

\begin{proof}
According to \cite[Chapter III, \S 3, Subsection~4]{ViGKo}, we have
\begin{equation*}
B_\varepsilon ^{-1/2}=\frac{1}{\pi}\int _0^\infty \nu ^{-1/2}(B_\varepsilon +\nu I)^{-1}\,d\nu .
\end{equation*}
The similar representation holds for the operator $(B^0)^{-1/2}$. So,
\begin{equation*}
\begin{split}
&(B_\varepsilon ^{-1/2}-(B^0)^{-1/2})(B^0)^{-1}
\\
&=\frac{1}{\pi}\int _0^\infty \nu ^{-1/2}\left( (B_\varepsilon +\nu I)^{-1}-(B^0+\nu I)^{-1}\right)(B^0)^{-1}\,d\nu.
\end{split}
\end{equation*}
By using the resolvent identity, we obtain
\begin{equation}
\label{B-1/2 start of estimating}
\begin{split}
&(B_\varepsilon ^{-1/2}-(B^0)^{-1/2})(B^0)^{-1}\\
&=\frac{1}{\pi}\int _0^\infty \nu ^{-1/2} B_\varepsilon (B_\varepsilon +\nu I)^{-1}(B_\varepsilon ^{-1}-(B^0)^{-1})(B^0+\nu I)^{-1}\,d\nu.
\end{split}
\end{equation}
We have
\begin{equation}
\Vert  B_\varepsilon (B_\varepsilon +\nu I)^{-1}\Vert _{L_2(\mathbb{R}^d)\rightarrow L_2(\mathbb{R}^d)}
\leqslant \sup _{x\geqslant 0}\frac{x}{x+\nu}\leqslant 1,\quad \nu \geqslant 0.
\end{equation}
Next, by \eqref{B0>=},
\begin{equation}
\label{B0+I nu -1}
\Vert (B^0+\nu I)^{-1}\Vert _{L_2(\mathbb{R}^d)\rightarrow L_2(\mathbb{R}^d)}
\leqslant\sup _{x\geqslant \check{c}}\frac{1}{x+\nu}\leqslant \frac{1}{\check{c} +\nu},\quad \nu \geqslant 0.
\end{equation}
Now from Theorem~\ref{Theorem elliptic} and \eqref{B-1/2 start of estimating}--\eqref{B0+I nu -1} it follows that
\begin{equation}
\label{sq root new est}
\begin{split}
\Vert &(B_\varepsilon ^{-1/2}-(B^0)^{-1/2})(B^0)^{-1}\Vert _{L_2(\mathbb{R}^d)\rightarrow L_2(\mathbb{R}^d)}
\leqslant \frac{1}{\pi}C_1\varepsilon \int _0^\infty \nu ^{-1/2}(\check{c}+\nu)^{-1}\,d\nu
\\
&\leqslant \frac{C_1\varepsilon}{\pi}\Bigl(\check{c}^{-1}\int _0^1\nu ^{-1/2}\,d\nu+\int _1^\infty \nu ^{-3/2}\,d\nu\Bigr)
=\frac{2C_1\varepsilon}{\pi}(\check{c}^{-1}+1).
\end{split}
\end{equation}
Together with \eqref{B0 on H2}, this implies estimate \eqref{lemma B_eps -1/2 appr} with the constant $\mathfrak{C}:=2\pi ^{-1}(\check{c}^{-1}+1)C_1C_L$.
\end{proof}

\begin{remark}
Estimate \eqref{lemma B_eps -1/2 appr} does not look optimal with respect to the type of the norm. It is natural to expect that the correct type of the norm is $(H^1\rightarrow L_2)$- one. But to prove such estimate we need to have approximation of the operator $(B_\varepsilon +\nu I)^{-1}$, $\nu\in\mathbb{R}_+$, in $(H^1\rightarrow L_2)$-norm with the error estimate of the form $C\varepsilon(1+\nu )^{-1}$. It is one of the results of the author's work in progress \textnormal{\cite{M+progress}}.
\end{remark}

\subsection{The case where the corrector is equal to zero}

Assume that relations \eqref{overline-g} and \eqref{sum Dj aj =0} hold. Then the corrector is equal to zero (see Subsection~\ref{Subsection K=0 in operator terms}), i. e., $ {\mathbf{v}}_\varepsilon (\cdot ,t)= {\mathbf{v}}_0(\cdot ,t)$. Proposition~\ref{Proposition K=0 for sin and cos} implies the following result.

\begin{proposition}
Let the assumptions of Subsections~\textnormal{\ref{Subsection lattices}--\ref{Subsection Pi_eps}} be satisfied. Suppose that relations \eqref{overline-g} and \eqref{sum Dj aj =0} hold. Let $ {\mathbf{v}}_\varepsilon (\cdot ,t)$ and $ {\mathbf{v}}_0(\cdot ,t)$ be the solutions of problems \eqref{u_eps tilde problem 4.2a} and \eqref{u_0 tilde problem intr}, respectively, where $\boldsymbol{\phi}\in H^{-1+r}(\mathbb{R}^d;\mathbb{C}^n)$, $\boldsymbol{\psi}\in H^{r}(\mathbb{R}^d;\mathbb{C}^n)$, and $\mathbf{F}\in L_1((0,T);H^{r}(\mathbb{R}^d;\mathbb{C}^n))$ for some $0<T\leqslant\infty$ and $0\leqslant r\leqslant 2$. Then for $0<\varepsilon\leqslant 1$ and $t\in\mathbb{R}$ we have
\begin{equation*}
\begin{split}
\Vert &  {\mathbf{v}}_\varepsilon (\cdot ,t)- {\mathbf{v}}_0 (\cdot ,t)\Vert _{H^1(\mathbb{R}^d)} 
\leqslant (c_6^2 + \check{c}^{-1})^{1-r/2} C_{11}^{r/2} \varepsilon ^{r/2} (1+\vert t\vert)^{r/2}\Vert \boldsymbol{\phi}\Vert _{H^{-1+r}(\mathbb{R}^d)}
\\
&+(c_6+\check{c}^{-1/2})^{1-r/2} C_{12}^{r/2}\varepsilon ^{r/2} (1+\vert t\vert)^{r/2}(\Vert \boldsymbol{\psi}\Vert _{H^{r}(\mathbb{R}^d)}+\Vert \mathbf{F}\Vert _{L_1((0,t);H^{r}(\mathbb{R}^d))}).
\end{split}
\end{equation*}
\end{proposition}

\subsection{Approximation of the flux}

\begin{theorem}
\label{Theorem flux}
Let the assumptions of Subsections~\textnormal{\ref{Subsection lattices}--\ref{Subsection Pi_eps}} and \textnormal{\ref{Subsection appr with corr for sol with F}} be satisfied. Let $0\leqslant q\leqslant 1$.  
Let $\widetilde{g}$ be the matrix-valued function \eqref{tilde g}. Then for $0<\varepsilon\leqslant 1$ and $t\in (0,T)$, for the flux $\mathbf{p}_\varepsilon (\cdot ,t):= g^\varepsilon b(\mathbf{D}) {\mathbf{v}}_\varepsilon (\cdot ,t)$ we have an approximation
\begin{equation}
\label{Th flux}
\begin{split}
\Vert &\mathbf{p}_\varepsilon (\cdot ,t)-\widetilde{g}^\varepsilon \Pi _\varepsilon b(\mathbf{D}) {\mathbf{v}}_0(\cdot ,t)-g^\varepsilon (b(\mathbf{D})\widetilde{\Lambda})^\varepsilon \Pi _\varepsilon  {\mathbf{v}}_0(\cdot ,t)\Vert _{L_2(\mathbb{R}^d)}
\\
&\leqslant C_{25}\varepsilon ^q (1+\vert t\vert )^q\Vert \boldsymbol{\phi}\Vert _{H^q(\mathbb{R}^d)}
\\
&+C_{26}\varepsilon ^q (1+\vert t\vert )^q(\Vert \boldsymbol{\psi}\Vert _{H^{1+q}(\mathbb{R}^d)}+\Vert \mathbf{F}\Vert _{L_1((0,t);H^{1+q}(\mathbb{R}^d))}).
\end{split}
\end{equation}
The constants $C_{25}$ and $C_{26}$ are controlled explicitly in terms of $q$ and the problem data \eqref{problem data}.
\end{theorem}

\begin{proof}
By \eqref{<b^*b<} and \eqref{appr sol in H1},
\begin{equation}
\label{proof fluxes start}
\begin{split}
\Vert &g^\varepsilon b(\mathbf{D}) {\mathbf{v}}_\varepsilon (\cdot ,t)- g^\varepsilon b(\mathbf{D}){\mathbf{w}}_\varepsilon (\cdot ,t)
\Vert _{L_2(\mathbb{R}^d)}
\leqslant \alpha _1^{1/2}\Vert g\Vert _{L_\infty}C_{14}\varepsilon ^q (1+\vert t\vert )^q\Vert \boldsymbol{\phi}\Vert _{H^q(\mathbb{R}^d)}
\\
&+\alpha _1^{1/2}\Vert g\Vert _{L_\infty}C_{15}\varepsilon ^q (1+\vert t\vert )^q(\Vert \boldsymbol{\psi}\Vert _{H^{1+q}(\mathbb{R}^d)}+\Vert \mathbf{F}\Vert _{L_1((0,t);H^{1+q}(\mathbb{R}^d))}).
\end{split}
\end{equation}
Combining \eqref{b(D)=} and \eqref{tilde v_eps=}, we obtain
\begin{equation}
\label{tozd g b(D) V-eps}
\begin{split}
g^\varepsilon b(\mathbf{D}){\mathbf{w}}_\varepsilon 
&=
g^\varepsilon b(\mathbf{D}) {\mathbf{v}}_0
+g^\varepsilon (b(\mathbf{D})\Lambda )^\varepsilon \Pi _\varepsilon b(\mathbf{D}) {\mathbf{v}}_0
\\
&+g^\varepsilon (b(\mathbf{D})\widetilde{\Lambda})^\varepsilon \Pi _\varepsilon  {\mathbf{v}}_0
+\varepsilon\sum _{j=1}^d g^\varepsilon b_j (\Lambda ^\varepsilon\Pi _\varepsilon b(\mathbf{D})+\widetilde{\Lambda}^\varepsilon \Pi _\varepsilon ) D_j  {\mathbf{v}}_0.
\end{split}
\end{equation}
The fourth term in the right-hand side of \eqref{tozd g b(D) V-eps} can be estimated with the help of \eqref{b_l <=} and \eqref{partof corr H1 to L2}:
\begin{equation}
\begin{split}
\Vert \varepsilon\sum _{j=1}^d g^\varepsilon b_j (\Lambda ^\varepsilon\Pi _\varepsilon b(\mathbf{D})+\widetilde{\Lambda}^\varepsilon \Pi _\varepsilon) D_j  {\mathbf{v}}_0(\cdot ,t)\Vert _{L_2(\mathbb{R}^d)}
\leqslant \varepsilon {C}_{27}
\Vert  {\mathbf{v}}_0 (\cdot ,t)\Vert _{H^2(\mathbb{R}^d)},
\end{split}
\end{equation}
where ${C}_{27}:=\Vert g\Vert _{L_\infty}(d\alpha _1)^{1/2}(M_1\alpha _1^{1/2}+\widetilde{M}_1)$. 
From \eqref{u-0 tilde with F=} and the estimate \eqref{L(xi)>=} for the symbol of the operator $B^0$ it follows that
\begin{equation}
\label{norm tilde u0 in H1}
\begin{split}
\Vert  {\mathbf{v}}_0 (\cdot ,t)\Vert _{H^2(\mathbb{R}^d)}
&\leqslant 
\Vert  {\mathbf{v}}_0 (\cdot ,t)\Vert _{H^{2+q}(\mathbb{R}^d)}\leqslant
\Vert (B^0)^{-1}\boldsymbol{\phi}\Vert _{H^{2+q}(\mathbb{R}^d)}
\\
&+\Vert (B^0)^{-1/2}\boldsymbol{\psi}\Vert _{H^{2+q}(\mathbb{R}^d)}
+\int _0^t\Vert (B^0)^{-1/2}\mathbf{F}(\cdot ,s)\Vert _{H^{2+q}(\mathbb{R}^d)}\,ds
\\
&\leqslant \check{c}^{-1}\Vert \boldsymbol{\phi}\Vert _{H^q(\mathbb{R}^d)}
+\check{c}^{-1/2}(\Vert \boldsymbol{\psi}\Vert _{H^{1+q}(\mathbb{R}^d)}+\Vert \mathbf{F}\Vert _{L_1((0,t);H^{1+q}(\mathbb{R}^d))}),
\end{split}
\end{equation}
for any $0\leqslant q\leqslant 1$. 
Next, by Proposition~\ref{Proposition Pi -I} and \eqref{<b^*b<},
\begin{equation}
\label{proof fluxes final}
\Vert g^\varepsilon b(\mathbf{D}) {\mathbf{v}}_0(\cdot ,t)-g^\varepsilon \Pi _\varepsilon b(\mathbf{D}) {\mathbf{v}}_0(\cdot ,t)
\Vert _{L_2(\mathbb{R}^d)}
\leqslant \varepsilon r_0^{-1}\alpha _1^{1/2}\Vert g\Vert _{L_\infty}\Vert \mathbf{D}^2 {\mathbf{v}}_0(\cdot ,t)\Vert _{L_2(\mathbb{R}^d)}.
\end{equation}
Now from \eqref{tilde g}, and \eqref{proof fluxes start}--\eqref{proof fluxes final} we derive required estimate \eqref{Th flux} with the~constants
\begin{align*}
C_{25}&:=\alpha _1^{1/2}\Vert g\Vert _{L_\infty}C_{14}+C_{27}\check{c}^{-1}+r_0^{-1}\alpha _1^{1/2}\Vert g\Vert _{L_\infty}\check{c}^{-1},
\\
C_{26}&:=\alpha _1^{1/2}\Vert g\Vert _{L_\infty}C_{15}+C_{27}\check{c}^{-1/2}+r_0^{-1}\alpha _1^{1/2}\Vert g\Vert _{L_\infty}\check{c}^{-1/2}.
\end{align*}
\end{proof}

\subsection{On the possibility to remove $\Pi_\varepsilon$ from approximation of the flux}

If $d\leqslant 4$, we can derive approximation of the flux $\mathbf{p}_\varepsilon$ from Theorem~\ref{Theorem solutions d<=4}. The proof repeats the proof of Theorem~\ref{Theorem flux} with some simplifications.

\begin{theorem}
\label{Theorem flux d<=4}
Under the assumptions of Theorem~\textnormal{\ref{Theorem flux}}, let $d\leqslant 4$. 
If $d=3$, we additionally assume that $1/2\leqslant q\leqslant 1$\textnormal{;} and, if $d=4$, we suppose that $q=1$. 
Then for $0<\varepsilon\leqslant 1$ and $t\in\mathbb{R}$, we have
\begin{equation}
\label{Th flux d<=4}
\begin{split}
\Vert &\mathbf{p}_\varepsilon (\cdot ,t)-\widetilde{g}^\varepsilon  b(\mathbf{D}) {\mathbf{v}}_0(\cdot ,t)-g^\varepsilon (b(\mathbf{D})\widetilde{\Lambda})^\varepsilon   {\mathbf{v}}_0(\cdot ,t)\Vert _{L_2(\mathbb{R}^d)}
\\
&\leqslant C_{28}\varepsilon ^q (1+\vert t\vert )^q\Vert \boldsymbol{\phi}\Vert _{H^q(\mathbb{R}^d)}
\\
&+C_{29}\varepsilon ^q (1+\vert t\vert )^q(\Vert \boldsymbol{\psi}\Vert _{H^{1+q}(\mathbb{R}^d)}+\Vert \mathbf{F}\Vert _{L_1((0,t);H^{1+q}(\mathbb{R}^d))}).
\end{split}
\end{equation}
The constants $C_{28}$ and $C_{29}$ depend on $q$ and the problem data \eqref{problem data}.
\end{theorem}

\begin{proof}
By \eqref{<b^*b<} and \eqref{Th sol corr d<=4},
\begin{equation}
\label{proof fluxes d<=4 start}
\begin{split}
\Vert &g^\varepsilon b(\mathbf{D}) {\mathbf{v}}_\varepsilon (\cdot ,t)- g^\varepsilon b(\mathbf{D})\check{\mathbf{w}}_\varepsilon (\cdot ,t)
\Vert _{L_2(\mathbb{R}^d)}
\leqslant \alpha _1^{1/2}\Vert g\Vert _{L_\infty}C_{23}\varepsilon ^q (1+\vert t\vert )^q\Vert \boldsymbol{\phi}\Vert _{H^q(\mathbb{R}^d)}
\\
&+\alpha _1^{1/2}\Vert g\Vert _{L_\infty}C_{24}\varepsilon ^q (1+\vert t\vert )^q(\Vert \boldsymbol{\psi}\Vert _{H^{1+q}(\mathbb{R}^d)}+\Vert \mathbf{F}\Vert _{L_1((0,t);H^{1+q}(\mathbb{R}^d))}).
\end{split}
\end{equation}
From \eqref{b(D)=} and \eqref{check v} it follows that
\begin{equation}
\label{tozd g b(D) V-eps d<=4}
\begin{split}
g^\varepsilon b(\mathbf{D})\check{\mathbf{w}}_\varepsilon 
&=
g^\varepsilon b(\mathbf{D}) {\mathbf{v}}_0
+g^\varepsilon (b(\mathbf{D})\Lambda )^\varepsilon  b(\mathbf{D}) {\mathbf{v}}_0
\\
&+g^\varepsilon (b(\mathbf{D})\widetilde{\Lambda})^\varepsilon  {\mathbf{v}}_0
+\varepsilon\sum _{j=1}^d g^\varepsilon b_j (\Lambda ^\varepsilon b(\mathbf{D})+\widetilde{\Lambda}^\varepsilon ) D_j  {\mathbf{v}}_0.
\end{split}
\end{equation}
Let us estimate the last summand in the right-hand side of \eqref{tozd g b(D) V-eps d<=4}. By \eqref{<b^*b<}, \eqref{b_l <=}, Propositions \ref{Proposition Lambda in L infty <=} and \ref{Proposition tilde Lambda in Lp if}, Lemmas~\ref{Lemma_Lambda_tilda3} and \ref{Lemma Lambda multiplicator properties}($1^\circ$),
\begin{equation}
\label{star or 4.29a}
\begin{split}
\Vert & \varepsilon\sum _{j=1}^d g^\varepsilon b_j (\Lambda ^\varepsilon b(\mathbf{D})+\widetilde{\Lambda}^\varepsilon ) D_j  {\mathbf{v}}_0(\cdot ,t)\Vert _{L_2(\mathbb{R}^d)}
\\
&\leqslant (d\alpha _1)^{1/2}\Vert g\Vert _{L_\infty}(C^{(0)}\alpha _1^{1/2}\Vert   {\mathbf{v}}_0 (\cdot ,t)\Vert _{H^{l+1}(\mathbb{R}^d)}+\Vert \widetilde{\Lambda}\Vert _{L_p(\Omega )}C_\Omega ^{(p)}\Vert   {\mathbf{v}}_0 (\cdot ,t)\Vert _{H^2(\mathbb{R}^d)}).
\end{split}
\end{equation}
(Here $C^{(0)}:=\Vert \Lambda\Vert _{L_\infty}$ and $l=1$ for $d=1,2$.) Note that $\Vert  {\mathbf{v}}_0 (\cdot ,t)\Vert _{H^{l+1}(\mathbb{R}^d)}\leqslant \Vert  {\mathbf{v}}_0 (\cdot ,t)\Vert _{H^{2+q}(\mathbb{R}^d)}$ for $d=3,4$. So, 
together with \eqref{tilde g}, \eqref{norm tilde u0 in H1}, and \eqref{proof fluxes d<=4 start}, estimate \eqref{star or 4.29a} implies inequality \eqref{Th flux d<=4} with the constants 
\begin{align*}
&C_{28}:=(d\alpha _1)^{1/2}\Vert g\Vert _{L_\infty}(C_{23}+C^{(0)}\alpha _1^{1/2}\check{c}^{-1}+\Vert \widetilde{\Lambda}\Vert _{L_p(\Omega )}C_\Omega ^{(p)}\check{c}^{-1}),
\\
&C_{29}:=(d\alpha _1)^{1/2}\Vert g\Vert _{L_\infty}(C_{24}+C^{(0)}\alpha _1^{1/2}\check{c}^{-1/2}+\Vert \widetilde{\Lambda}\Vert _{L_p(\Omega )}C_\Omega ^{(p)}\check{c}^{-1/2}).
\end{align*}
\end{proof}

If $d\geqslant 5$, under Conditions~\ref{Condition Lambda in L infty} and \ref{Condition tilde Lambda in Lp}, we can remove the smoothing operator $\Pi_\varepsilon$ from approximation of the flux. The proof of the following result is based on Theorem~\ref{Theorem solutions no Pi eps under assumptions} and is quite similar to that of Theorem~\ref{Theorem flux d<=4}. We omit the details.

\begin{theorem}
\label{Theorem flux under condidtions}
Under the assumptions of Theorem~\textnormal{\ref{Theorem flux}}, let  Conditions~\textnormal{\ref{Condition Lambda in L infty}} and \textnormal{\ref{Condition tilde Lambda in Lp}} hold. Then for $0<\varepsilon\leqslant 1$ and $t\in\mathbb{R}$, we have
\begin{equation*}
\begin{split}
\Vert &\mathbf{p}_\varepsilon (\cdot ,t)-\widetilde{g}^\varepsilon  b(\mathbf{D}) {\mathbf{v}}_0(\cdot ,t)-g^\varepsilon (b(\mathbf{D})\widetilde{\Lambda})^\varepsilon   {\mathbf{v}}_0(\cdot ,t)\Vert _{L_2(\mathbb{R}^d)}
\\
&\leqslant C_{30}\varepsilon ^q (1+\vert t\vert )^q\Vert \boldsymbol{\phi}\Vert _{H^q(\mathbb{R}^d)}
\\
&+C_{31}\varepsilon ^q (1+\vert t\vert )^q(\Vert \boldsymbol{\psi}\Vert _{H^{1+q}(\mathbb{R}^d)}+\Vert \mathbf{F}\Vert _{L_1((0,t);H^{1+q}(\mathbb{R}^d))}).
\end{split}
\end{equation*}
The constants 
\begin{align*}
&C_{30}:=(d\alpha _1)^{1/2}\Vert g\Vert _{L_\infty}(C_{19}+\Vert \Lambda\Vert _{L_\infty}\alpha _1^{1/2}\check{c}^{-1}+\Vert \widetilde{\Lambda}\Vert _{L_p(\Omega )}C_\Omega ^{(p)}\check{c}^{-1}),
\\
&C_{31}:=(d\alpha _1)^{1/2}\Vert g\Vert _{L_\infty}(C_{20}+\Vert \Lambda\Vert _{L_\infty}\alpha _1^{1/2}\check{c}^{-1/2}+\Vert \widetilde{\Lambda}\Vert _{L_p(\Omega )}C_\Omega ^{(p)}\check{c}^{-1/2}),
\end{align*}
 depend on $q$, $\Vert \Lambda\Vert _{L_\infty}$, $\Vert \widetilde{\Lambda}\Vert _{L_p(\Omega)}$, $p$, and the problem data \eqref{problem data}.
\end{theorem}

\subsection{The special case} 
Assume that $g^0=\underline{g}$, i.~e., relations \eqref{underline-g} are satisfied. Then, by Proposition~\ref{Proposition Lambda in L infty <=}($3^\circ$), Condition~\ref{Condition Lambda in L infty} holds. Herewith, according to \cite[Remark 3.5]{BSu05}, the matrix-valued function \eqref{tilde g} is constant and coincides with $g^0$, i.~e., $\widetilde{g}(\mathbf{x})=g^0=\underline{g}$. Thus, $\widetilde{g}^\varepsilon b(\mathbf{D}) {\mathbf{v}}_0(\cdot,t)=g^0b(\mathbf{D}) {\mathbf{v}}_0(\cdot,t)$.

Assume that condition \eqref{sum Dj aj =0} holds. 
Then the $\Gamma$-periodic solution of problem \eqref{tildeLambda_problem} also equals to zero:  $\widetilde{\Lambda}(\mathbf{x})=0$. 
 So, Theorem~\ref{Theorem flux under condidtions} implies the following result.

\begin{proposition}
Under the assumptions of Theorem~\textnormal{\ref{Theorem flux}}, suppose that relations \eqref{underline-g} and \eqref{sum Dj aj =0} hold. Then for $t\in\mathbb{R}$ and $0<\varepsilon \leqslant 1$ we have
\begin{equation*}
\begin{split}
\Vert &\mathbf{p}_\varepsilon (\cdot ,t)-{g}^0 b(\mathbf{D}) {\mathbf{v}}_0(\cdot ,t)\Vert _{L_2(\mathbb{R}^d)}
\leqslant C_{30}\varepsilon ^q (1+\vert t\vert )^q\Vert \boldsymbol{\phi}\Vert _{H^q(\mathbb{R}^d)}
\\
&+C_{31}\varepsilon ^q (1+\vert t\vert )^q(\Vert \boldsymbol{\psi}\Vert _{H^{1+q}(\mathbb{R}^d)}+\Vert \mathbf{F}\Vert _{L_1((0,t);H^{1+q}(\mathbb{R}^d))}).
\end{split}
\end{equation*}
\end{proposition}

\section{On approximation of the operator $e^{-itB_\varepsilon}$}

\label{Section Schrodinger}

\subsection{Principal term of approximation}

In the present subsection we give an alternative proof of Theorem~3.1.1 from \cite{D}.

\begin{theorem}[\cite{D}]
Let the assumptions of Subsections~\textnormal{\ref{Subsection lattices}--\ref{Subsection Effective operator}} be satisfied. Let $0\leqslant r\leqslant 3$. Then for $0<\varepsilon\leqslant 1$ and $t\in\mathbb{R}$ we have
\begin{equation}
\label{Th by Dorodnyi}
\Vert e^{-itB_\varepsilon}-e^{-itB^0}\Vert _{H^r(\mathbb{R}^d)\rightarrow L_2(\mathbb{R}^d)}
\leqslant C_{32} \varepsilon ^{r/3}(1+\vert t\vert)^{r/3}.
\end{equation}
The constant $C_{32}$ depends only on $r$ and the problem data \eqref{problem data}.
\end{theorem}

\begin{proof}
Denote $\Upsilon(t):=e^{-itB_\varepsilon}B_\varepsilon ^{-1}(B^0)^{-1}e^{itB^0}$. Then, according to \eqref{R1}, 
\begin{equation*}
\frac{d\Upsilon (t)}{dt}=ie^{-it B_\varepsilon}\mathcal{R}_1(\varepsilon)e^{itB^0}.
\end{equation*}
Thus,
\begin{equation*}
\begin{split}
\Upsilon (t)-\Upsilon (0)
&=e^{-it B_\varepsilon}B_\varepsilon ^{-1}(B^0)^{-1}e^{it B^0}
-B_\varepsilon ^{-1}(B^0)^{-1}
\\
&=
i\int _0^t e^{-isB_\varepsilon}\mathcal{R}_1(\varepsilon)e^{isB^0}\,ds.
\end{split}
\end{equation*}
Multiplying this equality by $e^{-itB^0}$ from the right, we obtain
\begin{equation*}
\begin{split}
&e^{-it B_\varepsilon}B_\varepsilon ^{-1}(B^0)^{-1}-B_\varepsilon ^{-1}(B^0)^{-1}e^{-it B^0}
\\
&=
i\int _0^t e^{-isB_\varepsilon}\mathcal{R}_1(\varepsilon)e^{i(s-t)B^0}\,ds.
\end{split}
\end{equation*}
The ranges of all operators in this identity lie in $H^1(\mathbb{R}^d;\mathbb{C}^n)$, so we can multiply it by $B_\varepsilon ^{1/2}$ from the left:
\begin{equation*}
\begin{split}
&e^{-it B_\varepsilon}B_\varepsilon ^{-1/2}(B^0)^{-1}-B_\varepsilon ^{-1/2}(B^0)^{-1}e^{-it B^0}
\\
&=
i\int _0^t e^{-isB_\varepsilon}B_\varepsilon ^{1/2}\mathcal{R}_1(\varepsilon)e^{i(s-t)B^0}\,ds.
\end{split}
\end{equation*}
Thus,
\begin{equation}
\label{imaginary exp identity}
\begin{split}
(&e^{-itB_\varepsilon}-e^{-itB^0})(B^0)^{-3/2}
=
-e^{-itB_\varepsilon}(B_\varepsilon ^{-1/2}-(B^0)^{-1/2})(B^0)^{-1}
\\
&
+(B_\varepsilon ^{-1/2}-(B^0)^{-1/2})(B^0)^{-1}e^{-itB^0}
+i\int _0^t e^{-isB_\varepsilon}B_\varepsilon ^{1/2}\mathcal{R}_1(\varepsilon)e^{i(s-t)B^0}\,ds.
\end{split}
\end{equation}

Let us consider the last summand in the right-hand side of \eqref{imaginary exp identity}. According to \eqref{R1} and \eqref{R2},
\begin{equation}
\label{int with R1 for exp =}
\begin{split}
i\int _0^t e^{-isB_\varepsilon}B_\varepsilon ^{1/2}\mathcal{R}_1(\varepsilon)e^{i(s-t)B^0}\,ds
&=
i\int _0^t e^{-isB_\varepsilon}B_\varepsilon ^{1/2}\mathcal{R}_2(\varepsilon)e^{i(s-t)B^0}\,ds
\\
&+i\int _0^t e^{-isB_\varepsilon}B_\varepsilon ^{1/2}\varepsilon K(\varepsilon)
e^{i(s-t)B^0}\,ds.
\end{split}
\end{equation}
By \eqref{Beps1/2 appr corr},
\begin{equation}
\label{5.5 new}
\Bigl\Vert \int _0^t e^{-isB_\varepsilon}B_\varepsilon ^{1/2}\mathcal{R}_2(\varepsilon)e^{i(s-t)B^0}\,ds\Bigr\Vert _{L_2(\mathbb{R}^d)\rightarrow L_2(\mathbb{R}^d)}
\leqslant c_5C_2\varepsilon\vert t\vert.
\end{equation}
Integrating by parts and taking \eqref{K(eps)} into account, we rewrite the last summand in the right-hand side of \eqref{int with R1 for exp =} as follows:
\begin{equation}
\label{last summand= exp }
\begin{split}
&i\int _0^t e^{-isB_\varepsilon}B_\varepsilon ^{1/2}\varepsilon K(\varepsilon)
e^{i(s-t)B^0}\,ds
=-\int _0^t \frac{d e^{-isB_\varepsilon}}{ds}B_\varepsilon ^{-1/2}\varepsilon K(\varepsilon)
e^{i(s-t)B^0}\,ds
\\
&=
-e^{-it B_\varepsilon }B_\varepsilon ^{-1/2}\varepsilon K(\varepsilon)
+B_\varepsilon ^{-1/2}\varepsilon K(\varepsilon )e^{-itB^0}
\\
&+i\int _0^t e^{-isB_\varepsilon}B_\varepsilon ^{-1/2}\varepsilon (\Lambda ^\varepsilon \Pi_\varepsilon b(\mathbf{D})+\widetilde{\Lambda}^\varepsilon \Pi _\varepsilon )
e^{i(s-t)B^0}\,ds .
\end{split}
\end{equation}
 
By \eqref{K(eps)<= L2 to L2} and \eqref{B eps -1/2 est},
\begin{equation}
\label{5.8 new}
\Vert -e^{-it B_\varepsilon }B_\varepsilon ^{-1/2}\varepsilon K(\varepsilon)
+B_\varepsilon ^{-1/2}\varepsilon K(\varepsilon )e^{-itB^0}\Vert _{L_2(\mathbb{R}^d)\rightarrow L_2(\mathbb{R}^d)}\leqslant 2\varepsilon\beta ^{-1/2}C_K.
\end{equation}

We proceed to estimation of the third summand in the right-hand side of \eqref{last summand= exp }. By \eqref{<b^*b<},
\begin{equation}
\label{b(D) to H-1}
\begin{split}
\Vert b(\mathbf{D})\Vert _{L_2(\mathbb{R}^d)\rightarrow H^{-1}(\mathbb{R}^d)}
&\leqslant \sup _{\boldsymbol{\xi}\in\mathbb{R}^d}(1+\vert \boldsymbol{\xi}\vert ^2)^{-1/2}\vert b(\boldsymbol{\xi})\vert 
\\
&\leqslant \alpha _1^{1/2}\sup _{\boldsymbol{\xi}\in\mathbb{R}^d}(1+\vert \boldsymbol{\xi}\vert ^2)^{-1/2}\vert \boldsymbol{\xi}\vert \leqslant \alpha _1^{1/2}.
\end{split}
\end{equation}
By \eqref{b_eps=> H1-norm} for a function $\mathbf{u}=B_\varepsilon ^{-1/2}\mathbf{v}$, $\mathbf{v}\in L_2(\mathbb{R}^d;\mathbb{C}^n)$,
\begin{equation*}
\Vert (\mathbf{D}^2+I)^{1/2}B_\varepsilon ^{-1/2}\mathbf{v}\Vert _{L_2(\mathbb{R}^d)}=\Vert B_\varepsilon ^{-1/2}\mathbf{v}\Vert _{H^1(\mathbb{R}^d)}
\leqslant c_6\Vert \mathbf{v}\Vert _{L_2(\mathbb{R}^d)}.
\end{equation*}
So,
\begin{equation*}
\Vert B_\varepsilon ^{-1/2}\Vert _{L_2(\mathbb{R}^d)\rightarrow H^1(\mathbb{R}^d)}
=\Vert (\mathbf{D}^2+I)^{1/2}B_\varepsilon ^{-1/2}\Vert  _{L_2(\mathbb{R}^d)\rightarrow L_2(\mathbb{R}^d)}\leqslant c_6.
\end{equation*}
By the duality arguments,
\begin{equation}
\label{B eps -1/2 from H-1 to L2}
\Vert B_\varepsilon ^{-1/2}\Vert _{H^{-1}(\mathbb{R}^d)\rightarrow L_2(\mathbb{R}^d)}\leqslant c_6.
\end{equation}
Now from Proposition~\ref{Proposition f Pi on H-kappa} and \eqref{Lambda <=}, \eqref{b(D) to H-1}, and \eqref{B eps -1/2 from H-1 to L2} it follows that
\begin{equation}
\label{first part in int corr exp}
\begin{split}
&\Vert B_\varepsilon ^{-1/2}\varepsilon \Lambda ^\varepsilon \Pi _\varepsilon b(\mathbf{D})\Vert _{L_2(\mathbb{R}^d)\rightarrow L_2(\mathbb{R}^d)}
\\
&\leqslant \varepsilon \Vert B_\varepsilon ^{-1/2}\Vert _{H^{-1}(\mathbb{R}^d)\rightarrow L_2(\mathbb{R}^d)}
\Vert [\Lambda ^\varepsilon]\Pi _\varepsilon\Vert _{H^{-1}(\mathbb{R}^d)\rightarrow H^{-1}(\mathbb{R}^d)}\Vert b(\mathbf{D})\Vert _{L_2(\mathbb{R}^d)\rightarrow H^{-1}(\mathbb{R}^d)}
\\
&\leqslant \varepsilon c_6 \alpha _1^{1/2}M_1.
\end{split}
\end{equation}
By Proposition~\ref{Proposition f Pi on H-kappa} and \eqref{tilde Lambda<=}, \eqref{B eps -1/2 est},
\begin{equation}
\label{second part in int corr exp}
\Vert B_\varepsilon ^{-1/2}\varepsilon [ \widetilde{\Lambda}^\varepsilon ]\Pi _\varepsilon\Vert _{L_2(\mathbb{R}^d)\rightarrow L_2(\mathbb{R}^d)}\leqslant \varepsilon
\beta ^{-1/2}\widetilde{M}_1.
\end{equation}
From \eqref{first part in int corr exp} and \eqref{second part in int corr exp} it follows that
\begin{equation}
\label{imaginary exp identity last summand estimate}
\begin{split}
\Bigl\Vert i\int _0^t e^{-isB_\varepsilon}B_\varepsilon ^{-1/2}\varepsilon (\Lambda ^\varepsilon \Pi_\varepsilon b(\mathbf{D})+\widetilde{\Lambda}^\varepsilon \Pi _\varepsilon )
e^{i(s-t)B^0}\,ds\Bigr\Vert _{L_2(\mathbb{R}^d)\rightarrow L_2(\mathbb{R}^d)}
\\
\leqslant (c_6 \alpha _1^{1/2}M_1+\beta ^{-1/2}\widetilde{M}_1)\varepsilon\vert t\vert .
\end{split}
\end{equation}

It remains to estimate the first two summands in the right-hand side of \eqref{imaginary exp identity}. 
Using  \eqref{sq root new est}, \eqref{imaginary exp identity}--\eqref{5.8 new}, and \eqref{imaginary exp identity last summand estimate}, we arrive at
\begin{equation*}
\begin{split}
\Vert (e^{-itB_\varepsilon}-e^{-itB^0})(B^0)^{-3/2}\Vert _{L_2(\mathbb{R}^d)\rightarrow L_2(\mathbb{R}^d)}
\leqslant {C}_{33}\varepsilon(1+\vert t\vert),
\end{split}
\end{equation*}
where 
$${C}_{33}:=\max\lbrace 4\pi ^{-1}(\check{c}^{-1}+1)C_1+2\beta ^{-1/2}C_K ;c_5C_2+c_6 \alpha _1^{1/2}M_1+\beta ^{-1/2}\widetilde{M}_1 \rbrace .$$
Together with \eqref{B0 3/2 on H3} this implies that
\begin{equation}
\label{im exp new est}
\begin{split}
\Vert e^{-itB_\varepsilon}-e^{-itB^0}\Vert _{H^3(\mathbb{R}^d)\rightarrow L_2(\mathbb{R}^d)}
\leqslant {C}_{33}C_L^{3/2}\varepsilon(1+\vert t\vert).
\end{split}
\end{equation}
Interpolating between the rough estimate
\begin{equation*}
\Vert e^{-itB_\varepsilon}-e^{-itB^0}\Vert _{L_2(\mathbb{R}^d)\rightarrow L_2(\mathbb{R}^d)}\leqslant 2
\end{equation*}
and  \eqref{im exp new est}, we derive the required estimate \eqref{Th by Dorodnyi} with the constant $C_{32}:=2^{1-r/3}{C}_{33}^{r/3}C_L^{r/2}$.
\end{proof}

\subsection{Approximation of the operator $e^{-itB_\varepsilon}B_\varepsilon ^{-1}$ in the energy norm}

As for hyperbolic problems, we can give approximation for the solution of the non-stationary Schr\"odinger equation  in the energy norm only for the very specific choice of the initial data. In operator terms, we deal with the operator $e^{-itB_\varepsilon}B_\varepsilon ^{-1}$. We give the corresponding approximation for completeness of the  presentation. 

\begin{theorem} 
\label{Theorem exp energy norm}
Let the assumptions of Subsections~\textnormal{\ref{Subsection lattices}--\ref{Subsection Effective operator}} be satisfied. 
Denote 
\begin{equation}
\label{K3(eps,t)}
\mathcal{K}_3(\varepsilon ;t):=(\Lambda ^\varepsilon \Pi _\varepsilon b(\mathbf{D})+\widetilde{\Lambda}^\varepsilon  \Pi _\varepsilon )e^{-itB^0}(B^0)^{-1}.
\end{equation}
Let $0\leqslant r\leqslant 2$. Then for $0<\varepsilon\leqslant 1$ and $t\in\mathbb{R}$ we have
\begin{equation}
\label{Th exp corr}
\begin{split}
\Vert  e^{-itB_\varepsilon}B_\varepsilon ^{-1}-e^{-it B^0}(B^0)^{-1}-\varepsilon\mathcal{K}_3(\varepsilon ;t)\Vert _{H^r(\mathbb{R}^d)\rightarrow H^1(\mathbb{R}^d)}
\leqslant C_{34}\varepsilon ^{r/2} (1+\vert t\vert)^{r/2}.
\end{split}
\end{equation}
Here the constant $C_{34}$ depends only on $r$ and the problem data \eqref{problem data}.
\end{theorem}

\begin{remark}
Lemma~\textnormal{\ref{Lemma square root homogenization}} together with Theorem~\textnormal{\ref{Theorem exp energy norm}} allows us to obtain  approximation for the operator $ e^{-itB_\varepsilon}B_\varepsilon ^{-1/2}$ in the $(H^{3}\rightarrow H^1)$-norm with the error estimate of the form $C\varepsilon (1+\vert t\vert)$.
\end{remark}

\begin{remark}
Using properties of the matrix-valued functions $\Lambda$ and $\widetilde{\Lambda}$ \textnormal{(}see Subsection~\textnormal{\ref{Subsection properties of Lambda and tilde Lambda}}, Lemma~\textnormal{\ref{Lemma Lambda multiplicator properties}}, and Lemma \textnormal{6.5} from \textnormal{\cite{MSuAA17}),} it is always possible to remove $\Pi _\varepsilon$ from the corrector \eqref{K3(eps,t)} for $d=1,2$. If $3\leqslant d\leqslant 6$, it is possible under the additional restriction on $r$\textnormal{:} $d/2-1\leqslant r\leqslant 2$. If Conditions~\textnormal{\ref{Condition Lambda in L infty}} and \textnormal{\ref{Condition tilde Lambda in Lp}} hold simultaneously, it is also possible to remove $\Pi _\varepsilon$ from the corrector.
\end{remark}

\begin{proof}[Proof of Theorem~\textnormal{\ref{Theorem exp energy norm}}]
Denote 
\begin{equation*}
\Psi (t):=e^{-it B_\varepsilon}B_\varepsilon ^{-1}((B^0)^{-1}+\varepsilon K(\varepsilon ))e^{itB^0}.
\end{equation*}
Then, according to \eqref{K(eps)} and \eqref{R2},
\begin{equation*}
\frac{d\Psi (t)}{dt}
=ie^{-itB_\varepsilon}\mathcal{R}_2(\varepsilon)e^{itB^0}+ie^{-it B_\varepsilon }B_\varepsilon ^{-1}(\varepsilon \Lambda ^\varepsilon \Pi _\varepsilon b(\mathbf{D})+\varepsilon\widetilde{\Lambda}^\varepsilon \Pi _\varepsilon )e^{itB^0}.
\end{equation*}
So,
\begin{equation*}
\begin{split}
&\Psi (t)-\Psi (0)=e^{-itB_\varepsilon}B_\varepsilon ^{-1}((B^0)^{-1}+\varepsilon K(\varepsilon ))e^{itB^0}
-B_\varepsilon ^{-1}((B^0)^{-1}+\varepsilon K(\varepsilon))
\\
&= i\int _0^t e^{-isB_\varepsilon}\mathcal{R}_2(\varepsilon)e^{isB^0}\,ds
+i\int _0^t e^{-isB_\varepsilon }B_\varepsilon ^{-1}(\varepsilon \Lambda ^\varepsilon \Pi _\varepsilon b(\mathbf{D})+\varepsilon\widetilde{\Lambda}^\varepsilon \Pi _\varepsilon)e^{isB^0}\,ds.
\end{split}
\end{equation*}
Multiplying this identity by $e^{-itB^0}$ from the right, we obtain
\begin{equation*}
\begin{split}
&e^{-itB_\varepsilon}B_\varepsilon ^{-1}(B^0)^{-1}
-B_\varepsilon ^{-1}((B^0)^{-1}+\varepsilon K(\varepsilon))e^{-itB^0}
\\
&= -e^{-itB_\varepsilon}B_\varepsilon ^{-1}\varepsilon K(\varepsilon )
+i\int _0^t e^{-isB_\varepsilon}\mathcal{R}_2(\varepsilon)e^{i(s-t)B^0}\,ds
\\
&+i\int _0^t e^{-isB_\varepsilon }B_\varepsilon ^{-1}(\varepsilon \Lambda ^\varepsilon \Pi _\varepsilon b(\mathbf{D})+\varepsilon\widetilde{\Lambda}^\varepsilon \Pi _\varepsilon)e^{i(s-t)B^0}\,ds.
\end{split}
\end{equation*}
Together with \eqref{R2} and \eqref{K3(eps,t)}, this implies
\begin{equation}
\label{5.21 new}
\begin{split}
&(e^{-itB_\varepsilon}B_\varepsilon ^{-1}
-e^{-itB^0} (B^0)^{-1}-\varepsilon\mathcal{K}_3(\varepsilon ;t))(B^0)^{-1}
\\
&= \mathcal{R}_2(\varepsilon)(B^0)^{-1}e^{-itB^0} +B_\varepsilon ^{-1}\varepsilon K(\varepsilon)e^{-itB^0}
-e^{-itB_\varepsilon}B_\varepsilon ^{-1}\varepsilon K(\varepsilon )
\\
&+i\int _0^t e^{-isB_\varepsilon}\mathcal{R}_2(\varepsilon)e^{i(s-t)B^0}\,ds
\\
&+i\int _0^t e^{-isB_\varepsilon }B_\varepsilon ^{-1}(\varepsilon \Lambda ^\varepsilon \Pi _\varepsilon b(\mathbf{D})+\varepsilon\widetilde{\Lambda}^\varepsilon \Pi _\varepsilon)e^{i(s-t)B^0}\,ds.
\end{split}
\end{equation}

Combining \eqref{B0^2 L2 to H2}, \eqref{K(eps)<= L2 to L2}, \eqref{Beps1/2 appr corr},  \eqref{B eps -1/2 est}, \eqref{5.8 new}, \eqref{imaginary exp identity last summand estimate}, and \eqref{5.21 new}, we arrive at the estimate
\begin{align}
\label{im exp corr not final}
\begin{split}
&\Bigl\Vert B_\varepsilon ^{1/2}\left( e^{-itB_\varepsilon}B_\varepsilon ^{-1}-e^{-it B^0}(B^0)^{-1}-\varepsilon\mathcal{K}_3(\varepsilon ;t)\right)(B^0)^{-1}\Bigr\Vert _{L_2(\mathbb{R}^d)\rightarrow L_2(\mathbb{R}^d)}
\\
&\leqslant {C}_{35}\varepsilon (1+\vert t\vert),
\end{split}
\\
&{C}_{35}:=\max\lbrace 
c_5C_2\check{c}^{-1}+2\beta ^{-1/2}C_K ;c_5C_2+c_6\alpha _1^{1/2}M_1+\beta ^{-1/2}\widetilde{M}_1 \rbrace .
\nonumber
\end{align}
Bringing together \eqref{b_eps=> H1-norm}, \eqref{B0 on H2}, and \eqref{im exp corr not final}, we get
\begin{equation}
\label{exp corr H2 to H1}
\begin{split}
\Vert  e^{-itB_\varepsilon}B_\varepsilon ^{-1}-e^{-it B^0}(B^0)^{-1}-\varepsilon\mathcal{K}_3(\varepsilon ;t)\Vert _{H^2(\mathbb{R}^d)\rightarrow H^1(\mathbb{R}^d)}
\leqslant C_{36}\varepsilon (1+\vert t\vert),
\end{split}
\end{equation}
where $C_{36}:=c_6 {C}_{35}C_L$. 

By \eqref{b_eps=> H1-norm} and \eqref{B eps -1/2 est},
\begin{equation}
\label{exp gr 1}
\Vert e^{-itB_\varepsilon}B_\varepsilon ^{-1}\Vert _{L_2(\mathbb{R}^d)\rightarrow H^1(\mathbb{R}^d)}\leqslant c_6\Vert B_\varepsilon ^{-1/2}\Vert _{L_2(\mathbb{R}^d)\rightarrow L_2(\mathbb{R}^d)}\leqslant c_6\beta ^{-1/2}.
\end{equation}
Next, according to \eqref{B0^2 L2 to H2},
\begin{equation}
\Vert e^{-it B^0}(B^0)^{-1}\Vert _{L_2(\mathbb{R}^d)\rightarrow H^1(\mathbb{R}^d)}\leqslant \check{c}^{-1}.
\end{equation}
Finally, by analogy with \eqref{K1 est start} and \eqref{K1 est not final},
\begin{equation}
\label{K3 grubo}
\Vert \mathcal{K}_3(\varepsilon ;t)\Vert _{L_2(\mathbb{R}^d)\rightarrow H^1(\mathbb{R}^d)} \leqslant (2\varepsilon M_1+M_2)\alpha _1^{1/2}\check{c}^{-1}+(2\varepsilon \widetilde{M}_1+\widetilde{M}_2)\check{c}^{-1/2}.
\end{equation}
Bringing together \eqref{exp gr 1}--\eqref{K3 grubo}, we obtain
\begin{equation}
\label{exp corr L2 to H1}
\Vert  e^{-itB_\varepsilon}B_\varepsilon ^{-1}-e^{-it B^0}(B^0)^{-1}-\varepsilon\mathcal{K}_3(\varepsilon ;t)\Vert _{L_2(\mathbb{R}^d)\rightarrow H^1(\mathbb{R}^d)}
\leqslant C_{37},
\end{equation}
where $C_{37}:=c_6\beta ^{-1/2}+ \check{c}^{-1}+(2 M_1+M_2)\alpha _1^{1/2}\check{c}^{-1}+(2 \widetilde{M}_1+\widetilde{M}_2)\check{c}^{-1/2}$. Interpolating between \eqref{exp corr L2 to H1} and \eqref{exp corr H2 to H1}, we arrive at estimate \eqref{Th exp corr} with the constant $C_{34}:=C_{37}^{1-r/2}C_{36}^{r/2}$.
\end{proof}

\end{document}